\documentclass[final,onefignum,onetabnum,letterpaper]{siamart250211}

\usepackage{xr-hyper} 
\usepackage{xr}
\usepackage{graphicx}
\usepackage{hyperref}

\usepackage{float}
\usepackage{titlesec}

\usepackage{url}
\usepackage{xr-hyper}
\usepackage{amsmath,amsfonts, multirow, graphicx}
\usepackage[caption=false]{subfig} 
\usepackage{algorithm,algpseudocode}
 \usepackage{booktabs}
 \usepackage{changes}

\newtheorem{remark}{Remark}
\usepackage{xcolor}
 \definecolor{mygrn}{RGB}{0,140,70}

\usepackage[mathscr]{eucal}
\usepackage{amsbsy,bm}
\usepackage{amssymb}
\usepackage{hyperref,cleveref}
\usepackage{tikz}
\usepackage{calc}
\usetikzlibrary{calc,arrows}
\usetikzlibrary{decorations.pathreplacing}
\usepackage{changes}
\usepackage[T1]{fontenc}

\makeatletter
\newcommand*{\addFileDependency}[1]{
  \typeout{(#1)}
  \@addtofilelist{#1}
  \IfFileExists{#1}{}{\typeout{No file #1.}}
}
\makeatother

\usepackage[caption=false]{subfig} 
\usepackage[utf8]{inputenc}
\usepackage{xr-hyper}
\usepackage{url}
\usepackage{bm}
\usepackage{soul}
\usepackage{lipsum}
\usepackage{amsfonts}

\usepackage{epstopdf}
\ifpdf
  \DeclareGraphicsExtensions{.eps,.pdf,.png,.jpg}
\else
  \DeclareGraphicsExtensions{.eps}
\fi

\usepackage{datetime}
\newdateformat{monthyeardate}{%
  \monthname[\THEMONTH] \THEDAY, \THEYEAR}

\usepackage{academicons}
\usepackage{xcolor}
\renewcommand{\orcid}[1]{\href{https://orcid.org/#1}{\textcolor[HTML]{A6CE39}{orcid.org/#1}}}

\usepackage{amsmath}
\allowdisplaybreaks
\usepackage{amssymb}
\usepackage{hyperref,cleveref}
\usepackage{commath}
\usepackage{mathtools}
\usepackage{bbm}

\usepackage{color}
\usepackage{graphicx}
\usepackage[small]{caption}
\usepackage{subcaption}

\usepackage{relsize}
\usepackage{adjustbox}
\usepackage{algorithm,algpseudocode}
\usepackage{booktabs}
\usepackage{tikz}

\usepackage{verbatim}

\usepackage{enumitem}
\setlist[enumerate]{leftmargin=.5in}
\setlist[itemize]{leftmargin=.5in}

\newsiamthm{problem}{Problem}

\newsiamremark{hypothesis}{Hypothesis} 
\crefname{hypothesis}{Hypothesis}{Hypotheses}
\newsiamremark{example}{Example}
\newsiamthm{claim}{Claim}
\newsiamthm{conjecture}{Conjecture}

\headers{Reduced Kalman Filter with Optical Flow}{A.\ Keating, M.\ Pasha}

\title{A Scalable Sequential Framework for Dynamic Inverse Problems via Model Parameter Estimation
\thanks{
\monthyeardate\today 
}}
\author{
Aryeh Keating\thanks{Department of Mathematics \& Academy of Data Science, Virginia Tech, USA, (\email{aryehkeating@vt.edu},\orcid{0009-0002-8262-0133}); \email{mpasha@vt.edu}, \orcid{0000-0003-4249-2421})} 
\and 
Mirjeta Pasha\footnotemark[2] 
}
\usepackage{amsopn}
\usepackage{ifthen}
\usepackage{float}
\floatstyle{ruled}
\newfloat{algorithm}{htb}{alg}
\floatname{algorithm}{Algorithm}
\newcounter{algo@row}
\newcounter{algo@rowindent}
\newcommand{\algofont}[1]{\textbf{#1}}
\newcommand{\algonumbersize}[1]{\scriptsize{#1}}
\newcommand{\algopreitem}[1][\arabic{algo@row}]{\texttt{\algonumbersize{#1}}}
\newcommand{\algoitemskip}{\hspace{\value{algo@rowindent}cc}}
\newcommand{\algonewnestedopen}[2]{
	\newcommand{#1}[1][]{%
		\ifthenelse{\equal{##1}{}}{\item}{\item[{\algopreitem[##1]}]}
		\algoitemskip\algofont{#2}%
		\addtocounter{algo@rowindent}{1}%
		\ignorespaces
	}
}
\newcommand{\algonewnestedaux}[2]{
	\newcommand{#1}[1][]{
		\addtocounter{algo@rowindent}{-1}
		\ifthenelse{\equal{##1}{}}{\item}{\item[{\algopreitem[##1]}]}
		\algoitemskip\algofont{#2}%
		\addtocounter{algo@rowindent}{+1}%
		\ignorespaces
	}
}
\newcommand{\algonewnestedclose}[2]{
	\newcommand{#1}[1][]{
		\addtocounter{algo@rowindent}{-1}
		\ifthenelse{\equal{##1}{}}{\item}{\item[{\algopreitem[##1]}]}
		\algoitemskip\algofont{#2}%
		\ignorespaces
	}
}
\newcommand{\algonewcommand}[2]{
	\newcommand{#1}[1][default]{
		\ifthenelse{\equal{##1}{default}}{\item}{\item[{\algopreitem[##1]}]}%
		\algoitemskip\algofont{#2}%
		\ignorespaces
	}%
}
\newcommand{\algonewkeyword}[2]{\newcommand{#1}{\algofont{#2}}}
\algonewcommand{\STATE}{\ignorespaces}
\algonewcommand{\INPUT}{Input: }
\algonewcommand{\pINPUT}{\phantom{Input: }}
\algonewcommand{\COMPUTE}{Compute: }
\algonewcommand{\OUTPUT}{Output: }
\algonewcommand{\pOUTPUT}{\phantom{Output: }}

\algonewnestedopen{\IF}{if }
\algonewnestedaux{\ELSEIF}{else if }
\algonewnestedaux{\ELSE}{else }
\algonewnestedclose{\ENDIF}{end if }
\algonewnestedopen{\FOR}{for }
\algonewnestedclose{\ENDFOR}{end for }
\algonewnestedopen{\WHILE}{while }
\algonewnestedclose{\ENDWHILE}{end while }
\algonewcommand{\BREAK}{break}%

\algonewkeyword{\To}{to }%
\algonewkeyword{\Do}{do }%
\algonewkeyword{\Then}{then }%
\algonewkeyword{\End}{end }%
\algonewkeyword{\AND}{and }%
\algonewkeyword{\True}{true }%
\algonewkeyword{\False}{false }%
\algonewkeyword{\irbleigs}{irbleigs }%
\algonewkeyword{\tridiag}{tridiag}%
\algonewkeyword{\reorth}{reorth}%

\newcommand{\vx}{\mathbf{x}}
\newcommand{\vy}{\mathbf{y}}
\newcommand{\vH}{\mathbf{H}}
\newcommand{\vR}{\mathbf{R}}
\newcommand{\vQ}{\mathbf{Q}}
\newcommand{\vM}{\mathbf{M}}
\newcommand{\vC}{\mathbf{C}}
\newcommand{\vP}{\mathbf{P}}
\newcommand{\vS}{\mathbf{S}}
\newcommand{\vU}{\mathbf{U}}
\newcommand{\vI}{\mathbf{I}}
\newcommand{\vB}{\mathbf{B}}
\newcommand{\vA}{\mathbf{A}}

\newcommand{\vxi}{\boldsymbol{\xi}}

\DeclareMathOperator*{\argmin}{arg\,min}

\newcommand{\R}{\mathbb{R}}

\newcommand{\bB}{{\bf B}}

\newcommand{\bH}{{\bf H}}
\newcommand{\bI}{{\bf I}}

\newcommand{\bM}{{\bf M}}

\newcommand{\bP}{{\bf P}}
\newcommand{\bQ}{{\bf Q}}
\newcommand{\bR}{{\bf R}}

\newcommand{\bT}{{\bf T}}
\newcommand{\bU}{{\bf U}}
\newcommand{\bV}{{\bf V}}
\newcommand{\bW}{{\bf W}}
\newcommand{\bX}{{\bf X}}

\newcommand{\bg}{{\bf g}}

\newcommand{\br}{{\bf r}}
\newcommand{\bs}{{\bf s}}

\newcommand{\bw}{{\bf w}}
\newcommand{\bx}{{\bf x}}
\newcommand{\by}{{\bf y}}
\newcommand{\bz}{{\bf z}}

\usepackage{verbatim}

\newcommand{\valpha}{\boldsymbol{\alpha}} 
\newcommand{\vepsilon}{\boldsymbol{\epsilon}} 
\newcommand{\vmu}{\boldsymbol{\mu}} 



\usepackage{lineno}
\linenumbers

\ifpdf
\hypersetup{
  pdftitle={Priorconditioned projection methods},
  pdfauthor={Lindbloom}
}
\fi
\begin{document}
\nolinenumbers
\maketitle

\begin{abstract}
Large-scale dynamic inverse problems are often ill-posed due to model complexity and the high dimensionality of the unknown parameters. Regularization is commonly employed to mitigate ill-posedness by incorporating prior information and structural constraints. However, classical regularization formulations are frequently infeasible in this setting due to prohibitive memory requirements, necessitating sequential methods that process data and state information online, eliminating the need to form the full space–time problem. In this work, we propose a memory-efficient framework for reconstructing dynamic sequences of undersampled images from computerized tomography data that requires minimal hyperparameter tuning. The approach is based on a prior-informed, dimension-reduced Kalman filter with smoothing. While well suited for dynamic image reconstruction, practical deployment is challenging when the state transition model and covariance parameters must be initialized without prior knowledge and estimated in a single pass. To address these limitations, we integrate regularized motion models with expectation–maximization strategies for the estimation of state transition dynamics and error covariances within the Kalman filtering framework.
We demonstrate the effectiveness of the proposed method through numerical experiments on limited-angle and single-shot computerized tomography problems, highlighting improvements in reconstruction accuracy, memory efficiency, and computational cost.
\end{abstract}




\begin{DOI}
    Not yet assigned
\end{DOI}
\section{Introduction}\label{sec: intro}

Many applications rely on recovering an image or a sequence of temporal images from corrupted and noisy measurements, typically collected at the boundary or the surface of the objects of interest \cite{yang2020ground, pasha2023computational, lan2023spatiotemporal}. These data may be inhomogeneous, meaning that some regions differ significantly from others. Not only can these complex data be spatiotemporal, but the target solutions may also span both space and time. In such scenarios we seek to solve a sequence of linear inverse problems of the form
\begin{equation}\label{eq: IP}
\bH_i\bx_i + \epsilon_i = \by_i,
\end{equation}
where $\bH_i \in \R^{m_t \times n_s}$ is the known operator that models the physical properties of the measuring device, $\bx_i \in \R^{n_s}$ denotes the desired solution, $\epsilon_i$ is known as noise, and $\by_i$ represents the measured data across time instances $i=0, 1, \dots, T$. 
The goal of the inverse problem is to recover the sequence of unknown parameters $\bx_i$ from the available data $\by_i$ and full or partial knowledge of the forward operators $\bH_i$.

The inverse problem \eqref{eq: IP} is often \emph{ill-posed}, reflecting non-uniqueness, lack of existence, or instability of solutions under data perturbations. Regularization addresses these challenges by reformulating the problem as a nearby inverse problem that incorporates additional information 
\begin{equation}\label{eq: regIP}
    \argmin_{\bx \in \R^{n_s}} \left\{ 
        \mathcal{F}(\bx {; \by})
        + \mathcal{R}(\mathbf{x})
    \right\},
\end{equation}
-- a technique known as \emph{regularization}
with $\mathcal{F}$ and $\mathcal{R}$ denoting the data-fidelity and regularization term, respectively.
The choice of $\mathcal{F}$ is informed by the assumptions about the data-generating process and noise characteristics. 
For simplicity, we assume that $\epsilon$ is a realization of a standard normal noise process, i.e., $\epsilon \sim \mathcal{N}(\mathbf{0},\mathbf{I})$, yielding $\mathcal{F}(\bx_i {; \by_i}) = \|\bH_i\bx_i - \by_i\|_2^2$.\footnote{
If $\epsilon \sim \mathcal{N}(\mathbf{0},\boldsymbol{\bR})$ with a symmetric positive definite covariance matrix $\boldsymbol{\bR} \neq \mathbf{I}$, there exists a Cholesky decomposition $\boldsymbol{\bR} = \mathbf{C} \mathbf{C}^T$, and the problem can be whitened by substituting $\mathbf{H} \gets \mathbf{C}^{-1} \mathbf{H}$ and $\mathbf{y} \gets \mathbf{C}^{-1} \mathbf{y}$.}
The regularization term $\mathcal{R}$ encodes prior beliefs about the structure of the otherwise unknown parameter vector $\bx_i$, for example promoting spatial smoothness, temporal coherence, sparsity, or low-rank behavior \cite{buccini2020modulus, pasha2023computational, gazzola2019flexible, buccini2021linearized}.

An alternative to the above deterministic regularization setting is the \emph{Bayesian approach} to inverse problems \cite{stuart2010inverse,calvetti2023bayesian}
where we treat the unknown parameters as random variables and impose a prior distribution on them. 
For instance, assuming additive standard normal noise $\epsilon_i \sim \mathcal{N}(\mathbf{0},\mathbf{I})$ in \eqref{eq: IP} corresponds to a likelihood function $p( \mathbf{y}_i \vert \bx_i) \ \propto \  \exp \left (- \frac{1}{2}\| \bH_i \bx_i - \mathbf{y}_i \|^2_2 \right )$. 
Assuming further a prior distribution $p(\bx_i)$ for the random variable of interest $\bx_i$, \emph{Bayes'} rule prescribes a formula for the posterior density $p_{\rm pos}(\bx_i|\by_i) \ \propto \ p(\by_i \vert \bx_i) p(\bx_i)$. 
Furthermore, the role of the regularization term is now taken by a prior distribution $p(\bx_i)$, encoding our structural belief about $\bx_i$. 
This perspective is particularly relevant for \emph{dynamic inverse problems}, where the unknown evolves in time and the resulting inference problem is severely ill-posed. In imaging applications, such as \textit{computerized tomography} (CT) \cite{CT1,CT2, buccini2025krylov, lindbloom2025priorconditioned}, ill-posedness is aggravated by the high dimensionality of the space–time model, severe undersampling in both angular and temporal domains, and clinical constraints that limit data acquisition \cite{lindbloom2025priorconditioned, buccini2025krylov}. Consequently, classical reconstruction methods such as \textit{filtered back projection} (FBP) \cite{FBP1,FBP2}, which is known to perform well for static imaging with dense measurements, become unreliable in dynamic settings characterized by sparse data.

To address this limitation, a common strategy for solving dynamic inverse problems is the so-called \emph{all-at-once} (AAO) approach \cite{pasha2023computational, chung2018efficient}. In this framework, the entire space--time problem is formulated as a single, large-scale inverse problem, in which the state variables at all time steps are recovered simultaneously. Specifically, the forward operator $\bH$ is represented as a block-diagonal matrix,
$
\bH = \mathrm{BlockDiag}(\bH_0, \bH_1, \ldots, \bH_T),
$
where each block $\bH_i \in \mathbb{R}^{m_t \times n_s}$ corresponds to the forward operator at time step $i$. When the forward operator is identical across time steps, this structure simplifies to $\bH = \bI_{T+1} \otimes \bH_0$, where $\otimes$ denotes the Kronecker product. The vector $\by ={\rm{vec}}([\by_{0},\dots,\by_{T}])\in \R^{m}$ represents measured data that are contaminated by an unknown error (or noise) $\epsilon \in \R^{m}$ 
that may stem from measurement errors. Other recent developments include alternating approaches that recover object motion via \textit{optical flow} \cite{OF1,OF2,OF3} and simultaneously reconstruct all images in the sequence. These methods incorporate temporal information either by using the current optical flow estimate as a regularization term \cite{okunola2025efficient} or, within a Bayesian framework, by encoding temporal prior information through the covariance matrix \cite{lan2023spatiotemporal}. 

However, for large-scale problems, such formulations become impractical due to the extreme dimensionality involved, leading to prohibitive computational and memory requirements. An alternative strategy for dynamic imaging is to solve the inverse problem sequentially, for instance using \emph{Kalman filtering and smoothing}. In addition, approaches such as streaming and recycling/restarting \cite{pasha2023recycling} have been proposed to reduce computational and memory costs by reusing information from previous time steps. Despite these advances, many challenges remain, particularly in handling extremely high-dimensional states and long temporal sequences efficiently.
One effective strategy to mitigate these memory and computational bottlenecks is to solve the inverse problem sequentially. For example, \textit{Kalman filtering and smoothing} (KFS) \cite{anderson2005optimal,Kalman1960} is capable of performing image estimation from measurement data, acquired one time step at a time, avoiding the need to process the full space--time system simultaneously.

For the purposes of this paper, where the state dimension is extremely large due to image vectorization, we employ a reduced-order model within the \emph{Kalman filter and smoother} (RKFS) framework that restricts computations to low-dimensional operations \cite{hakkarainen2019undersampled,dim_red1,RKFS}. In this sequential approach, the state estimate is first \emph{predicted} forward in time using the system dynamics, and then \emph{updated} with new measurements to incorporate observation information. The reduced-order formulation leverages low-rank structure in both the state and covariance, thereby enabling these prediction and update steps to be carried out efficiently, significantly reducing the computational and memory demands compared to full-dimensional Kalman filtering.

A key advantage of the sequential approach is its reduced computational and memory demands compared to the AAO method. It also enables near real-time processing, as state estimates can be updated using current observations without restarting the reconstruction. This capability is particularly valuable in time-dependent imaging applications, such as CT, where patients may undergo long sequences of scans. Moreover, the RKFS framework can accommodate dynamic observations, including variations in patient position across separate scans. Related strategies addressing similar challenges have also been proposed in \cite{cite12,cite13,cite14}.

However, while the RKFS method is theoretically capable, accurately estimating the error covariance and update operator for the motion model remains a significant challenge, particularly when sequences exhibit non-stationary behavior. Although motion estimation through optical flow has been shown to produce high-quality reconstructions \cite{okunola2025efficient}, it is often computationally expensive. This can pose a significant restriction for users with limited computational resources or in applications that require real-time solutions, where fast and efficient methods are essential. Previous work addressing similar challenges can be found in \cite{emkal1,emkal2,emkal3}.

To address these issues, in this paper we propose three complementary techniques for estimating the motion model within the RKFS framework. Among these, we develop a method based on dynamic mode decomposition (DMD) to provide an efficient approximation of the motion model. We further employ an expectation–maximization (EM) procedure to refine parameter estimates and improve overall accuracy. Our goal is to develop an efficient RKFS method that is robust across a wide range of hyperparameter initializations while maintaining high-quality image reconstruction.

\paragraph{Overview of Main Contributions}
Our primary contributions are as follows: 

\begin{enumerate}
 \item[$\diamond$] We develop a scalable sequential framework for dynamic inverse problems that corrects model parameters using the \emph{Expectation--Maximization} (EM) algorithm \cite{shumway1982approach,EM2} and incorporates an automated, regularized motion model. The proposed approach leverages previous image estimates to construct improved approximations of object motion at each time step, while simultaneously solving a regularized inverse problem to mitigate noise and reconstruction artifacts. The combined use of EM and motion estimation enables successive refinement of both parameter matrices and image estimates.

\item[$\diamond$] We introduce several strategies for motion estimation within the proposed framework. In particular, we first consider an optical-flow-based motion model, which yields high-quality reconstructions but incurs substantial computational and memory costs. To address these limitations, we develop efficient reduced-order motion models based on \emph{dynamic mode decomposition} (DMD), including a regularized DMD variant. These approaches significantly reduce computational complexity and memory usage while maintaining reconstruction accuracy.

\item[$\diamond$] We retain dimensionality reduction throughout a modified variant of the \emph{reduced Kalman filter and smoother} (RKFS), enabling scalability to massive-scale problems arising from high-dimensional image vectorization and long temporal sequences.

\item[$\diamond$] We demonstrate the effectiveness of the proposed methods through a wide range of numerical experiments and imaging applications. Comparisons are provided against existing approaches, including all-at-once (AAO) methods—which rely on spatio-temporal regularization but suffer from prohibitive memory requirements—and standard RKFS approaches, which struggle to recover high-quality reconstructions in challenging dynamic settings.
\end{enumerate}

\paragraph{Outline} The paper is organized as follows. Section \ref{sec: backgroundddd} introduces the measurement model and reviews the dimension-reduced Kalman filter and smoother. Section \ref{sec: OF} describes techniques for estimating the motion model and Section \ref{sec: Proposed} describes the proposed method and its variants along with an analysis of the memory strain of the discussed methods. 
Section \ref{sec:experiments} presents numerical results and Section \ref{sec: Conclusions} discusses concluding remarks and outlook.

\section{Setup and Preliminaries}\label{sec: backgroundddd}

\subsection{Problem Setup -- Dynamic Inverse Problem}\label{sub: problem_formulation}
Consider a sequence of two-dimensional images $\bX(t_i) \in \R^{n_x \times n_y}$, where $n_x$ and $n_y$ denote the number of rows and columns (i.e., the number of pixels in the horizontal and vertical directions, respectively), for time indices $i = 0,1,\dots,T$ corresponding to times $t_0 < t_1 < \dots < t_T$. We use the operator $\text{vec}$ to vectorize an image $\bX(t)$ into a vector $\bx(t) = \text{vec}(\bX(t))\in \R^{n_s}$, where $n_s = n_x n_y$ represents the total number of pixels in the image at time $t$. For completeness, we define the vector $\bx = \text{vec}([\bx_0, \bx_1, \dots, \bx_{T}]) \in \R^{n}$, where $n = n_s\cdot (T+1)$ represents the total number of pixels across all time steps, and for simplicity set $\mathbf{x}_i=\mathbf{x}(t_i)$. We consider the forward model $\mathcal{H}(\vx)=\vy$ (also known as the data model), where $\mathcal{H}$ represents the \textit{Radon transform} and $\bx$ the image sequence that we seek to reconstruct given the observed data (also known as \textit{sinogram}) $\vy \in \R^{m}$. In the dynamic setting, for a given sequence of images $\{\vx_i\}_{i=0}^{T}$, we have access to only the observations $\{\vy_i\}_{i=0}^{T}$ and the corresponding Radon transform operators $\{\vH_i\}_{i=0}^{T} \in \R^{m_t \times n_s}$  (also known as design matrices). We then consider the noisy linearized forward model:
\begin{equation}\label{eq: dataModel}
    \vH_i\vx_i+\vepsilon_i=\vy_i,
\end{equation}
where $\{\vepsilon_i\}_{i=0}^{T} \in \R^{m_t}$ are our observation noise vectors which we assume to have zero-mean Gaussian distribution $\vepsilon_i \sim \mathcal{N}(\mathbf{0},\vR_i)$, where $\bR_i \in \R^{m_t \times m_t}$ denotes the covariance matrix. Our data model \eqref{eq: dataModel} in the dynamic setting is complemented by the stochastic model 
\begin{equation}\label{eq: stocModel}
    \vM_i\vx_{i-1}+\vxi_i=\vx_i
\end{equation}
where  $\bM_{i} \in \R^{n_s \times n_s}$ is the dynamic model used to encode temporal knowledge and $\{\vxi_i\}_{i=1}^{T}\in\mathbb{R}^{n_s}$ are process noise vectors associated with the state transition model which we assume to have Gaussian distribution $\vxi_i \sim \mathcal{N}(\mathbf{0},\vQ_i)$, where $\vQ_i$ is the covariance matrix. Our end goal of the inverse problem \eqref{eq: dataModel}-\eqref{eq: stocModel} is to recover a sequence of high quality images along with the model parameters that also involve parameters that define covariance matrices $\vR_i$ and $\vQ_i$.

\subsection{Kalman Filtering and Smoothing with Prior-Based Dimension Reduction}\label{sub: prior_dim_red}

Assuming Gaussian observation error for the noise vectors $\{\vepsilon_i\}_{i=1}^T$ with $\vepsilon_i \sim N(\mathbf{0},\vR_i)$ and assuming further that we have a Gaussian prior for $\{\vx_i\}_{i=1}^T$, i.e., $\vx_i\sim N(\vmu_i,\mathbf{\Sigma})$ with mean $\vmu_i\in\mathbb{R}^{n_s}$ and shared covariance matrix $\mathbf{\Sigma}\in\mathbb{R}^{n_s\times n_s}$, then the posterior density is given by 
\begin{equation}
p(\vx_i|\vy_i) \propto \exp \left(-\frac{1}{2} \left(\| \vy_i - \vH_i \vx_i \|_{\vR_i}^2 + \| \vx_i - \vmu_i \|_{\mathbf{\Sigma}}^2 \right) \right).
\end{equation}

Here we use the notation $\|\mathbf{z}\|_{\vR_i}^2 := \mathbf{z}^\top \vR_i^{-1} \mathbf{z}$. The vector $\vmu_i$ represents the expected value of the corresponding parameter $\vx_i$ given by the prior distribution and $\mathbf{\Sigma}$ is the error covariance matrix (where $\mathbf{\Sigma}^{-1}$ is also known as the \emph{precision matrix}).
Assuming further a linear forward operator $\mathbf{H}_i$, the resulting posterior distribution is Gaussian with mean $\vx_i^{\text{\text{est}}}$ and covariance $\vC_i^{\text{\text{est}}}$ given by \eqref{eq: mean} and \eqref{eq: cov}, respectively
\begin{align}\label{eq: mean}
\vx_i^{\text{\text{est}}} &= \vmu_i + \vC_i^{\text{\text{est}}} \vH_i^\top \vR_i^{-1}(\vy_i - \vH_i \vmu_i)  \\
\label{eq: cov}
\vC_i^{\text{\text{est}}} &= (\vH_i^\top \vR_i^{-1} \vH_i+ \mathbf{\Sigma}^{-1})^{-1}.
\end{align}
\paragraph{Remark} Both the deterministic and the Bayesian setting require inverting or processing very large-scale matrices. In this work we seek to reduce the size of the problem through a well known technique as reduced order modeling (ROM) which we describe below.

Given the prior mean $\vmu_i$ and a projection matrix $\mathbf{P}_r\in\mathbb{R}^{n_s\times r}$ that projects parameters $\boldsymbol{\alpha}_i\in\mathbb{R}^r$ from the low-dimensional subspace onto the full state space, then we have the following relation
\begin{equation}\label{eq: ROM}
    \vx_i = \vmu_i + \mathbf{P}_r \boldsymbol{\alpha}_i,
\end{equation}
where $r\leq n_s$ is the dimension of the reduced space for the parameters. For large state dimension we wish to choose $r\ll n_s$. 

Given the reduced-order model (ROM) relation in \eqref{eq: ROM}, we consider inference in a reduced parameter space. Specifically, we represent the unknown state in terms of a low-dimensional set of coefficients $\boldsymbol{\alpha}i$ defined in a reduced subspace. Under this representation, the posterior distribution of the reduced parameters $\boldsymbol{\alpha}i$ conditioned on the observations $\vy_i$ is given by
\begin{equation}\label{eq: posterior1}
p(\valpha_i \mid \vy_i) \propto
\exp \left(
-\frac{1}{2} \left(
\| \vy_i - \vH_i (\vmu_i + \vP\valpha_i) \|_{\vR_i}^2
+
\| \vP\valpha_i \|^2_{\boldsymbol{\Sigma}}
\right)
\right).
\end{equation}

The first term in the exponent corresponds to the data misfit, measuring the discrepancy between the observations and the predicted measurements obtained from the reduced-order representation. The second term acts as a regularization term  \cite{hansen1994regularization} induced by the prior, penalizing deviations from the prior mean in directions weighted by the prior covariance $\boldsymbol{\Sigma}$. Together, these terms define a Gaussian posterior distribution over the reduced parameters.

To construct the projection matrix $\vP$, we exploit the structure of the prior covariance matrix $\boldsymbol{\Sigma}$. In particular, we compute a low-rank approximation of $\boldsymbol{\Sigma}$ by retaining the eigenvectors (or left singular vectors) associated with its leading singular values. These dominant modes capture the directions of largest prior variance and therefore define an effective low-dimensional subspace in which the inference problem can be solved efficiently. The columns of $\vP$ are given by these leading modes, yielding a reduced representation that preserves the most informative components of the prior while significantly reducing computational complexity.

Specifically, using MATLAB notation and given reduction dimension $r$, we compute the SVD: $$\boldsymbol{\Sigma} = \vU\vS\vU^\top $$
and take $\vP_r:=\vU(:,1:r)\vS(1:r,1:r)^{1/2}$ as our projection matrix $\vP_r$ which yields the following
\begin{equation}\label{eq: norm_proj}
\| \vP_r \valpha_i \|_{\boldsymbol{\Sigma}}^2 = \valpha_i^\top \vP_r^\top \boldsymbol{\Sigma}^{-1} \vP_r \valpha_i = \valpha_i^\top \valpha_i=\|\valpha_i\|^2
\end{equation}
and the projection matrix has limited our search for $\valpha_i$ to consider the following posterior instead of \eqref{eq: posterior1}
\begin{equation}\label{eq: post_proj}
p(\valpha_i|\vy_i) \propto \exp \left(-\frac{1}{2} \left(\| \vy_i- \vH_i (\vmu_i + \vP_r\valpha_i) \|_{\vR_i}^2 + \| \valpha_i \|^2 \right) \right).
\end{equation}

The corresponding equations to \eqref{eq: mean} and \eqref{eq: cov} for the new mean and covariance estimation are then given by 
\begin{align}\label{eq: proj_mean}
\valpha_i^{\text{\text{est}}} &= \boldsymbol{\Psi}_i^{\text{\text{est}}} (\vH_i \vP_r)^\top \vR_i^{-1} (\vy_i - \vH_i \vmu_i), \\
\label{eq: proj_cov}
\boldsymbol{\Psi}_i^{\text{\text{est}}} &= ((\vH_i \vP_r)^\top \vR_i^{-1} (\vH_i \vP_r) + \vI)^{-1}.
\end{align}

\paragraph{Contributions From Linear Update Model}
Up to this point we have not taken into consideration the linear update model $\vM_i\vx_{i-1}+\vxi_i=\vx_i$, which is key to the method of Kalman filtering. We summarize this known as the Kalman filtering step: Given a time step  $i\in\{1,\dots,T\}$, we update the estimations of the mean and error covariance from the previous times step $i-1$, $\vx_{i-1}^{\text{\text{est}}},\vC_{i-1}^{\text{\text{est}}}$ to the current time step $i$ via:
\begin{align}\label{eq: mean_KF}
\vx_i^p &= \vM_i \vx_{i-1}^{\text{\text{est}}}, \\
\label{eq: cov_KF}
\vC_i^p &= \vM_i \vC_{i-1}^{\text{\text{est}}} \vM_i^\top + \vQ_i.
\end{align}
Using this prior in the update step, the posterior distribution $p(\vx_i|\vy_{1:i})$ is given by:

\begin{equation}\label{eq: pos_linear_update}
p(\vx_i|\vy_{1:i}) \propto \exp \left(-\frac{1}{2} \left(\| \vy_i - \vH_i \vx_i \|_{\vR_i}^2 + \| \vx_i - \vx_i^p \|_{\vC_i^p}^2 \right) \right).
\end{equation}
Hence, we introduce a new notation and denote the posterior distribution as $p(\vx_i|\vy_{1:i})$ rather than $p(\vx_i|\vy_{i})$, as in the filtering step, we have accessed all of our observations up to the given time step $i$, rather than the single observation $\vy_i$.

\subsection{Explicit Construction of ROM in Kalman Filtering}\label{sub: ROM}

We now incorporate the reduced order model within the Kalman filter method. 
Let $\vP_r$ be the fixed projection, then we parametrized our state via the reduced order model as: $\vx_i=\vx_i^p+\vP_r\valpha_i$. Thus using the fact that $\vC_i^{\text{\text{est}}} =\vP_r\boldsymbol{\Psi}_i^{\text{\text{est}}}\vP_r^\top$ along with equations (\ref{eq: proj_mean}) and (\ref{eq: proj_cov}), we set our mean and error covariance estimations as follows in the Kalman filter step as: \begin{align}
\vx_i^p &= \vM_i(\vx_{i-1}^p + \vP_r \valpha_{i-1}^{\text{\text{est}}})\label{eq: x_p} \\
\vC_i^p &= (\vM_i \vP_r) (\boldsymbol{\Psi}_{i-1}^{\text{\text{est}}}) (\vM_i \vP_r)^\top + \vQ_i.
\end{align}

The new posterior distribution for the parameters $\valpha_i$ can be estimated as \begin{equation}
p(\valpha_i|\vy_i) \propto \exp \left(-\frac{1}{2} \left(\| \vy_i - \vH_i \vx_i^p - \vH_i \vP_r \valpha_i \|_{\vR_i}^2 + \| \vP_r \valpha_i \|_{\vC_i^p}^2 \right)
\right).\label{eq: p(a|y)}
\end{equation}

This allows us to estimate our parameters $\valpha_i$ by taking the posterior mean with corresponding error covariance matrix from the posterior distribution (\ref{eq: x_p}) as follows:
\begin{align}
\valpha_i^{\text{\text{est}}} &= \boldsymbol{\Psi}_i^{\text{\text{est}}} (\vH_i \vP_r)^\top \vR_i^{-1} (\vy_i - \vH_i \vx_i^p), \\
\boldsymbol{\Psi}_i^{\text{\text{est}}} &= ((\vH_i \vP_r)^\top \vR_i^{-1} (\vH_i \vP_r) + \vP_r^\top (\vC_i^p)^{-1} \vP_r)^{-1}.
\end{align}

\subsection{Computational Efficiency for Kalman Filter}\label{sub: EKF}
In equation (\ref{eq: p(a|y)}) we must compute $(\vC_i^p)^{-1}\vP_r$. Directly computing the inverse is expensive or even prohibitive. To circumvent this, we introduce the following scheme. First, note that $$\vC_i^p= (\vM_i \vP_r) (\boldsymbol{\Psi}_{i-1}^{\text{\text{est}}}) (\vM_i \vP_r)^\top + \vQ_i = \vB_i\vB_i^\top+\vQ_i$$ where we have defined $\vB_i=\vM_i\vP_r\vA_i$ with $\vA_i$ being the square root matrix satisfying the factorization $\vA_i\vA_i^\top = \boldsymbol{\Psi}_{i-1}^{\text{\text{est}}}.$ We compute $(\vC_i^p)^{-1}\vP_r$ via the  Sherman-Morrison-Woodbury matrix inversion lemma as follows
\begin{equation}
(\vC_i^p)^{-1}\vP_r = \vQ_i^{-1}\vP_r - \vQ_i^{-1}\vB_i(\vB_i^\top\vQ_i^{-1}\vB_i + \vI)^{-1}\vB_i^\top\vQ_i^{-1}\vP_r.
\end{equation} which is possible given that direct computation of $\vQ_i^{-1}\vP_r$ is feasible.

For the purposes of this project, we initialize the covariance matrix $\boldsymbol{\Sigma}$ using the commonly employed \emph{squared exponential} (SE) covariance function, which is widely used in computed tomography and related imaging applications \cite{Covariance}. Specifically, the covariance entries are defined as
$
\boldsymbol{\Sigma}_{ij}
=
\alpha^2 \exp\!\left(-\frac{d(x_i,x_j)^2}{2\ell^2}\right),
$
where $\alpha^2$ denotes the marginal variance, $\ell$ is the correlation length, and $d(x_i,x_j)$ is the Euclidean distance between pixel locations $x_i$ and $x_j$. This choice promotes spatial smoothness in the reconstructed images through rapidly decaying correlations.

While the squared exponential covariance is well suited for CT applications, alternative covariance models may be more appropriate in other settings. For example, the Mat\'ern covariance class provides additional flexibility in controlling smoothness and is commonly used in many applications.

For a summary of the reduced Kalman filtering algorithm we refer the reader in Algorithm \ref{Alg:RKF} in Section~\ref{sub: App_A}.

\subsection{Reduced Order Rauch–Tung–Striebel Smoother}\label{sub: RORTS}
We now turn our attention to Rauch–Tung–Striebel smoother that is seen as a 
post processing tool for the outputs of the Kalman filter. In this approach we seek to estimate the probability density $p(\vx_i|\vy_{1:T})$ for all $i\in\{1,\dots, T\} $ given all the observations up to time $T$, i.e. $\vy_{1:T}$. To estimate this distribution we first set $\vx_{T}^{\text{sm}} = \vx_{T}^{\text{est}}$ which is given by our Kalman Filter and recursively obtain the smoothed states: $$\vx_{i-1}^{\text{sm}} = \vx_{i-1}^{\text{est}} + \mathbf{C}_{i-1}^{\text{est}} \mathbf{M}_i^{\top} (\mathbf{C}_i^{p})^{-1} (\vx_i^{\text{sm}} - \vx_i^{p})$$

We note that the smoothed covariance matrices $\vC_{i-1}^{\text{sm}}$ and the smoothed reduced-dimension covariance matrices $\boldsymbol{\Psi}_{i-1}^{\text{sm}}$ are not required for estimating $\vx_{i-1}^{\text{sm}}$ alone. However, since these quantities are needed to update the error covariance matrices ${\mathbf{Q}_i}$ and ${\mathbf{R}_i}$ in the subsequent sections, they must be computed and stored for this purpose. With this in mind, the reduced order calculations are provided with the further utilization of the Sherman-Morrison-Woodbury matrix inversion lemma for $(\mathbf{C}_i^{ {p}})^{-1}$ as follows: 

\[
\mathbf{C}_{ i-1}^{ {\text{sm}}} = \mathbf{P}_r \boldsymbol{\Psi}_{ i-1}^{ {\text{sm}}} \mathbf{P}_r^\top
\]

\begin{align}
\mathbf{P}_r \boldsymbol{\Psi}_{ i-1}^{ {\text{sm}}} \mathbf{P}_r^\top 
&= \mathbf{P}_r \boldsymbol{\Psi}_{ i-1}^{ {\text{est}}} \mathbf{P}_r^\top 
+ \mathbf{P}_r \boldsymbol{\Psi}_{ i-1}^{ {\text{est}}} \mathbf{P}_r^\top \mathbf{M}_i^\top (\mathbf{C}_i^{ {p}})^{-1} \nonumber \\
&\quad \times \left( \mathbf{P}_r \boldsymbol{\Psi}_i^{ {\text{sm}}} \mathbf{P}_r^\top - \mathbf{C}_i^{ {p}} \right) \nonumber \\
&\quad \times (\mathbf{C}_i^{ {p}})^{-1} \mathbf{M}_i \mathbf{P}_r \boldsymbol{\Psi}_{ i-1}^{ {\text{est}}} \mathbf{P}_r^\top, 
\end{align}

and furthermore
\begin{align}
\boldsymbol{\Psi}_{ i-1}^{{\text{sm}}} 
&= \boldsymbol{\Psi}_{ i-1}^{ {\text{est}}} 
+ \boldsymbol{\Psi}_{ i-1}^{ {\text{est}}} \mathbf{P}_r^\top \mathbf{M}_i^\top (\mathbf{C}_i^{ {p}})^{-1} \nonumber \\
&\quad \times \left( \mathbf{P}_r \boldsymbol{\Psi}_i^{ {\text{sm}}} \mathbf{P}_r^\top - \mathbf{C}_i^{ {p}} \right) \nonumber \\
&\quad \times (\mathbf{C}_i^{ {p}})^{-1} \mathbf{M}_i \mathbf{P}_r \boldsymbol{\Psi}_{ i-1}^{ {\text{est}}}. 
\end{align}\label{eq: reduced_cov}

All details of the smoothing procedure are provided in Algorithm~\ref{Alg:RKS} in Section~\ref{sub: App_B}.
\section{Estimating the Motion Model}\label{sec: OF}
In this section we describe three different ways to estimate the stochastic model $\bM_i$ which in this paper is also referred to as motion model. The first model -- so called \textit{motion 1} (M1) is a well-known optical flow model \cite{okunola2025efficient} with the disadvantage of high computational cost. To overcome the cost, we propose the second and the third models known as \textit{motion 2} (M2) and \textit{motion 3} (M3) that both rely on building a DMD model with the advantage of a low computational cost and competitive reconstruction quality. Details and discussion for all three motion models is provided in the following. 
\subsection{Motion Model 1} \label{sub: M1}
Now equipped with Kalman filtering and smoothing schemes, it is important to note the dependence both algorithms have on the motion update model with motion operators $\{\mathbf{M}_i\}_{i=1}^{T}$. To obtain such operators we seek to estimate a model that captures the feature evolution with time. For this, we utilize the optical flow model that rests on the assumption that the so called optical flow constraint (OFC) is satisfied. This constraint enforces that as the pixels in the image move or change position with time, their intensities remain the same. 

Given any single pixel position $(x^i(t),y^i(t)),\,i\in\{1,\dots,n_s\}$ at time $t$ of our true image $\vx(t)$ we assume:

\[
    \frac{d\vx(t,x^i(t),y^i(t))}{dt} = 0
    \]
With the OFC and given constant velocity $\bm s^i(t):= (s_x^i(t),s_y^i(t))$ over the interval $[t,t+\Delta t)$, as shown in \cite{okunola2025efficient}, we have the relation: 
\begin{equation} 
   \forall i\in\{1,\dots,n_s\}: s^i_x(t) \frac{\partial \mathbf{x}}{\partial x}(t,x^i(t), y^i(t)) + s^i_y(t) \frac{\partial \mathbf{x}}{\partial y}(t,x^i(t), y^i(t)) + \frac{\partial \mathbf{x}}{\partial t}(t,x^i(t), y^i(t)) =0
\end{equation}
 We write this relation as a linear system of equations: 
\begin{small}
\[
\resizebox{\linewidth}{!}{$
\underbrace{\begin{pmatrix}
\frac{\partial \mathbf{x}}{\partial x}(t,x^1(t), y^1(t)) & \cdots & 0 &
\frac{\partial \mathbf{x}}{\partial y}(t,x^1(t), y^1(t)) & \cdots & 0 \\
\vdots & \ddots & \vdots & \vdots & \ddots & \vdots \\
0 & \cdots & \frac{\partial \mathbf{x}}{\partial x}(t,x^{n_s}(t), y^{n_s}(t)) &
0 & \cdots & \frac{\partial \mathbf{x}}{\partial y}(t,x^{n_s}(t), y^{n_s}(t))
\end{pmatrix}}_{\mathbf{V}(\mathbf{x}(t))}
\underbrace{\begin{pmatrix}
s_x^1(t) \\ s_y^1(t) \\ \vdots \\ s_x^{n_s}(t) \\ s_y^{n_s}(t)
\end{pmatrix}}_{\mathbf{s}(t)}
= -
\underbrace{\begin{pmatrix}
\frac{\partial \mathbf{x}}{\partial t}(t,x^1(t), y^1(t)) \\[2pt]
\vdots \\[2pt]
\frac{\partial \mathbf{x}}{\partial t}(t,x^{n_s}(t), y^{n_s}(t))
\end{pmatrix}}_{\mathbf{T}(t)}
$}
\]
\end{small}
\normalsize
and thus we are able to solve for our velocity $\mathbf{s}(t)$ by solving the following inverse problem: 
\begin{equation}
\label{eq:2.7}
\min_{\mathbf{s} \in \mathbb{R}^{2n_s}} \left\|\mathbf{V}(\vx(t))\mathbf{s}(t) + \mathbf{T}(t) \right\|_p^p + \gamma \left\| \begin{pmatrix} \nabla \bs_x(t) \\ \nabla \bs_y(t) \end{pmatrix} \right\|_q^q,
\end{equation}
where
$$ \mathbf{s}_x(t) = \begin{pmatrix} s_x^1(t) & \dots & s_x^{n_s}(t) \end{pmatrix}^\top, \quad \mathbf{s}_y(t) = \begin{pmatrix} s_y^1(t) & \dots & s_y^{n_s}(t) \end{pmatrix}^\top, \quad \gamma \in \mathbb{R}.$$
Further we have $p,q\in\{1,2\}$ and $\nabla \bg$ denotes the gradient of $\bg$.
A choice of $p = 2$ can be interpreted as a constraint imposing smoothness on the velocity vector, while choosing $q = 1$ imposes sparsity/piece-wise constancy of the velocity vectors.
\paragraph{Construction of Linear Update Operators}\label{sub: LUO}

With the pixel velocity solve given in (\ref{eq:2.7}), and two images $\vx_{i-1}=\vx(t_{i-1}),\vx_i=\vx(t_i)$ with times $t_{i-1}<t_i$ and $i\in\{1,\dots,T\}$, we can now construct our linear update operator: $\vM_i$. To do so, as in \cite{okunola2025efficient}, we first consider the generalized (nonlinear) reverse motion function: $\mathcal{M}:\mathbb{R}^{n_s}\times \mathbb{R}^{2n_s}\rightarrow\mathbb{R}^{n_s}$ as follows: \begin{equation}
    \mathcal{M}(\vx_i,\mathbf{s}(t_{i-1})) := \vx_{i-1}
\end{equation}
We linearize this model via matrix $\vM(\mathbf{s}(t_{i-1}))\in\mathbb{R}^{n_s\times n_s}$ dependent on the array of pixel velocities $\mathbf{s}(t)$ such that: 
\begin{equation}
    \vM(\mathbf{s}(t_{i-1}))\vx_{i} = \vx_{i-1}
\end{equation}

Here, $\mathbf{M}(\mathbf{s}(t_{i-1}))$ denotes a bilinear warping operator constructed from the motion field $\mathbf{s}(t_{i-1})$; see \cite{okunola2025efficient} for details.
It is important to note that we have been only considering the reverse motion model. However, it is a simple calculation to obtain our desired forward operator: 
$ \vM(-\mathbf{s}(t_{i-1}))\vx_{i-1} = \vx_{i}$
and thus we set $\vM_{i} =  \vM(-\mathbf{s}(t_{i-1}))$ for $i\in\{1,\dots,T\}$ that is referred to as the \textit{motion model 1} (M1). 

One of the major computational burdens is on estimating the velocity parameters of the optical flow model, i.e., in estimating $\bs \in \R^{2n_s}$ in (\ref{eq:2.7}). To overcome this issue, we use the well-known and efficient majorization–minimization on a generalized Krylov subspace (MMGKS) method \cite{lanza2015generalized}, which we briefly describe in the following for completeness.
We consider the optimization problem: 
\begin{equation}\label{eq:rec_problem_s_fixed_joint}
    \bs^{(k+1)} = \min_{\mathbf{s} \in \mathbb{R}^{2n_s}} \left\|\mathbf{V}(\vx(t))\mathbf{s}(t) + \mathbf{T}(t) \right\|_2^2 + \lambda \|\boldsymbol{\Theta}\bs(t)\|_1,
\end{equation}
where 
\[\boldsymbol{\Theta} = \begin{pmatrix} \nabla s_x(t) \\ \nabla s_y(t) \end{pmatrix}.\]
Let $\bs^{(k)}(t)$ be an approximate solution of \eqref{eq:rec_problem_s_fixed_joint}, we define $\bz^{(k)}(t) = \boldsymbol{\Theta}^{(k)}\bs^{(k)}$ and the weighting matrix $\bP_{\epsilon}^{(k)} = \left(\text{diag}\left( \frac{1}{\phi_\epsilon(\bz^{(k)})}\right)\right)^{1/2},$

where operations on the right are performed element-wise. Then we obtain the quadratic tangent majorant $\mathcal{Q}(\bs, \bs^{(k)})$ for $\mathcal{J}_{\epsilon,\lambda}(\bs)$,  
with weighting matrix 
$\bP_{\epsilon}^{(k)}$ as 

\begin{equation}\label{eq: QuadraticMajorantQ}
\begin{array}{rcl}
\mathcal{Q}(\bs, \bs^{(k)})  &:=&
\|\mathbf{V}(\vx(t))\mathbf{s}(t) + \mathbf{T}(t)\|^{2}_{2}
+\lambda \left(\|\bP_{\epsilon}^{(k)}\boldsymbol{\Theta}^{(k)}\bs\|^{2}_{2}\right)+c,
\end{array}
\end{equation}
where $c$ denotes a suitable constant\footnote{A non-negative value of $c$ is technically required for the second condition in the definition to hold. However, since the minimizer is independent of the value of $c$, we will not discuss it further.} that is independent of $\bs^{(k)}$.
Minimizing the quadratic tangent majorant \eqref{eq: QuadraticMajorantQ} yields the next iterate $\bs^{(k+1)}$.
We proceed by setting the gradient of 
\eqref{eq: QuadraticMajorantQ}
to zero which gives the normal equations
\begin{equation}\label{eq: normaleqQuadMajorant}
(\left(\mathbf{V}(\vx(t))\right)^T\left(\mathbf{V}(\vx(t))\right) + \lambda  {\boldsymbol{\Theta}^{(k)}}^{T} (\bP_{\epsilon}^{(k)})^2\boldsymbol{\Theta}^{(k)})\bs(t) = \left(\mathbf{V}(\vx(t))\right)^T\bT(t).
\end{equation}
The condition \(\mathcal{N}(\left(\mathbf{V}(\vx(t))\right)^T\left(\mathbf{V}(\vx(t))\right)\cap \mathcal{N}({\boldsymbol{\Theta}^{(k)}}^{T} (\bP_{\epsilon}^{(k)})^2\boldsymbol{\Theta}^{(k)})=\{0\}\), which guarantees that the system \eqref{eq: normaleqQuadMajorant} has a unique solution, typically holds in practice. In that case, for fixed $\lambda$,
the solution 

of \eqref{eq: normaleqQuadMajorant} is the unique minimizer
of the quadratic tangent majorant $\mathcal{Q}(\bs, \bs^{(k)})$.

In practice, $\lambda$ is usually not known in advance but can be efficiently estimated at each iteration $k$ \cite{lanza2015generalized, pasha2023computational}.

Unfortunately, solving \eqref{eq: normaleqQuadMajorant} for large $\mathbf{V}(\vx(t))$ and  $\boldsymbol{\Theta}^{(k)}$ is typically be computationally expensive or even prohibitive.  Therefore, 
in \cite{huang2017majorization} the authors propose the 
MMGKS method that computes approximations by projection onto low dimensional subspaces. 
The method starts with a few steps of Golub-Kahan
bidiagonalization (GKB) 
with $\mathbf{V}(\vx(t))$ and $\bw_1 = \bT(t)/\| \bT(t) \|_2$ to determine the initial subspace $\bW_\ell$ such that 
\(\mathbf{V}(\vx(t))\bW_{\ell}=\bU_{\ell+1}\bB_{\ell}\),

and $\bW_\ell$ and
$\bU_{\ell+1}$ have orthonormal columns.

Given $\bs^{(k)}$ and 
$\bP_{\epsilon}^{(k)}$ with $k=0$ in the first step),

MMGKS computes the thin QR factorizations 

\begin{align}\label{eq: QR}
\mathbf{V}(\vx(t))\bW_{\ell+k} = \bQ_{\bV}\bR_{\bV}, \quad
\bP_{\epsilon}^{(k)}\Theta^{(k)}\bW_{\ell+k} = \bQ_{\Theta}\bR_{\Theta}.
\end{align}

which are now inexpensive to compute because of the smaller dimension. Restricting 
$\bs^{(k+1)}$ to 
${\rm range}(\bW_{\ell})$, 
$\bs^{(k+1)} = \bW_{\ell} \bz^{(k+1)}$, 
and adjusting \eqref{eq: normaleqQuadMajorant} accordingly leads to the following small system of equations for $\bz^{(k+1)}$, 

\begin{equation} \label{eq: minKryov2}
  (\bR_{\bV}^T\bR_{\bV} + 
  \lambda \bR_{\Theta}^T\bR_{\Theta})\bz^{(k+1)} =
  \bR_{\bV}^T\bQ_{\bV}^T\bT(t) .
\end{equation}

Its solution gives for
the residual vector of the (full) normal equations:
\begin{equation*}\label{eq: residual}
\br^{(k+1)}=\left(\mathbf{V}(\vx(t))\right)^T(\left(\mathbf{V}(\vx(t))\right)\bW_{\ell}\bz^{(k+1)} -\bT(t))+\lambda  {\boldsymbol{\Theta}^{(k)}}^T(\bP_{\epsilon}^{(k)})^2\boldsymbol{\Theta}^{(k)}\bW_{\ell}\bz^{(k+1)}.
\end{equation*}
We expand the solution subspace with a new normalized residual
$\bw_{\rm new}=\frac{\br^{(k+1)}}{\|\br^{(k+1)}\|_2}$,

\begin{equation}\label{eq: enlargeMM-GKS}
\bW_{\ell+1}=[\bW_{\ell},\bw_{\rm new}]\in\R^{n\times(\ell+1)}.\end{equation} 

In exact arithmetic,  $\bw_{\rm new}$ is orthogonal to the columns of $\bW_{\ell}$, but in computer arithmetic reorthogonalization of $\bW_{\ell+1}$ is typically needed. Next, $\bP_{\epsilon}^{(k)}$ is updated for the new solution estimate, and the process in steps \eqref{eq: QR} - \eqref{eq: enlargeMM-GKS} is repeated, expanding the solution space until a sufficiently accurate solution is reached. An overview of the linear update operator solver is given in Section~\ref{sub: App_C}.

\subsection{Motion Model 2} \label{sub: M2}

There are several inherent limitations associated with \textit{motion~1}. In constructing the motion matrices $\{\mathbf{M}_i\}_{i=1}^T$, we follow a two-step procedure:
\begin{itemize}
\item First, the pixel velocity fields ${\mathbf{s}(t_{i-1})}_{i=1}^T$ are estimated.
\item Second, the corresponding forward motion operators ${\mathbf{M}_i}_{i=1}^T$ are constructed.
\end{itemize}

This procedure leads to a loss of accuracy due to the accumulation of approximation errors, as well as increased computational cost resulting from repeated operations over time. To address these limitations, we propose an alternative approach based on a more efficient representation of the motion dynamics.

One straightforward way to construct such a motion operator is through Dynamic Mode Decomposition (DMD). In the general DMD setting, we consider a collection of linearly related input–output data pairs,

\begin{equation}
    \{(\mathbf{x}_i,\mathbf{y}_i)\}_{i=0}^T \text{ with } \mathbf{y}_i=\mathbf{M}_{\text{true}}\mathbf{x}_i
\end{equation}
where we define the DMD operator $\mathbf{M}_{\text{DMD}}$ as the least squares fit of the data: \begin{equation}\mathbf{M}^{\text{DMD}}\coloneqq\mathbf{Y}\mathbf{X}^+=\text{argmin}_{\mathbf{M}}\|\mathbf{M}\mathbf{X}-\mathbf{Y}\|_F^2
\end{equation}
where $\mathbf{X},\mathbf{Y}$ are the snapshot matrices formed from the provided input/output pairs: \begin{equation}\mathbf{X} \coloneqq\begin{bmatrix}
\vert & \vert &        & \vert \\
\mathbf{x}_0 & \mathbf{x}_1 & \cdots & \mathbf{x}_{T} \\
\vert & \vert &        & \vert
\end{bmatrix}
,
\quad
\mathbf{Y} \coloneqq 
\begin{bmatrix}
\vert & \vert &        & \vert \\
\mathbf{y}_0 & \mathbf{y}_1 & \cdots & \mathbf{y}_T \\
\vert & \vert &        & \vert
\end{bmatrix}\end{equation}

Now in the filtering setting, we only have access to the single data pair $(\mathbf{x}_{i-1}, \mathbf{x}_i)$ to form the motion matrix $\mathbf{M}_i$ in which case the DMD setting simplifies: $\mathbf{X}= \mathbf{x}_{i-1}$,  $\mathbf{Y}= \mathbf{x}_{i}$, and the DMD operator becomes the \textit{rank-1} matrix 
\begin{equation}\mathbf{M}_i^{\text{DMD}}\coloneqq\mathbf{x}_{i}\mathbf{x}_{i-1}^+ = \frac{\mathbf{x}_{i}\mathbf{x}_{i-1}^T}{\left\| \mathbf{x}_{i-1}\right\|_2^2}
\end{equation}

Now this formula gives an exact fit ($\mathbf{M}^{\text{DMD}}_i\mathbf{x}_{i-1}=\mathbf{x}_{i}$) given that $\mathbf{x}_{i-1}\neq \mathbf{0}$. However, this isn't necessarily what is desired. In particular, given any errors in the input $\mathbf{x}_{i-1}$ or the output $\mathbf{x}_{i}$, the operator will match these inconsistencies exactly, causing a divergence from the desired true solution. For this reason we must enforce some form of regularization to prevent this overfitting. 

The most straightforward method to regularize this approach is via \textit{Tikhonov regularization} in which we look to solve the modified expression: 

\begin{equation}\mathbf{M}_i^{\text{DMD},\zeta}\coloneqq\frac{\mathbf{x}_{i}\mathbf{x}_{i-1}^\top}{\left\| \mathbf{x}_{i-1}\right\|_2^2+\zeta}=\text{argmin}_{\mathbf{M}}\|\mathbf{M}\mathbf{x}_{i-1}-\mathbf{x}_{i}\|_F^2+\zeta\|\mathbf{M}\|_{F}^2
\end{equation} where $\zeta>0$ is our regularization parameter. This final regularized forward operator is referred to as \textit{motion model 2} (M2).

\subsection{Motion Model 3} \label{sub: M3}

While the regularization term, $\zeta$, in the second motion model, adds stability in the rank-1 least squares solves, there is an inherent lack of flexibility due to the regularization being applied to the entire forward operator. This global application can limit the effect the regularization has on pixel values closer to noise level (which we wish to reduce to 0). To counter this limitation, we isolate $k$ local fixed areas across both images $\mathbf{x}_{i-1},\mathbf{x}_{i}$, as non-overlapping and independent patches of fixed dimension $z_x\times z_y$. Here $z_x$ and $z_y$ divide $n_x$ and $n_y$, respectively, and we apply the regularized motion model M2 along each of these patches independently to form \textit{motion model 3} (M3):  
\begin{equation}\mathbf{M}_i^{\text{PDMD},\zeta}\coloneqq\sum_{j=1}^k \frac{\mathbf{S}_j{\mathbf{x}}^j_{i}(\mathbf{S}_j{\mathbf{x}}^j_{i-1})^\top}{\left\|\mathbf{S}_j{\mathbf{x}}^j_{i-1}\right\|_2^2+\zeta}
\end{equation}
where $\mathbf{S}_j\in\{0,1\}^{n_xn_y\times z_xz_y}$ represents the projection of the pixel values of the $j$th patch of both images, $\mathbf{x}^j_{i},\mathbf{x}^j_{i-1}\in\mathbb{R}^{ z_xz_y}$ (the superscript notation denotes the $j$th patch of an arbitrary image), to the full space. Hence we have the property: $\forall   
  \mathbf{z}\in\mathbb{R}^{n_xn_y}:\sum_{j=1}^k {\mathbf{S}_j{\mathbf{z}^j}}=\mathbf{z}$. This construction of the non-overlapping patchwise regularized DMD operators is referred to as \textit{motion model 3} (M3) and is our final motion model we will consider. 

\begin{remark}
While we focus on three motion models in this work, a number of promising alternatives merit further investigation. In particular, temporally enriched variants of DMD that exploit multiple consecutive state pairs may significantly improve robustness to noise and model mismatch. Additionally, estimating motion operators in a learned or reduced subspace could substantially lower computational and memory costs while preserving dominant dynamics. Such approaches may offer improved scalability and stability for large-scale dynamic inverse problems and represent compelling directions for future work.
\end{remark}

\section{The proposed method}\label{sec: Proposed}
\subsection{Expectation Maximization in Time Varying Reduced Kalman Filtering and Smoothing}
Building on the previous background and discussions, we now introduce the proposed method and describe the incorporation of expectation maximization into time-varying reduced Kalman filtering and smoothing.
We consider the joint log-likelihood of the complete data  
$\mathbf{x}_0, \mathbf{x}_1, \ldots, \vx_T, \mathbf{y}_1, \ldots, \mathbf{y}_T$.  
For the time-varying case, where the noise covariances can change at each step, it can be written in a form similar to that of \cite{shumway1982approach}:

\begin{align}
\log L &\doteq -\frac{1}{2} \log |\mathbf{\Sigma}|
- \frac{1}{2} (\mathbf{x}_0 - \vmu_0)^\top \mathbf{\Sigma}^{-1} (\mathbf{x}_0 - \vmu_0) \notag \\
&\quad - \frac{1}{2} \sum_{i=1}^{T} \log |\mathbf{Q}_i|
- \frac{1}{2} \sum_{i=1}^{T} (\mathbf{x}_i - \mathbf{M}_i \mathbf{x}_{i-1})^\top \mathbf{Q}_i^{-1} (\mathbf{x}_i - \mathbf{M}_i \mathbf{x}_{i-1}) \notag \\
&\quad - \frac{1}{2} \sum_{i=1}^{T} \log |\mathbf{R}_i|
- \frac{1}{2} \sum_{i=1}^{T} (\mathbf{y}_i - \mathbf{H}_i \mathbf{x}_i)^\top \mathbf{R}_i^{-1} (\mathbf{y}_i - \mathbf{H}_i \mathbf{x}_i).
\end{align}

Here, $\log L$ is to be maximized with respect to the parameter sequences $\{\mathbf{Q}_i\}_{i=1}^T$ and $\{\mathbf{R}_i\}_{i=1}^T$. Since the log likelihood is dependent on unobserved data series $\mathbf{x}_i$, $i = 0, 1, \ldots, T$, we apply EM algorithm conditionally with respect to the observed series $\mathbf{y}_1, \mathbf{y}_2, \ldots, \mathbf{y}_T$. Specifically, we update the estimated parameters at the $(j+1)$-st iteration such that $\{\mathbf{Q}_i^{(j+1)}\}$ and $\{\mathbf{R}_i^{(j+1)}\}$ maximize

\begin{equation}
G(\{\mathbf{Q}_i\}, \{\mathbf{R}_i\}) = \mathbb{E}^{(j)}[\log L \mid \mathbf{y}_1, \ldots, \mathbf{y}_T]\label{eq: G(Q,R)}
\end{equation}

where $\mathbb{E}^{(j)}$ denotes the conditional expectation relative to a density containing the $j$th iterate values $ \{\mathbf{Q}_i^{(j)}\}$ and $\{\mathbf{R}_i^{(j)}\}$. Any $j$th iterate element will be denoted via a superscript $(j)$ if such distinction is necessary. In order to calculate the conditional expectation defined in (\ref{eq: G(Q,R)}), it is convenient to define the conditional mean

With the smoothed estimates:
\begin{equation}
\mathbf{x}_i^{\text{sm}}= \mathbb{E}[\mathbf{x}_i \mid \mathbf{y}_1, \ldots, \mathbf{y}_T]
\end{equation}
and smoothed covariances
\begin{equation}
\mathbf{C}_i^{\text{sm}} = \text{cov}(\mathbf{x}_i \mid \mathbf{y}_1, \ldots, \mathbf{y}_T),\quad \mathbf{C}_{i,i-1}^{\text{sm}} = \text{cov}(\mathbf{x}_i, \mathbf{x}_{i-1} \mid \mathbf{y}_1, \ldots, \mathbf{y}_T)
\end{equation}

we express the function $G$ explicitly as in \cite{shumway1982approach}:

\small{
\begin{equation}
\begin{split}
G(\{\mathbf{Q}_i\}, \{\mathbf{R}_i\}) 
&= - \frac{1}{2} \log |\mathbf{\Sigma}| - \frac{1}{2} \, \text{tr} \left\{ \mathbf{\Sigma}^{-1} \left(\mathbf{C}_0^{\text{sm}} + (\mathbf{x}_0^{\text{sm}} - \vmu_0)(\mathbf{x}_0^{\text{sm}} - \vmu_0)^\top\right) \right\} \\
&\quad - \frac{1}{2} \sum_{i=1}^T \log |\mathbf{R}_i| - \frac{1}{2} \sum_{i=1}^T\text{tr} \left\{ \mathbf{R}_i^{-1}\left[ (\mathbf{y}_i - \mathbf{H}_i \mathbf{x}_i^{\text{sm}})(\mathbf{y}_i - \mathbf{H}_i \mathbf{x}_i^{\text{sm}})^\top + \mathbf{H}_i \mathbf{C}_i^{\text{sm}}\mathbf{H}_i^\top \right] \right\} \\
&\quad - \frac{1}{2} \sum_{i=1}^T \log |\mathbf{Q}_i| \\
&\quad -\frac{1}{2} \sum_{i=1}^T\text{tr} \left\{\mathbf{Q}_i^{-1} \left[ \left(\mathbf{C}_i^{\text{sm}} + \mathbf{x}_i^{\text{sm}} (\mathbf{x}_i^{\text{sm}})^\top\right) - \left(\mathbf{C}_{i,i-1}^{\text{sm}} + \mathbf{x}_i^{\text{sm}} (\mathbf{x}_{i-1}^{\text{sm}})^\top \right) \mathbf{M}_i^\top \right. \right. \\
&\left. \left. \quad - \mathbf{M}_i\left(\mathbf{C}_{i,i-1}^{\text{sm}} + \mathbf{x}_i^{\text{sm}} (\mathbf{x}_{i-1}^{\text{sm}})^\top \right)^\top + \mathbf{M}_i \left( \mathbf{C}_{i-1}^{\text{sm}} + \mathbf{x}_{i-1}^{\text{sm}} (\mathbf{x}_{i-1}^{\text{sm}})^\top \right) \mathbf{M}_i^\top \right] \right\}
\end{split}
\end{equation}
}
\normalsize
where $\text{tr(\bM)}$ denotes the trace of $\bM$. It follows that the choices

\begin{small}
\begin{align}
\mathbf{R}_i^{(j+1)} 
&= (\mathbf{y}_i - \mathbf{H}_i \mathbf{x}_i^{\text{sm}(j)})
   (\mathbf{y}_i - \mathbf{H}_i \mathbf{x}_i^{\text{sm}(j)})^\top
   + \mathbf{H}_i \mathbf{C}_i^{\text{sm}(j)} \mathbf{H}_i^\top
   \label{eq:R^i+1} \\[0.5em]
\mathbf{Q}_i^{(j+1)} 
&= \big(
      \mathbf{C}_i^{\text{sm}(j)} 
      + \mathbf{x}_i^{\text{sm}(j)}(\mathbf{x}_i^{\text{sm}(j)})^\top
   \big)
   - \big(
      \mathbf{C}_{i,i-1}^{\text{sm}(j)} 
      + \mathbf{x}_i^{\text{sm}(j)}(\mathbf{x}_{i-1}^{\text{sm}(j)})^\top
   \big)\mathbf{M}_i^{(j)\top} \nonumber \\[3pt]
&\quad
   - \mathbf{M}_i^{(j)}
     \big(
        \mathbf{C}_{i,i-1}^{\text{sm}(j)}
        + \mathbf{x}_i^{\text{sm}(j)}(\mathbf{x}_{i-1}^{\text{sm}(j)})^\top
     \big)^\top
   + \mathbf{M}_i^{(j)}
     \big(
        \mathbf{C}_{i-1}^{\text{sm}(j)}
        + \mathbf{x}_{i-1}^{\text{sm}(j)}(\mathbf{x}_{i-1}^{\text{sm}(j)})^\top
     \big)\mathbf{M}_i^{(j)\top} 
\label{eq:Q^i+1}
\end{align}
\end{small}

\normalsize
maximize $G(\{\mathbf{Q}_i\},\{\mathbf{R}_i\})$ during the M-step.

All smoothed terms on the right-hand side are computed during the E-step using the $j$th iterate parameters. We add superscript to linear update operators $\mathbf{M}_i$ to allow instances where they have been updated via optical flow.

\bigskip
In order to obtain the terms $\mathbf{C}_i^{\text{sm}} $, and $\mathbf{C}_{i,i-1}^{\text{sm}} $, we utilize Reduced Kalman Smoother via equations:
\begin{align}
\mathbf{C}_i^{\text{sm}} &= \mathbf{P}_r \boldsymbol{\Psi}_i^{\text{sm}} \mathbf{P}_r^\top \label{eq:C_i}\\
\mathbf{C}_{i,i-1}^{\text{sm}} &= (\mathbf{P}_r \boldsymbol{\Psi}_i^{\text{sm}} \mathbf{P}_r^\top) (\mathbf{C}_i^p)^{-1} \mathbf{M}_i (\mathbf{P}_r \boldsymbol{\Psi}_{i-1}^{\text{est}} \mathbf{P}_r^\top)\label{eq:C_i,i-1}
\end{align}
\subsection{Iterative Reduced Kalman Filter with Expectation Maximization and Motion Modeling (EMIRKFS-M)}

Equipped with the motion methods M1,M2, and M3, and Expectation Maximization algorithm, we propose an approach to utilize both methodologies for improved parameter estimation. A restriction is that image estimates at all time steps are required for implementation and thus, we structure our algorithm via an iteratively updating approach. Given initialized parameters at the $(j-1)$th iteration, we run a pass of the Reduced Kalman Filter and Smoother via (\ref{sub: App_A}) and (\ref{sub: App_B}) and pass image estimates $\{\mathbf{x}_i^{\text{sm}(j-1)}\}_{i=0}^T$ to obtain the updated parameters $\{\mathbf{R}_i^{(j)}\}_{i=1}^T,\{\mathbf{Q}_i^{(j)}\}_{i=1}^T$, and $\{\mathbf{M}_i^{(j)}\}_{i=1}^T$, which we use as parameters for the next iteration. 

We summarize the procedure in the Algorithm \ref{Alg: EMIRKFS-M}. 
\begin{algorithm}[H]
\caption{\textbf{EMIRKFS-M: Iterative Reduced Kalman Filter with Expectation Maximization and Motion Modeling}}
\label{Alg: EMIRKFS-M}
\normalsize
\begin{algorithmic}[1]
\Require Covariance parameters $\alpha$, $\ell$; reduction dimension $r$; iterations $n_{\text{iter}}$; initial estimates $\mathbf{x}_0^{{\text{\text{est}}}}$, $\bm{\Psi}_0^{{\text{\text{est}}}}$; motion operators $\{\mathbf{M}^{(0)}_i\}_{i=1}^{T}$; noise covariances $\{\mathbf{Q}^{(0)}_i\}_{i=1}^{T}$, $\{\mathbf{R}^{(0)}_i\}_{i=1}^{T}$; observations $\{\mathbf{y}_i\}_{i=1}^{T}$ 
\Statex

\For{$j = 1$ to $n_{\text{iter}}$}

    \State Utilizing $\{\mathbf{M}_i^{(j-1)}\}_{i=1}^{T},\{\mathbf{Q}_i^{(j-1)}\}_{i=1}^T, \{\mathbf{R}_i^{(j-1)}\}_{i=1}^T$
    \State \quad Compute $\{\mathbf{x}_i^{\text{\text{est}}(j-1)}\}_{i=0}^{T}$ using (\ref{Alg:RKF})
    \State \quad Compute $\{\mathbf{x}_i^{\text{\text{sm}}(j-1)}\}_{i=0}^{T}$ using (\ref{Alg:RKS}) and compute equations (\ref{eq:C_i}) and (\ref{eq:C_i,i-1})
    \For{$z = 1$ to $T$}
        \State Update $\mathbf{M}_z^{(j)}$ from $\mathbf{x}_{z-1}^{\text{\text{sm}}(j-1)}$ and $\mathbf{x}_z^{\text{\text{sm}}(j-1)}$
    \EndFor
    \State Update $\{\mathbf{Q}_i^{(j)}\}_{i=1}^T, \{\mathbf{R}_i^{(j)}\}_{i=1}^T$ using (\ref{eq:R^i+1}) and (\ref{eq:Q^i+1})
\EndFor
\end{algorithmic}
\end{algorithm}

\subsection{Computation and Memory Issues}\label{sec: CMI}
Before discussing performance, we first highlight the key computational and memory challenges associated with all-at-once formulations. We note that the biggest memory strain of the all-at-once approach is the storing of matrix $\mathbf{H} = \mathrm{BlockDiag}(\mathbf{H}_0,\dots,\mathbf{H}_T)$. Assuming our images are of dimension $n_x\times n_y$ with vectorized form $\mathbf{x}_i\in\mathbb{R}^{n_s}$, and our radon operators are discretized via a total of $s$ views, we have  $\mathbf{H}_i\in\mathbb{R}^{l\times n_s}$ where $l= O(s\sqrt{n_x^2+n_y^2})=O(s\sqrt{n_s})$. 
After utilization of sparsity, we have to store approximately $ \sqrt{n_s}$ total nonzero elements per each of the $l$ total rays, requiring in totality $O(sTn_s)$ total nonzero elements of $\vH$.

Regarding the EMIRKFS-M algorithm, for usability, we restrict that matrices $\{\mathbf{Q}_i\}_{i=1}^T$ are both easy to invert and aren't stored densely. It is important to note that under the assumption that pixel errors are uncorrelated from our motion model, we have that our error covariance matrices $\{\mathbf{Q}_i\}_{i=1}^T$ will have zero off-diagonal entries, and thus storage of the matrices' main diagonal is both necessary and sufficient for our purposes. This requires storage of $O(Tn_s)$ elements. Also, we must form several dense $n_s\times r$ matrices throughout the algorithm,  which require the storage of an additional $O(n_sr)$ elements. This gives a total storage requirement of $O(n_s(r+T))$ elements.

\section{Numerical Experiments}
\label{sec:experiments}
Despite these storage requirements, we now demonstrate the capability of the EMIRKFS-M algorithm via three numerical experiments. The first of which utilizes a simulated moving blocks dataset with limited motion over a modestly sized image for testing of the algorithm's accuracy and memory strain. The second dataset features a moving MNIST dataset, also simulated, to test the algorithm's effectiveness with significant non-stationary behavior. The third considers real data of an “emoji” phantom obtained at the University of Helsinki \cite{emoji,pasha2024trips}. We will also assess the asymptotic characteristics when many iterations are used. In our experiments, we determine the quality of various reconstructions using the relative residual error (RRE) between the an estimate $\tilde{\vx}$ and the ground truth $\vx$ defined via the formula $\text{RRE}(\tilde{\vx}, \vx) = \frac{\|\tilde{\vx} - {\vx}\|_2}{\|{\vx}\|_2}$\label{RRE}.   
We also define the noise level as $\sigma_{\text{NL}} = \|\mathbf{y}-\mathbf{\bH} \mathbf{x}_{\text{true}} \|_2/ \| \mathbf{\bH} \mathbf{x}_{\text{true}} \|_2$.

The code associated with this paper will be made publicly available in the following GitHub repository upon acceptance of the manuscript: \url{https://github.com/aryk4747/Sequential_IP_with_RKFS} \label{github}.

\subsection{Overview of Algorithms Tested}
\label{algorithms}
As a benchmark, we compare the image estimates of the algorithm with the following methods: 

\begin{enumerate}

\item An iterative reduced Kalman filtering and smoothing (IRKFS) method, which omits the parameter update from both the motions models and expectation maximization methods (and hence performs no update at all and gives identical reconstructions to the original RKFS method \cite{hakkarainen2019undersampled}).

\item An iterative reduced Kalman filtering and smoothing method with motion modeling (IRKFS-M1, IRKFS-M2, and IRKFS-M3, dependent on which of motion estimation is used from Subsections \ref{sub: M1},  \ref{sub: M2}, and \ref{sub: M3}, respectively) which omits the parameter update from only the expectation maximization method (line 8 in Algorithm \ref{Alg: EMIRKFS-M}).
\item An iterative reduced Kalman filtering and smoothing method with expectation maximization (EMIRKFS) which omits the parameter update from only the motion models (lines 5-7 in Algorithm \ref{Alg: EMIRKFS-M}).

\item The full iterative reduced Kalman filtering and smoothing method with expectation maximization and motion modeling method (\ref{Alg: EMIRKFS-M}) (EMIRKFS-M1, EMIRKFS-M2, and EMIRKFS-M3, dependent on which of motion estimation is used from Subsections \ref{sub: M1}, \ref{sub: M2}, and \ref{sub: M3}, respectively).

\item Three variants of all-at-once  implementation of the Majorization Minimization on a Generalized Krylov Subspace (MMGKS) method \cite{okunola2025efficient}: 
    \begin{enumerate}
        \item \textbf{AAO:} The all-at-once (AAO) formulation simultaneously solves the regularized minimization problem
$
\|\bH \bx - \by\|_2^2 + \lambda \|\bx\|_1.$
In the absence of temporal regularization, this problem decouples across time and is equivalent to solving $T+1$ independent minimization problems of the form
$
\|\bH_i \bx_i - \by_i\|_2^2 + \lambda \|\bx_i\|_1,
$ each corresponding to a static inverse problem at time step $i$ \cite{pasha2023computational}.
\item \textbf{AAO-ST:} The all-at-once formulation with spatial--temporal regularization (AAO-ST) solves the dynamic inverse problem
$
\|\bH \bx - \by\|_2^2 + \lambda \|\Psi \bx\|_1,
$
where $\Psi$ is a regularization operator that enforces both spatial and temporal structure. In particular, this formulation corresponds to the anisotropic total variation (AnisoTV) model introduced in \cite{pasha2023computational}.
\item \textbf{AAO-OF:} The all-at-once formulation with optical-flow-based regularization (AAO-OF) solves the dynamic inverse problem
$
\|\bH \bx - \by\|_2^2 + \lambda \,\mathcal{R}_{\mathrm{OF}}(\bx),
$ where $\mathcal{R}_{\mathrm{OF}}$ incorporates both spatial and temporal regularization through an optical flow model that enforces temporal coherence between consecutive time steps. In particular, we refer to the method proposed in \cite{okunola2025efficient} as AAO-OF.
\end{enumerate}
\end{enumerate}

 To ensure fairness, all variations of AAO ((a)-(c)) were run for 100 iterations and utilized generalized cross-validation (GCV) for selection of the parameter $\lambda$.

Experiments were performed using Virginia Tech's TinkerCliffs machine \cite{ARCVT}, a general-purpose CPU and GPU cluster. This cluster has 308 nodes with 39,424 AMD Rome CPU cores, 16 nodes with 1,536 Intel Xeon CPUs, 8 high-memory nodes with 1 TB RAM each, 14 nodes with 8 NVIDIA A100-80GB GPUs each, and 7 nodes with 8 NVIDIA H200-141GB GPUs each. Nodes are connected via HDR InfiniBand offering 100 Gbps throughput.
\subsection{Experiment 1: Simulated Dynamic Blocks Dataset}\label{ssec:exp_phantom}
\begin{figure}[ht]
\centering
\begin{tabular}{lccccccc}
\footnotesize{\textbf{Method}} & \footnotesize{\textbf{i=3}} & \footnotesize{\textbf{i=7}} & \footnotesize{\textbf{i=11}} & \footnotesize{\textbf{i=15}} & \footnotesize{\textbf{i=19}} & \footnotesize{\textbf{i=22}} &\footnotesize{ \textbf{i=29}} \\
\hline
\small{\tiny{$\vx_{\text{true}}$}} & 
\includegraphics[width=1.15cm]{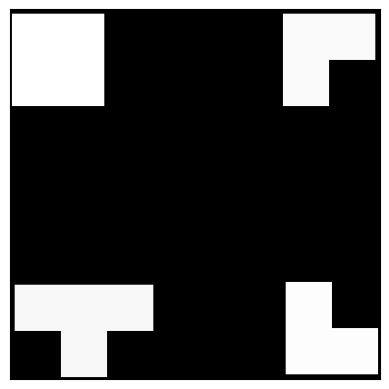} & 
\includegraphics[width=1.15cm]{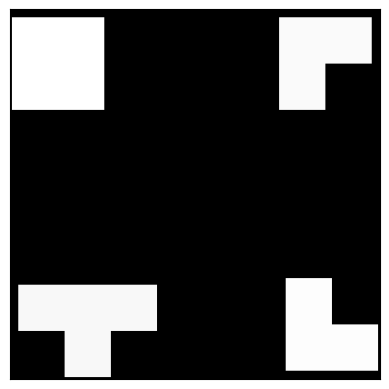} & 
\includegraphics[width=1.15cm]{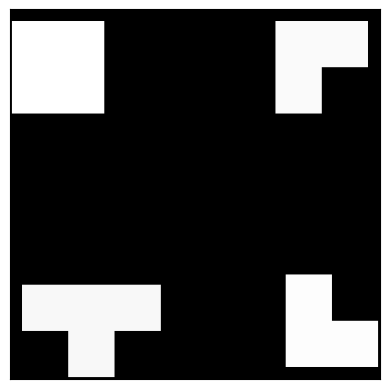} & 
\includegraphics[width=1.15cm]{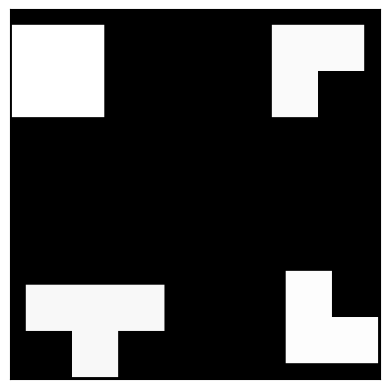} & 
\includegraphics[width=1.15cm]{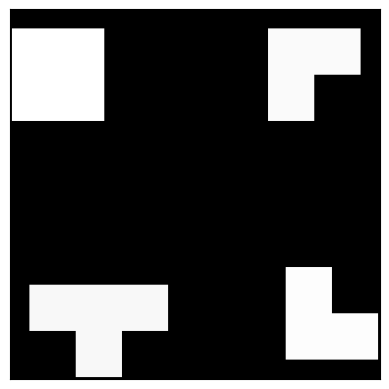} & 
\includegraphics[width=1.15cm]{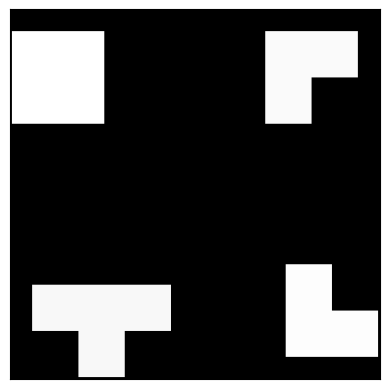} & 
\includegraphics[width=1.15cm]{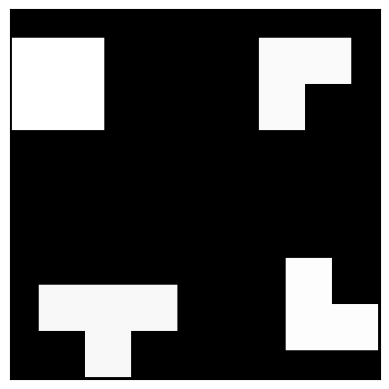} \\
\tiny{\textbf{AAO}} & 
\includegraphics[width=1.15cm]{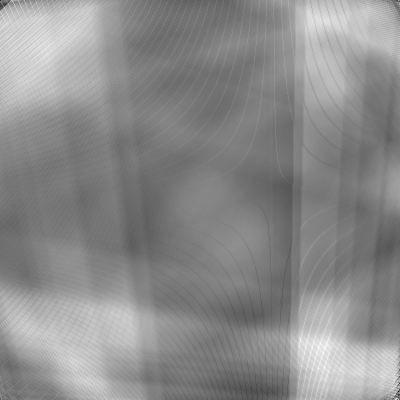} & 
\includegraphics[width=1.15cm]{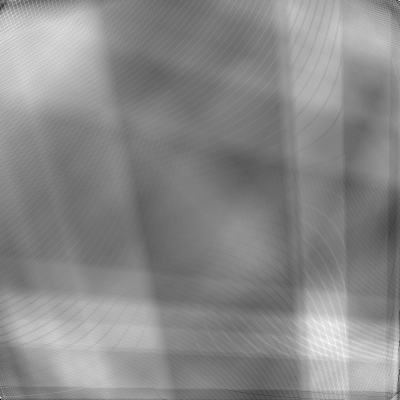} & 
\includegraphics[width=1.15cm]{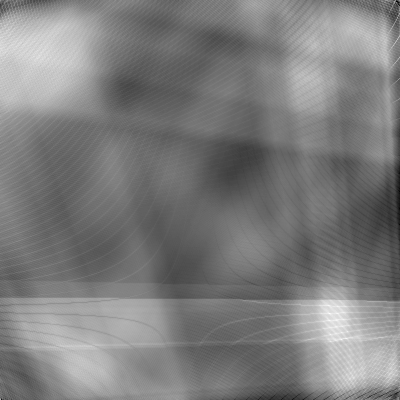} & 
\includegraphics[width=1.15cm]{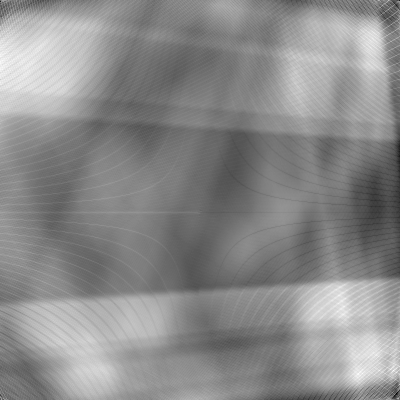} & 
\includegraphics[width=1.15cm]{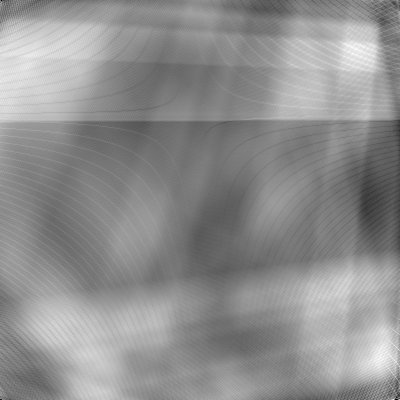} & 
\includegraphics[width=1.15cm]{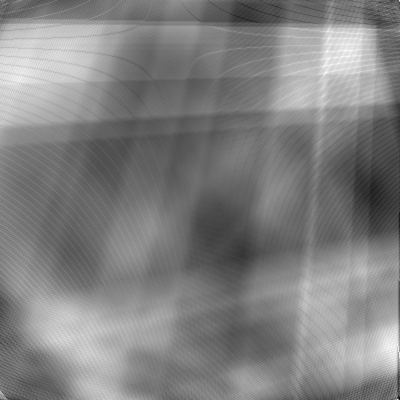} & 
\includegraphics[width=1.15cm]{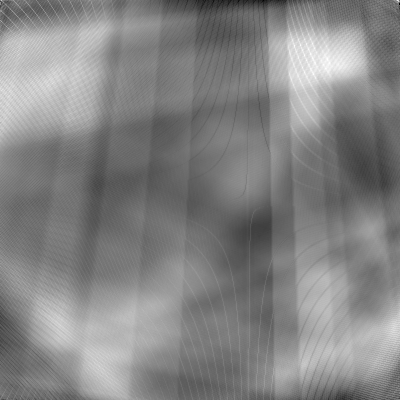} \\

\tiny{\textbf{AAO-ST}} & 
\includegraphics[width=1.15cm]{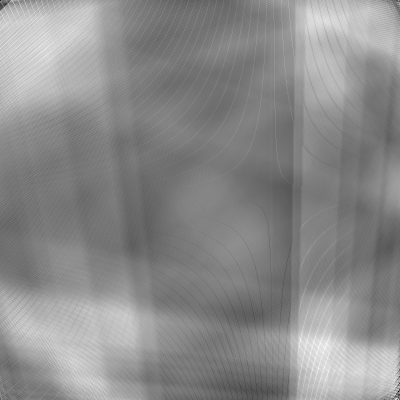} & 
\includegraphics[width=1.15cm]{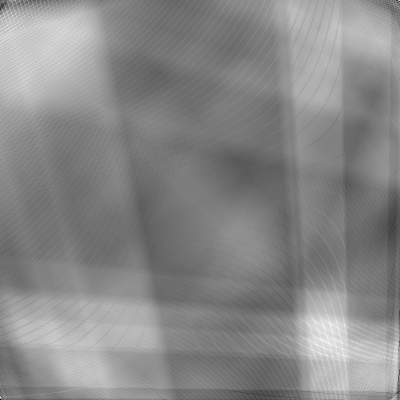} & 
\includegraphics[width=1.15cm]{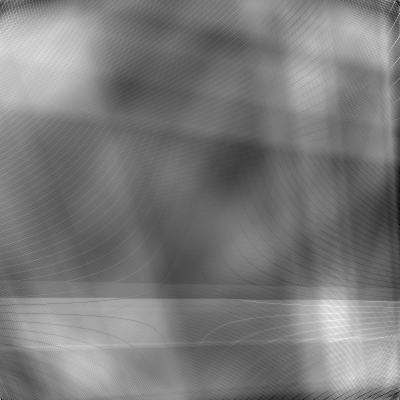} & 
\includegraphics[width=1.15cm]{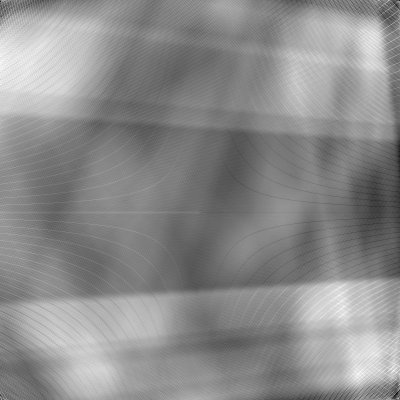} & 
\includegraphics[width=1.15cm]{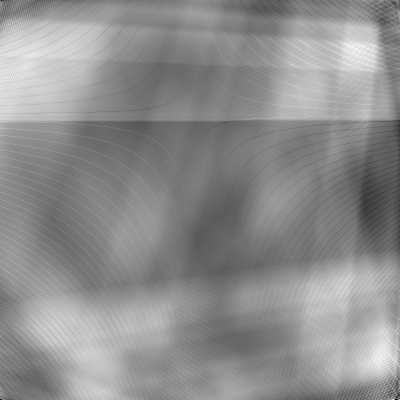} & 
\includegraphics[width=1.15cm]{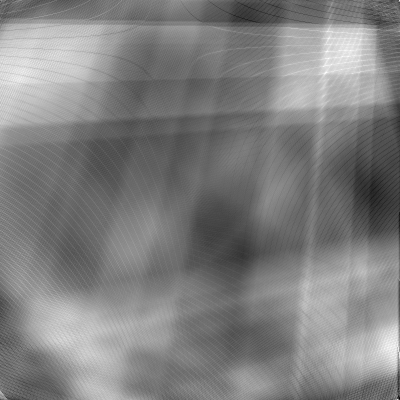} & 
\includegraphics[width=1.15cm]{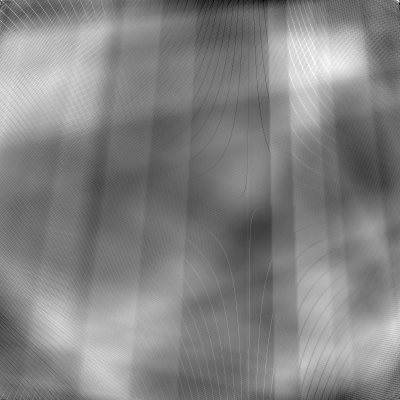} \\
\tiny{\textbf{AAO-OF}} & 
\includegraphics[width=1.15cm]{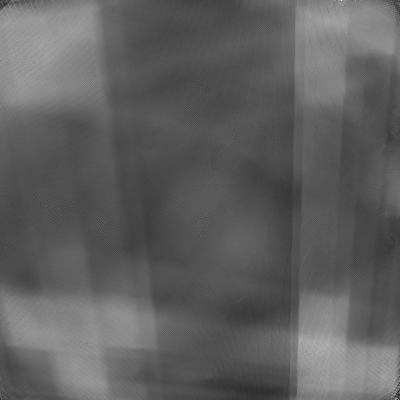} & 
\includegraphics[width=1.15cm]{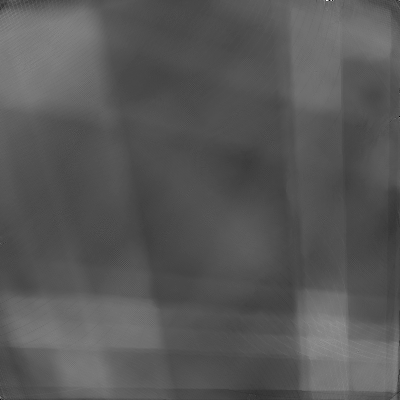} & 
\includegraphics[width=1.15cm]{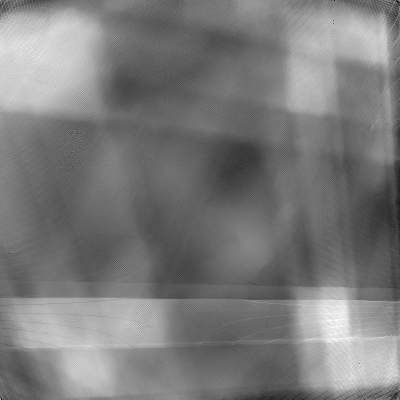} & 
\includegraphics[width=1.15cm]{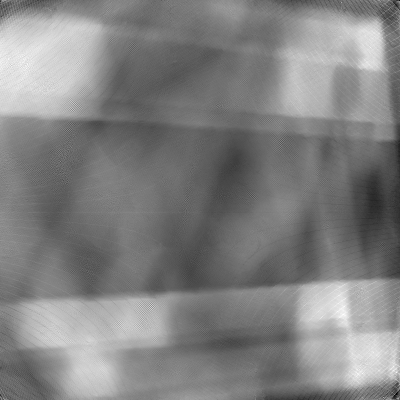} & 
\includegraphics[width=1.15cm]{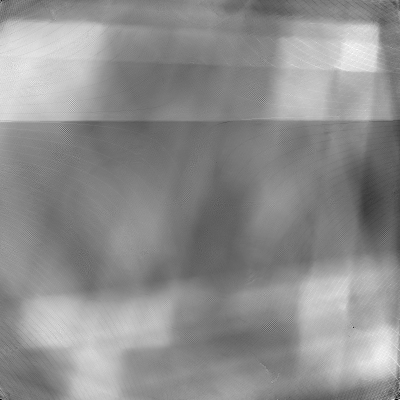} & 
\includegraphics[width=1.15cm]{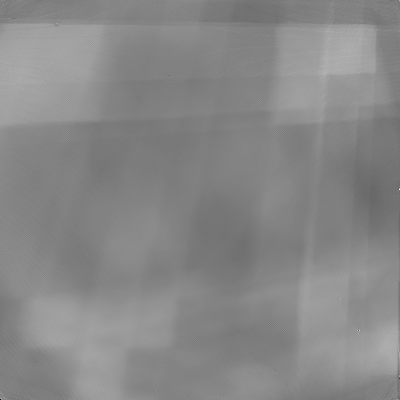} & 
\includegraphics[width=1.15cm]{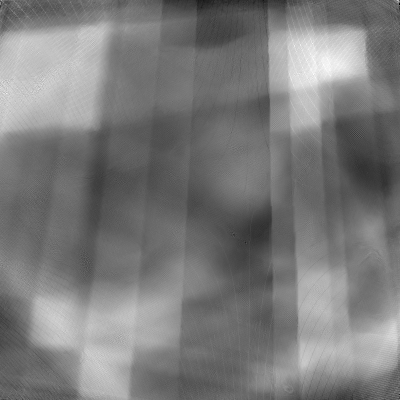} \\

\tiny{\textbf{IRKFS}} & 
\includegraphics[width=1.15cm]{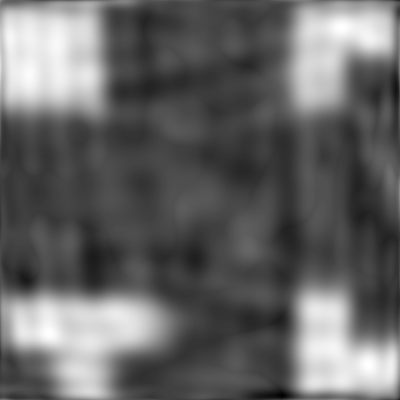} & 
\includegraphics[width=1.15cm]{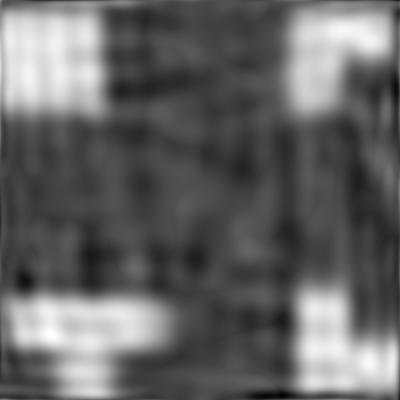} & 
\includegraphics[width=1.15cm]{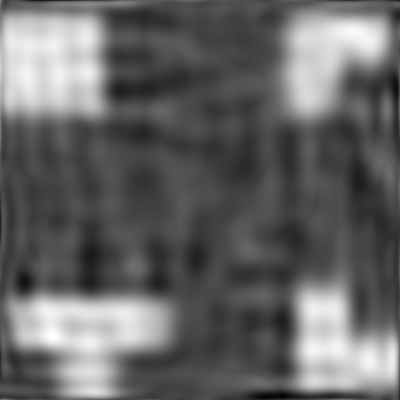} & 
\includegraphics[width=1.15cm]{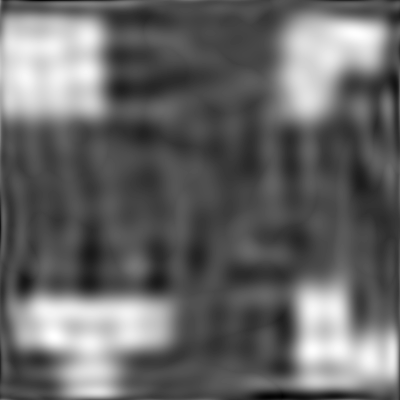} & 
\includegraphics[width=1.15cm]{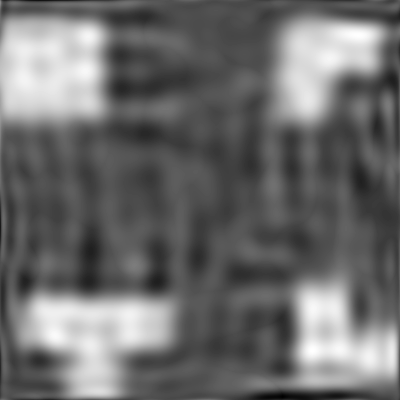} & 
\includegraphics[width=1.15cm]{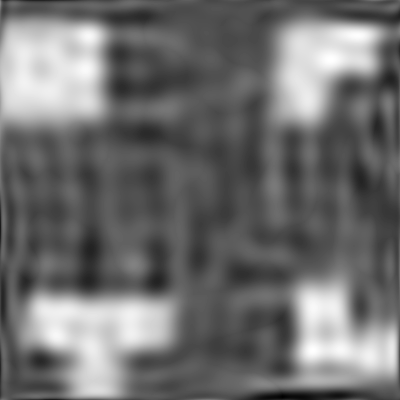} & 
\includegraphics[width=1.15cm]{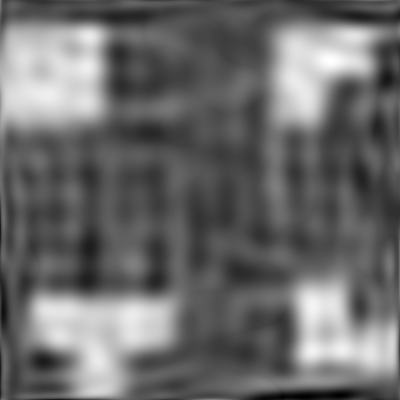} \\

\tiny{\textbf{IRKFS-M1}} & 
\includegraphics[width=1.15cm]{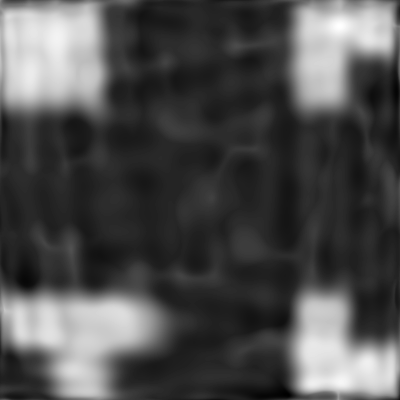} & 
\includegraphics[width=1.15cm]{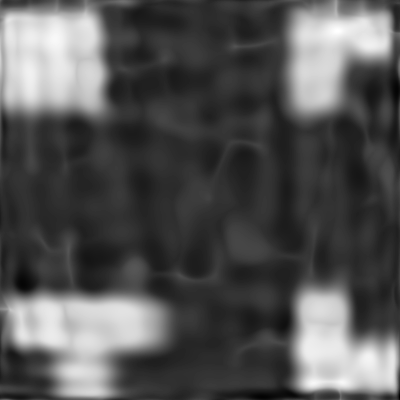} & 
\includegraphics[width=1.15cm]{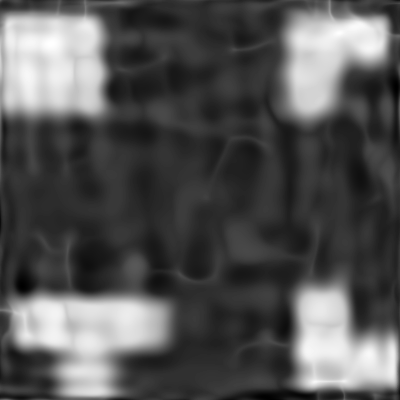} & 
\includegraphics[width=1.15cm]{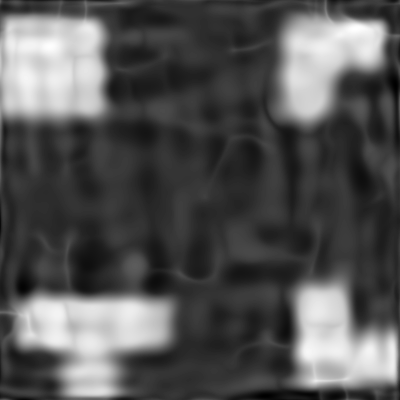} & 
\includegraphics[width=1.15cm]{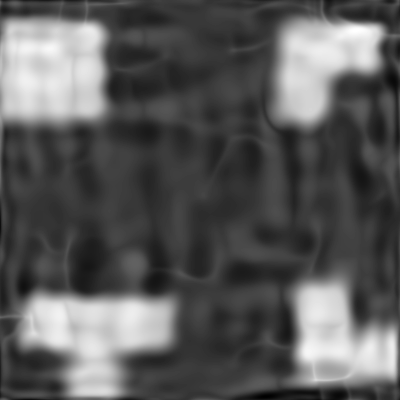} & 
\includegraphics[width=1.15cm]{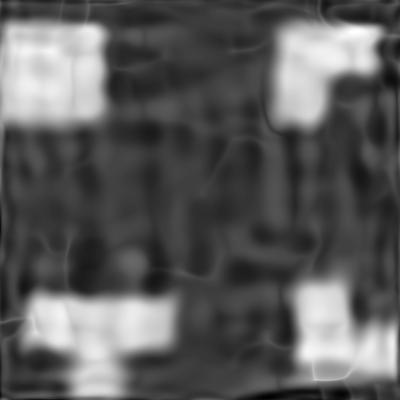} & 
\includegraphics[width=1.15cm]{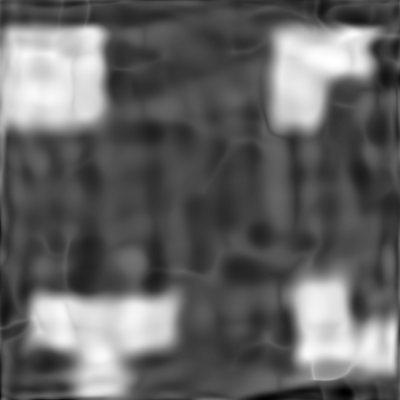} \\

\tiny{\textbf{IRKFS-M2}} & 
\includegraphics[width=1.15cm]{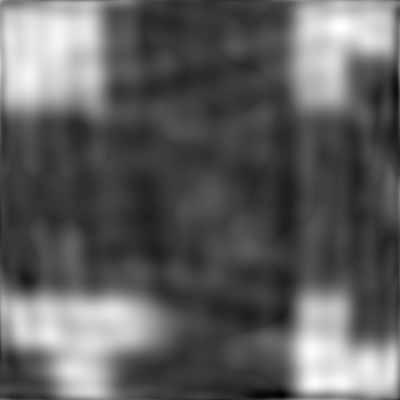} & 
\includegraphics[width=1.15cm]{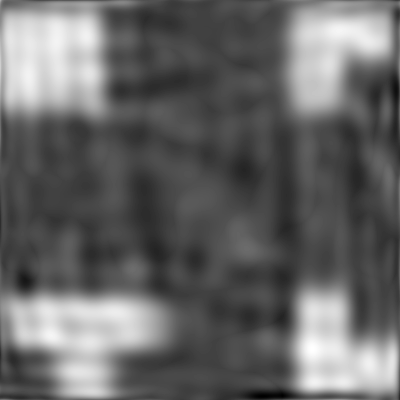} & 
\includegraphics[width=1.15cm]{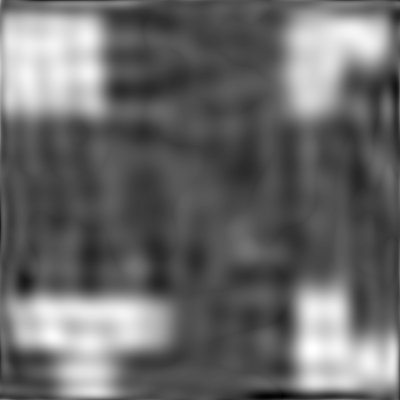} & 
\includegraphics[width=1.15cm]{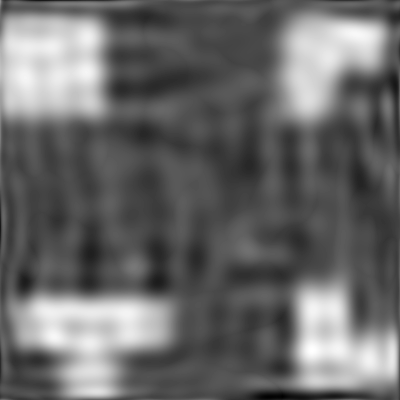} & 
\includegraphics[width=1.15cm]{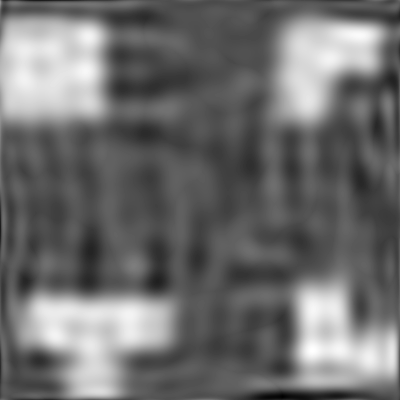} & 
\includegraphics[width=1.15cm]{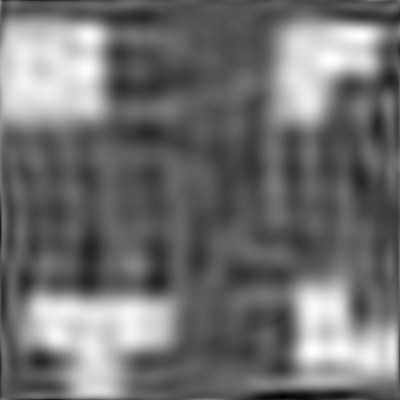} & 
\includegraphics[width=1.15cm]{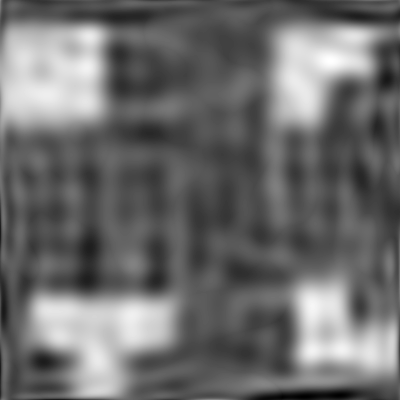} \\

\tiny{\textbf{IRKFS-M3}} & 
\includegraphics[width=1.15cm]{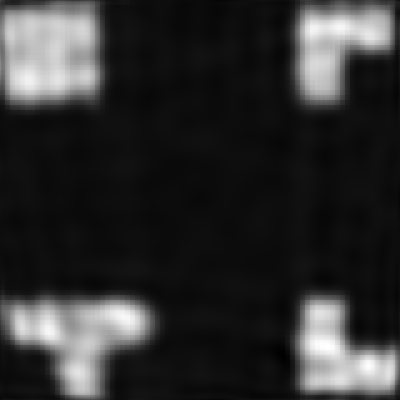} & 
\includegraphics[width=1.15cm]{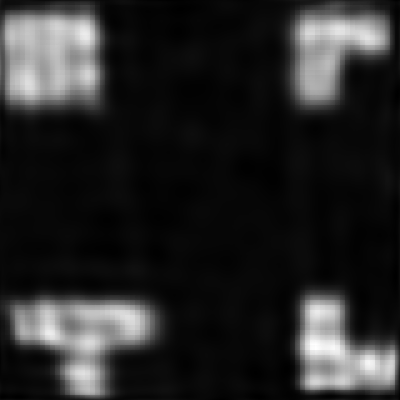} & 
\includegraphics[width=1.15cm]{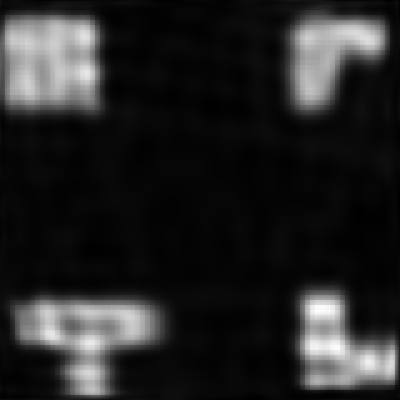} & 
\includegraphics[width=1.15cm]{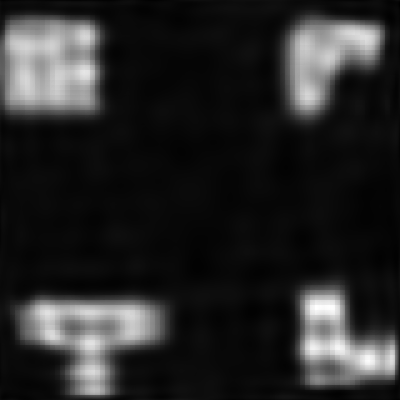} & 
\includegraphics[width=1.15cm]{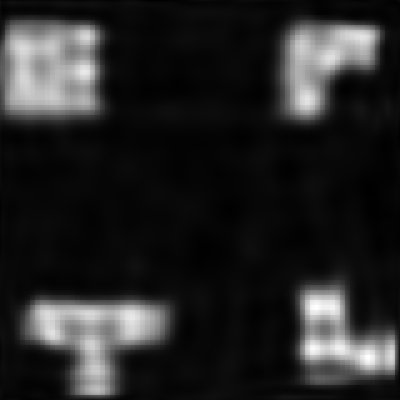} & 
\includegraphics[width=1.15cm]{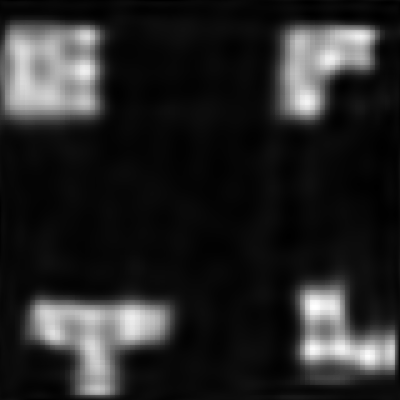} & 
\includegraphics[width=1.15cm]{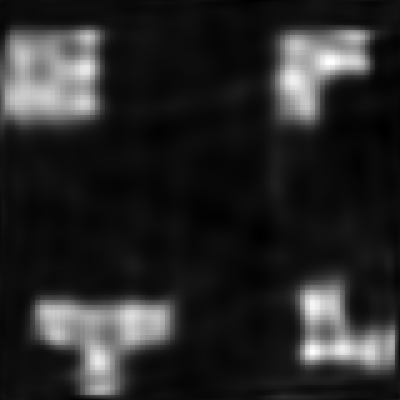} \\

\tiny{\textbf{EMIRKFS}} & 
\includegraphics[width=1.15cm]{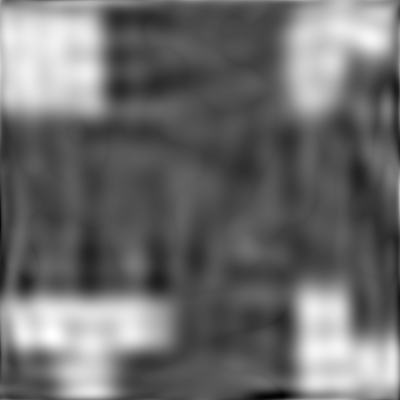} & 
\includegraphics[width=1.15cm]{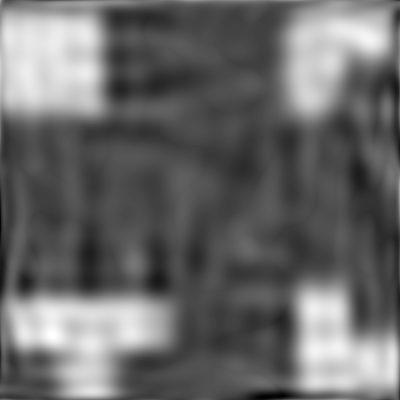} & 
\includegraphics[width=1.15cm]{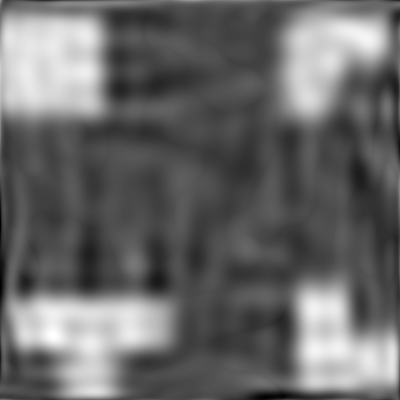} & 
\includegraphics[width=1.15cm]{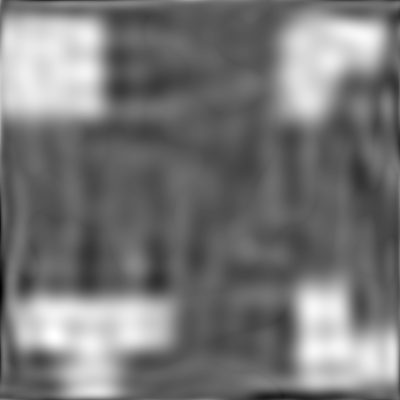} & 
\includegraphics[width=1.15cm]{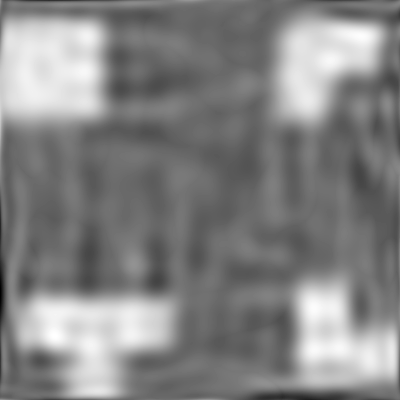} & 
\includegraphics[width=1.15cm]{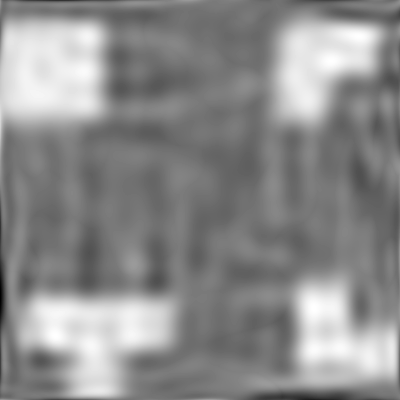} & 
\includegraphics[width=1.15cm]{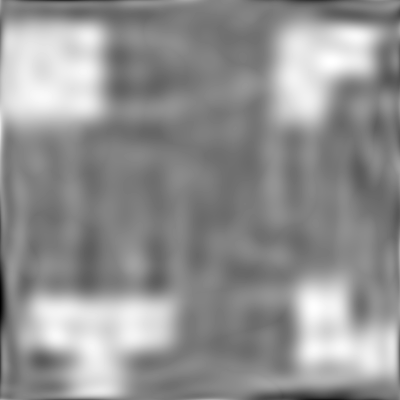} \\

\tiny{\textbf{EMIRKFS-M1}} & 
\includegraphics[width=1.15cm]{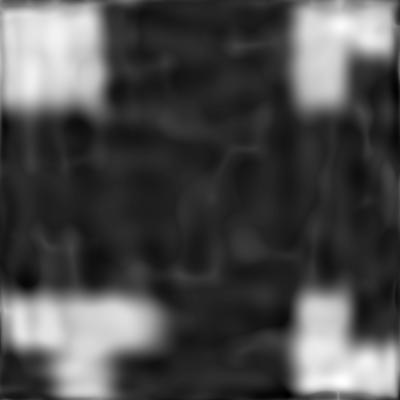} & 
\includegraphics[width=1.15cm]{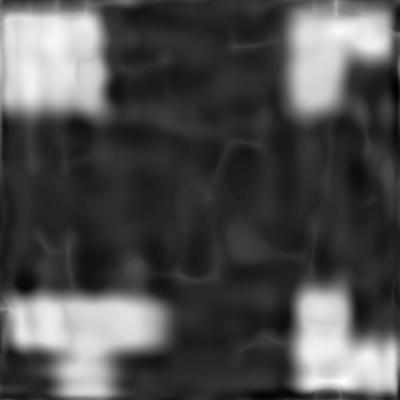} & 
\includegraphics[width=1.15cm]{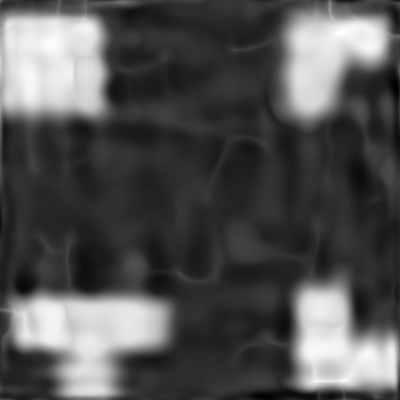} & 
\includegraphics[width=1.15cm]{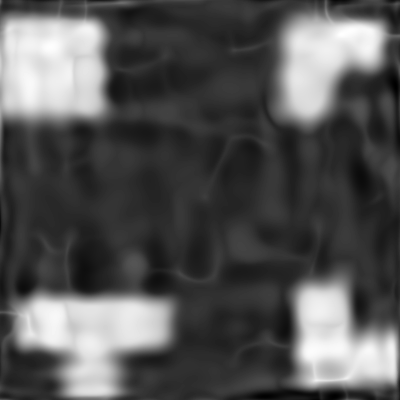} & 
\includegraphics[width=1.15cm]{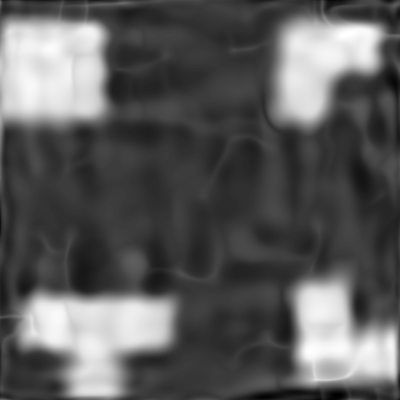} & 
\includegraphics[width=1.15cm]{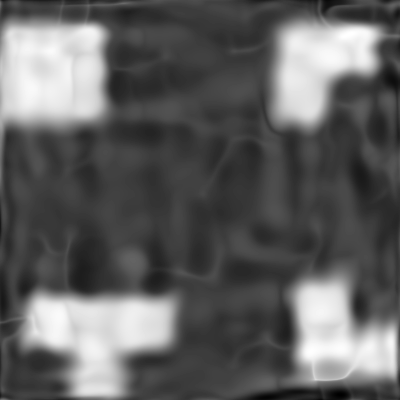} & 
\includegraphics[width=1.15cm]{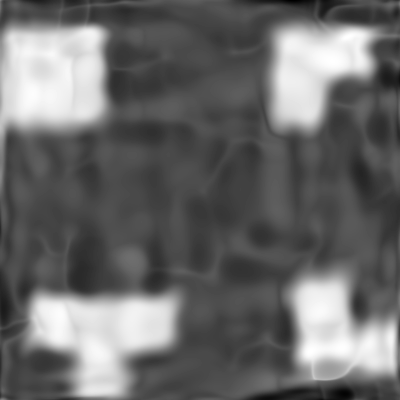} \\

\tiny{\textbf{EMIRKFS-M2}} & 
\includegraphics[width=1.15cm]{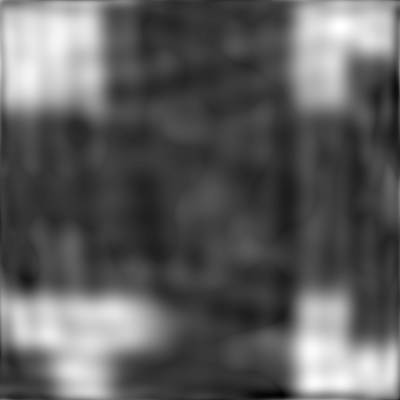} & 
\includegraphics[width=1.15cm]{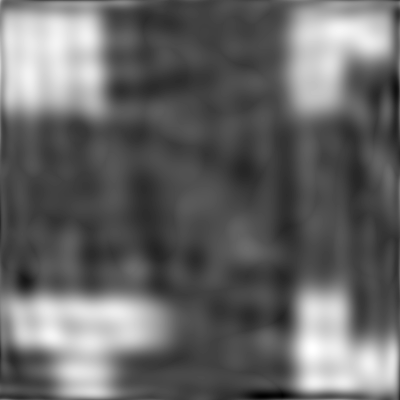} & 
\includegraphics[width=1.15cm]{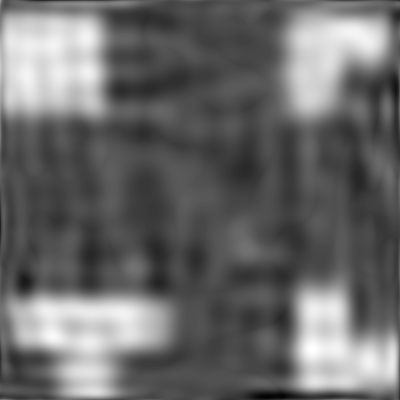} & 
\includegraphics[width=1.15cm]{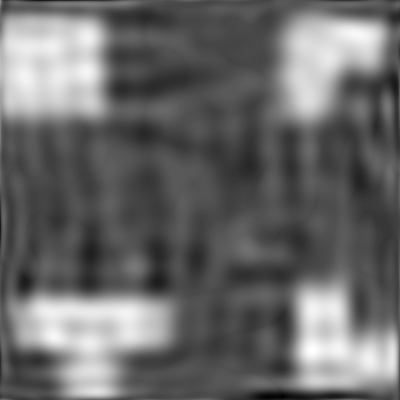} & 
\includegraphics[width=1.15cm]{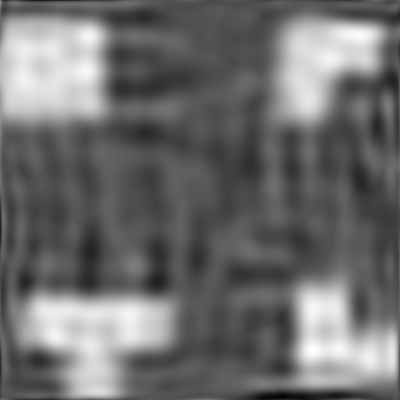} & 
\includegraphics[width=1.15cm]{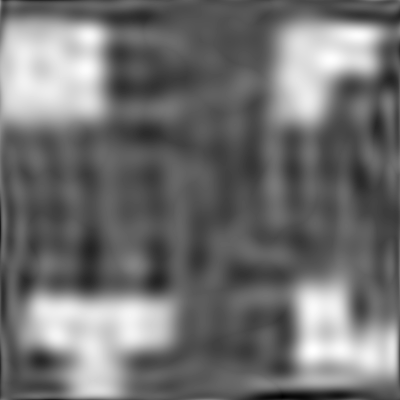} & 
\includegraphics[width=1.15cm]{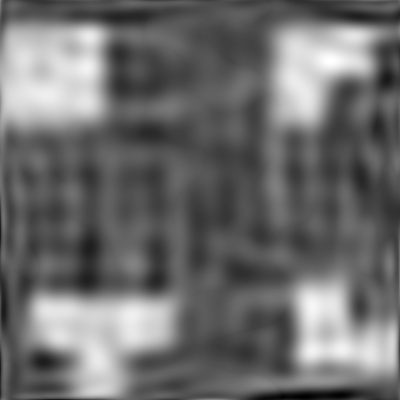} \\

\tiny{\textbf{EMIRKFS-M3}} & 
\includegraphics[width=1.15cm]{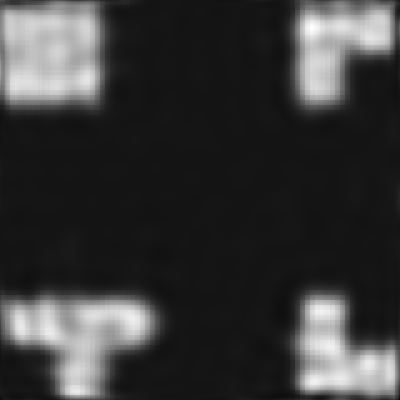} & 
\includegraphics[width=1.15cm]{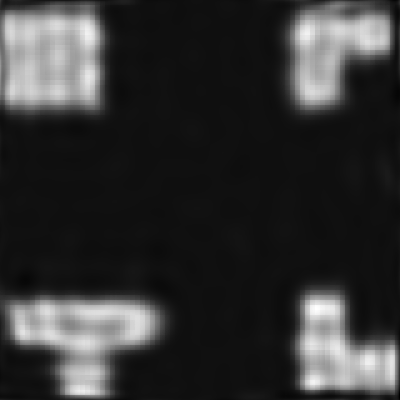} & 
\includegraphics[width=1.15cm]{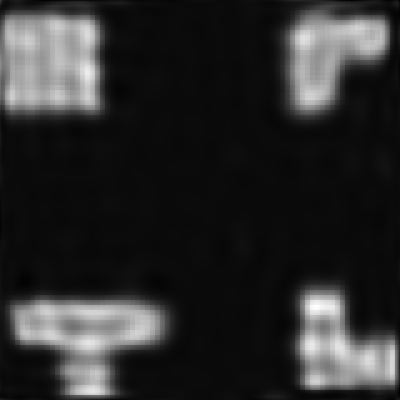} & 
\includegraphics[width=1.15cm]{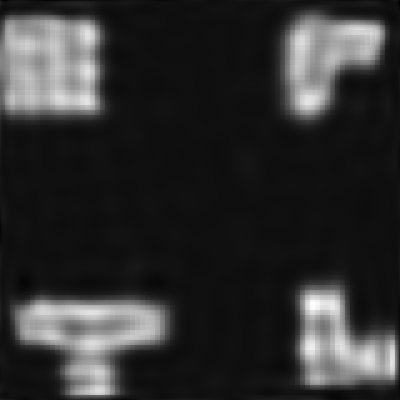} & 
\includegraphics[width=1.15cm]{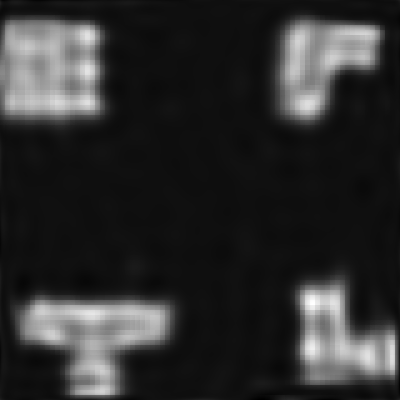} & 
\includegraphics[width=1.15cm]{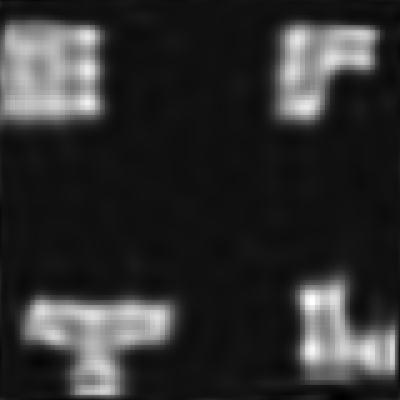} & 
\includegraphics[width=1.15cm]{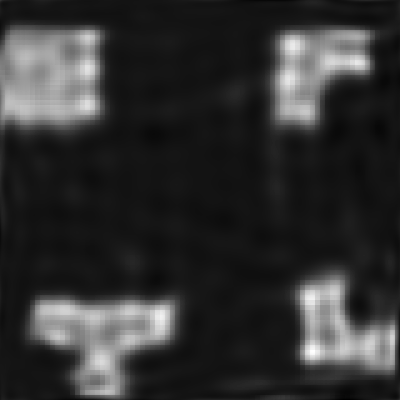} \\

\end{tabular}
\caption{Experiment 1: Truth image (Row 1) comparison to the image reconstruction of methods: AAO, AAO-ST, AAO-OF, IRKFS, IRKFS-M1, IRKFS-M2, IRFKS-M3, EMIRFKS, EMIRFKS-M1, EMIRKFS-M2, and EMIRFKS-M3 (Rows 2-12 from top to bottom) at time-steps $i\in\{3,7,11,15,19,22,29\}$ (Columns 1-7 from left to right).} 
\label{fig:method_comparison}
\end{figure}

For our first experiment, we simulate a 2D dynamic phantom with four blocks moving with speed within a bounded space
for 31 frames. Two of the four blocks move faster. The resulting phantom is a
3D array with dimensions of $400\times 400\times 31$. The phantom at time steps $t = 3, 7, 11, 15, 19, 22, 29$ are shown in the first row of Figure \ref{fig:method_comparison}, from left to right, respectively.

\paragraph{Forward Model and Projection Geometry} To generate the forward operator at each timestep, we use the Trips-Py toolbox \cite{pasha2024trips} with 5 projection angles with the methodology shown in \cite{pasha2023computational}. Each individual system matrix $\mathbf{H}_i$ has dimension $2825 \times 160000$ with observations $\mathbf{y}_i\in\mathbb{R}^{1695}$, and the complete forward model is represented by the global matrix $\mathbf{H}$  with size $87575 \times 4960000$ and corresponding vectorized sinograms from all time-steps $\mathbf{y}\in\mathbb{R}^{87575}$. To simulate realistic acquisition conditions, the respective sinogram data is corrupted with additive noise level $\sigma_{\text{NL}}=0.01$. For both M2 and M3 motion models we fix $\zeta=5$, and for M3 we select a patchsize of $z_x=z_y=5$.  

\begin{figure}[ht]
\centering
\begin{tabular}{lccccccc}
\footnotesize{\textbf{Method}} & \footnotesize{\textbf{i=3}} & \footnotesize{\textbf{i=7}} & \footnotesize{\textbf{i=11}} & \footnotesize{\textbf{i=15}} & \footnotesize{\textbf{i=19}} & \footnotesize{\textbf{i=22}} & \footnotesize{\textbf{i=29}} \\
\hline

\tiny{\textbf{AAO}} & 
\includegraphics[width=1.15cm]{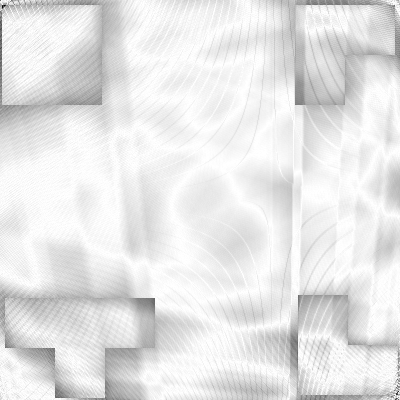} & 
\includegraphics[width=1.15cm]{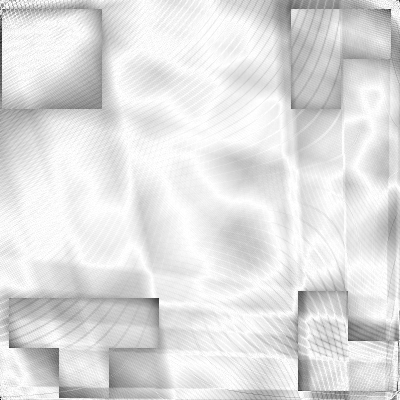} & 
\includegraphics[width=1.15cm]{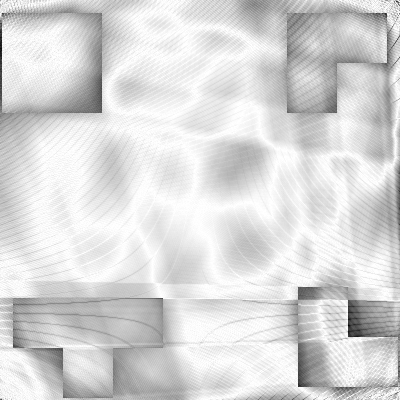} & 
\includegraphics[width=1.15cm]{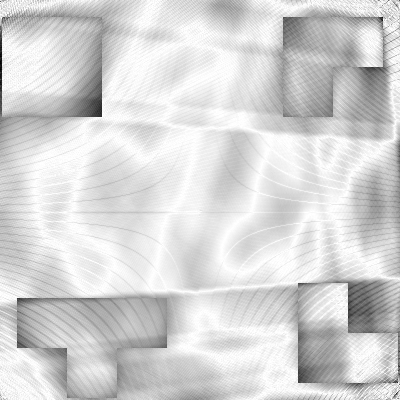} & 
\includegraphics[width=1.15cm]{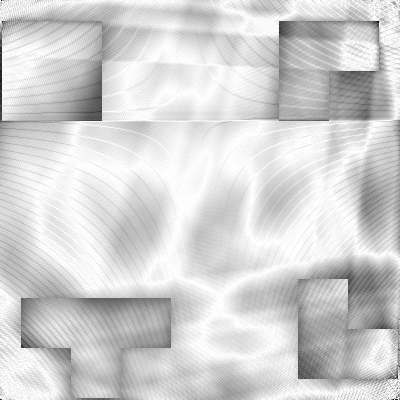} & 
\includegraphics[width=1.15cm]{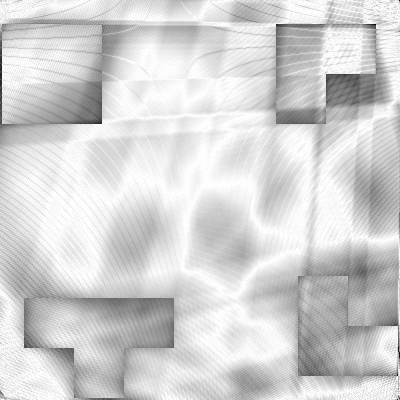} & 
\includegraphics[width=1.15cm]{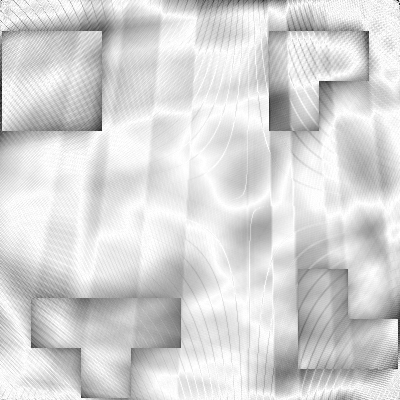} \\

\tiny{\textbf{AAO-ST}} & 
\includegraphics[width=1.15cm]{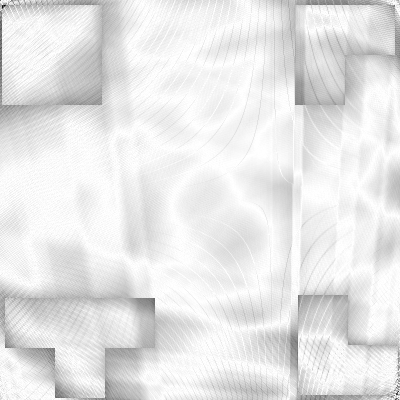} & 
\includegraphics[width=1.15cm]{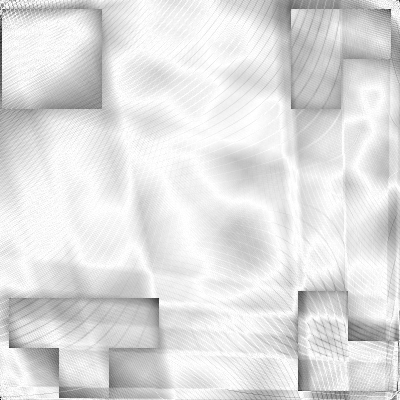} & 
\includegraphics[width=1.15cm]{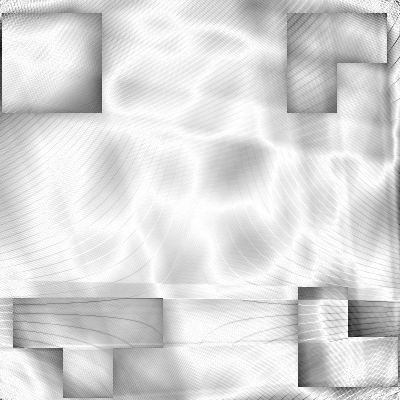} & 
\includegraphics[width=1.15cm]{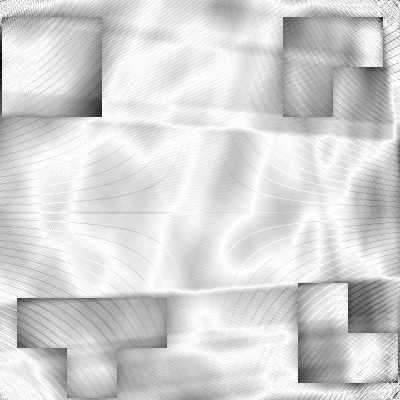} & 
\includegraphics[width=1.15cm]{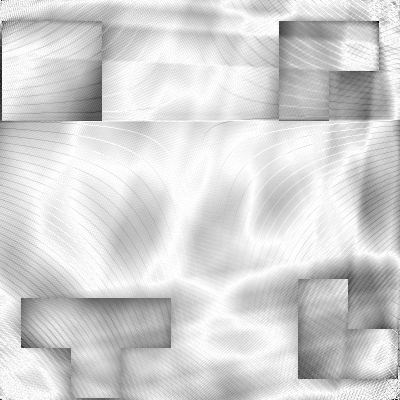} & 
\includegraphics[width=1.15cm]{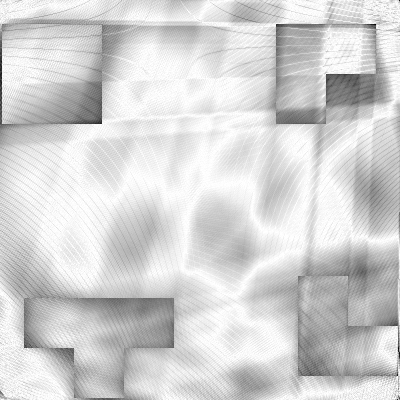} & 
\includegraphics[width=1.15cm]{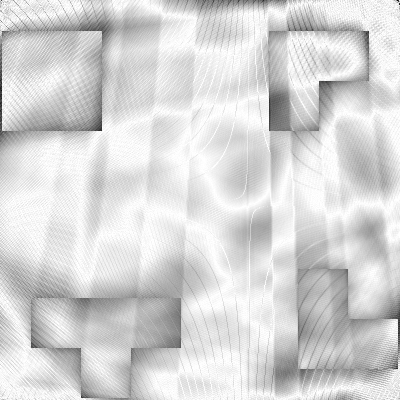} \\

\tiny{\textbf{AAO-OF}} & 
\includegraphics[width=1.15cm]{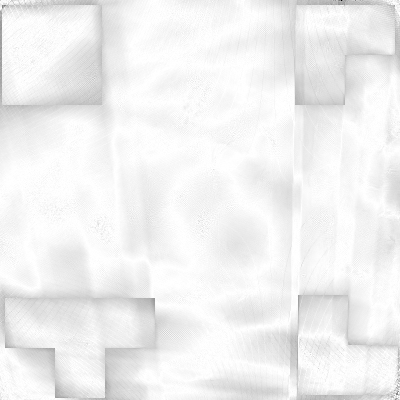} & 
\includegraphics[width=1.15cm]{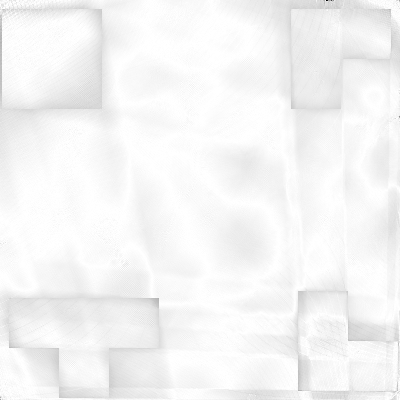} & 
\includegraphics[width=1.15cm]{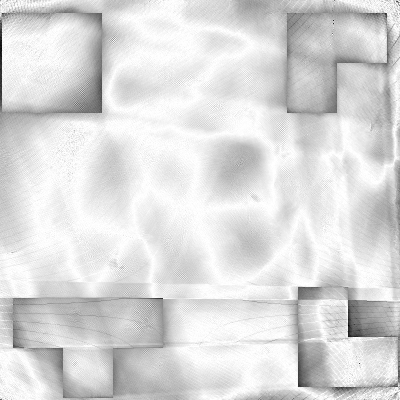} & 
\includegraphics[width=1.15cm]{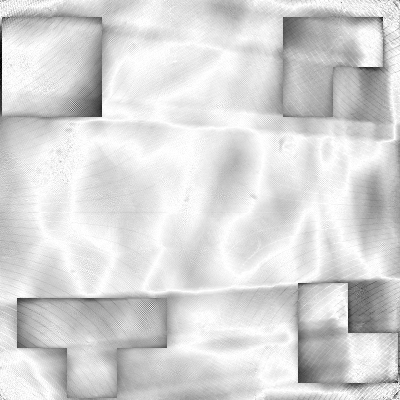} & 
\includegraphics[width=1.15cm]{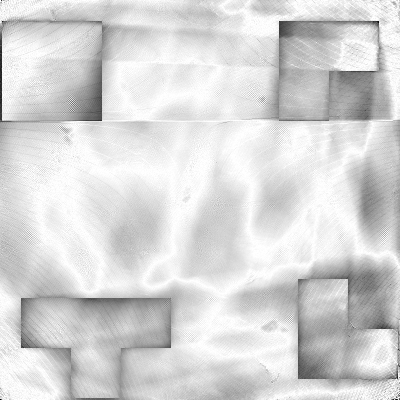} & 
\includegraphics[width=1.15cm]{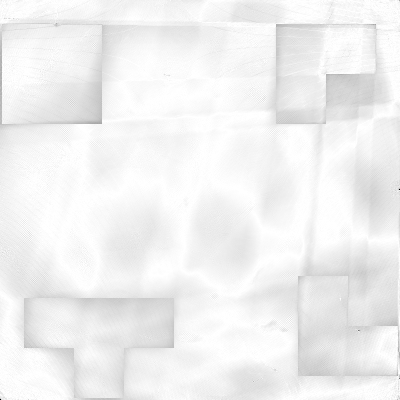} & 
\includegraphics[width=1.15cm]{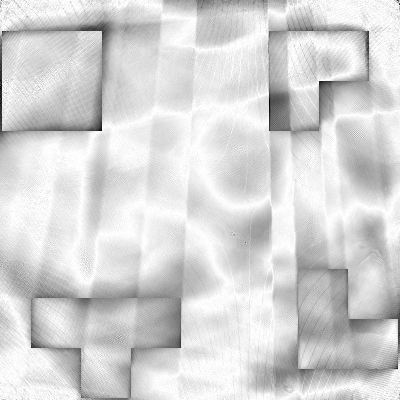} \\

\tiny{\textbf{IRKFS}} & 
\includegraphics[width=1.15cm]{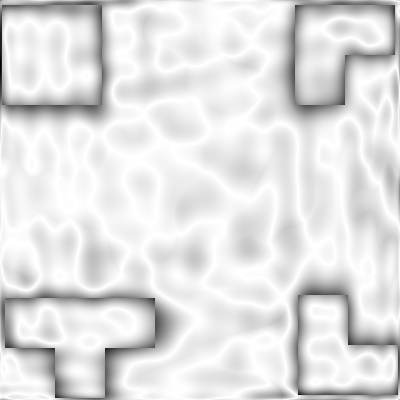} & 
\includegraphics[width=1.15cm]{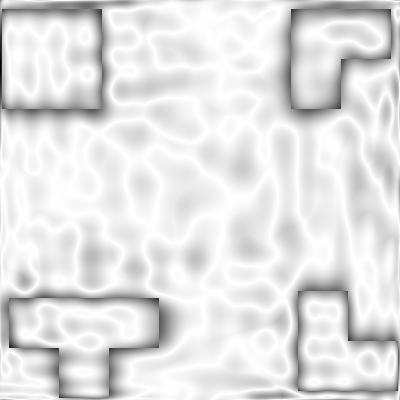} & 
\includegraphics[width=1.15cm]{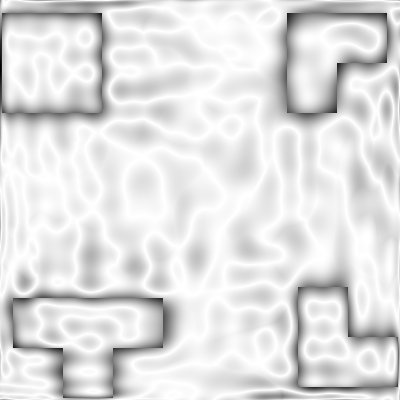} & 
\includegraphics[width=1.15cm]{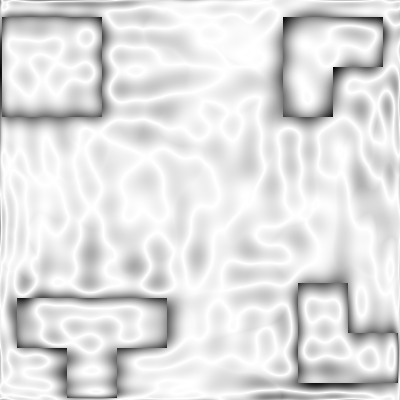} & 
\includegraphics[width=1.15cm]{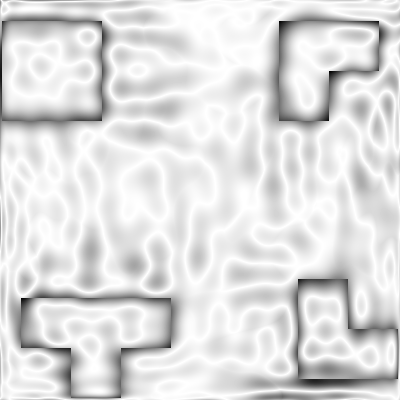} & 
\includegraphics[width=1.15cm]{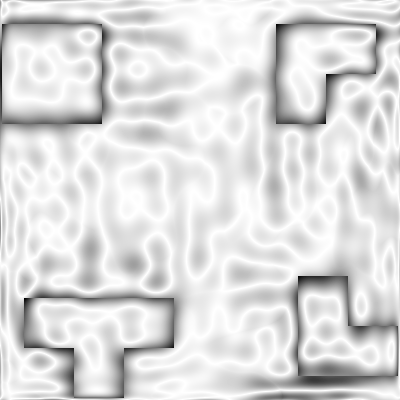} & 
\includegraphics[width=1.15cm]{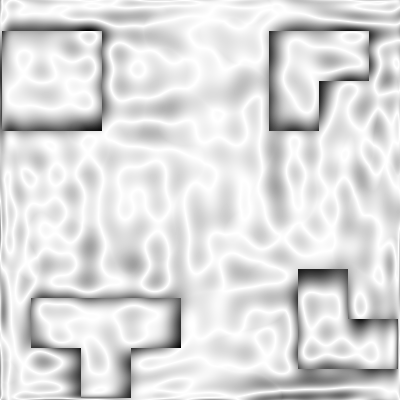} \\

\tiny{\textbf{IRKFS-M1}} & 
\includegraphics[width=1.15cm]{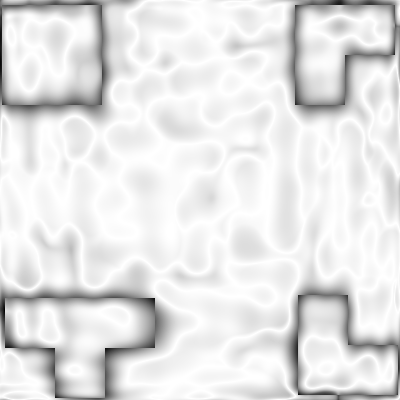} & 
\includegraphics[width=1.15cm]{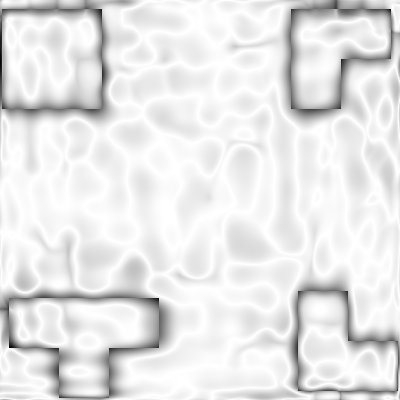} & 
\includegraphics[width=1.15cm]{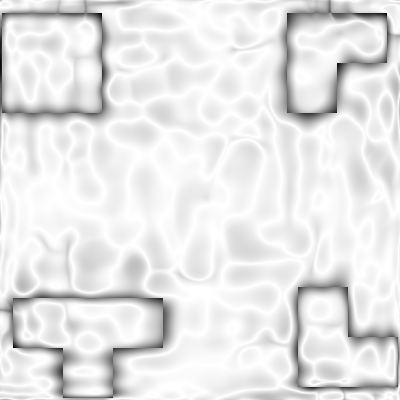} & 
\includegraphics[width=1.15cm]{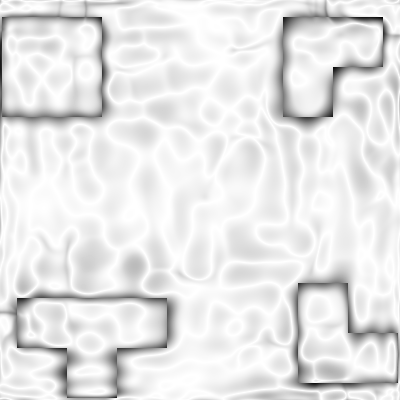} & 
\includegraphics[width=1.15cm]{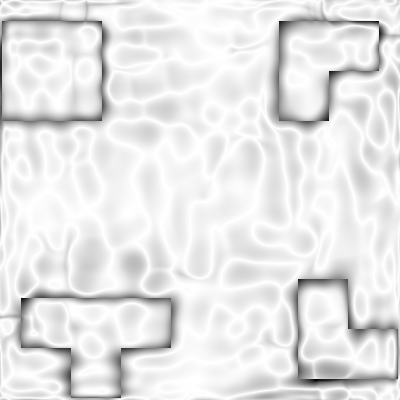} & 
\includegraphics[width=1.15cm]{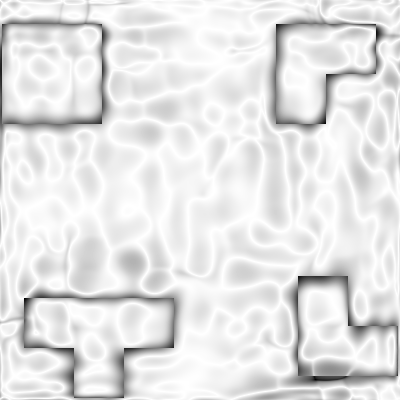} & 
\includegraphics[width=1.15cm]{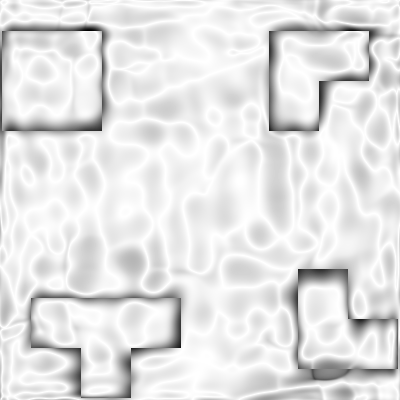} \\

\tiny{\textbf{IRKFS-M2}} & 
\includegraphics[width=1.15cm]{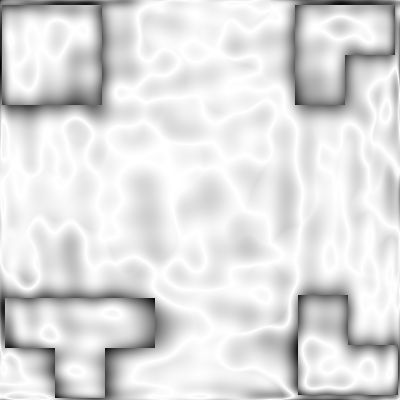} & 
\includegraphics[width=1.15cm]{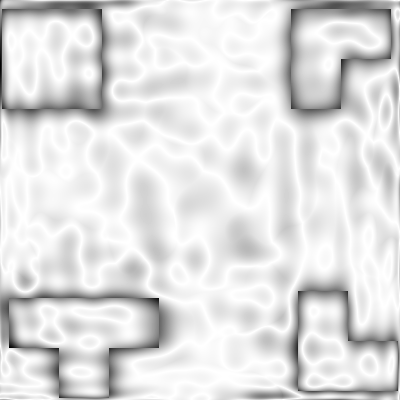} & 
\includegraphics[width=1.15cm]{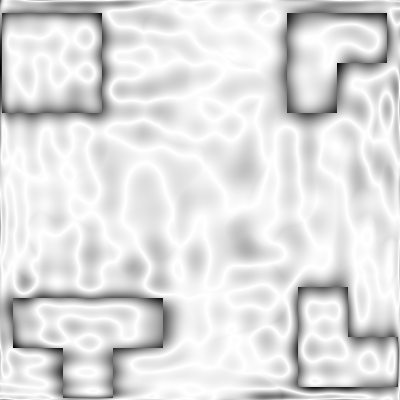} & 
\includegraphics[width=1.15cm]{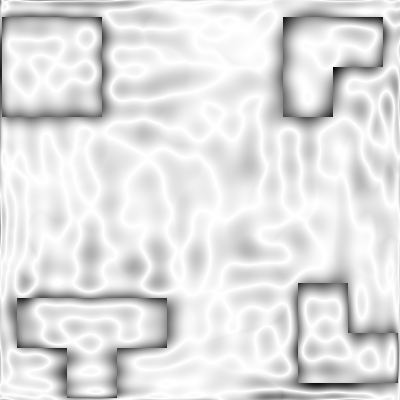} & 
\includegraphics[width=1.15cm]{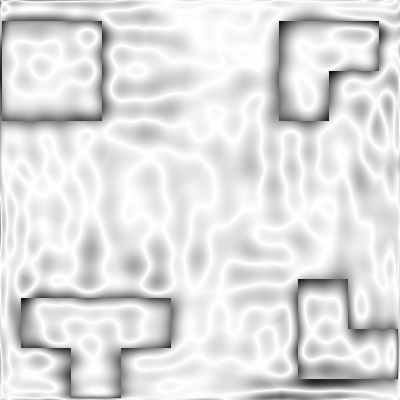} & 
\includegraphics[width=1.15cm]{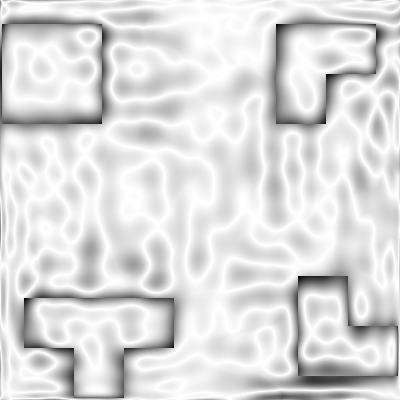} & 
\includegraphics[width=1.15cm]{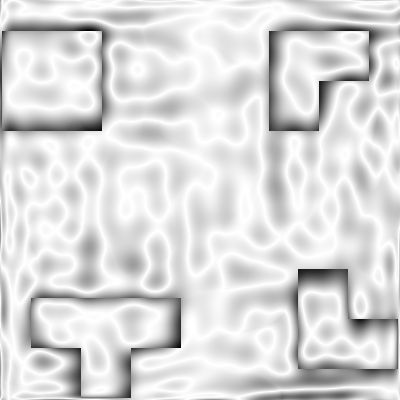} \\

\tiny{\textbf{IRKFS-M3}} & 
\includegraphics[width=1.15cm]{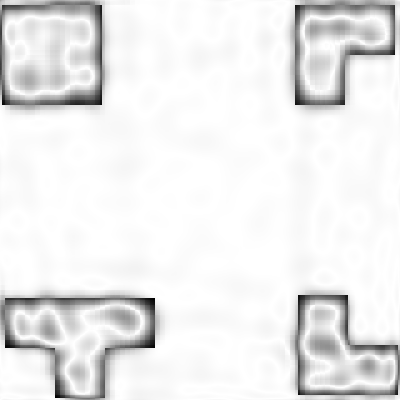} & 
\includegraphics[width=1.15cm]{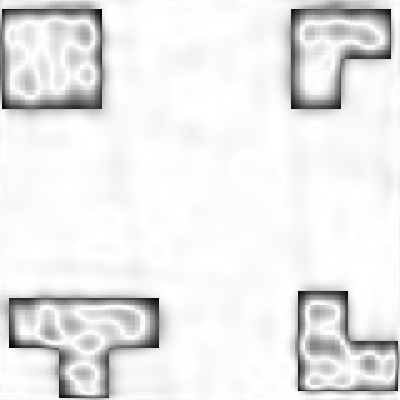} & 
\includegraphics[width=1.15cm]{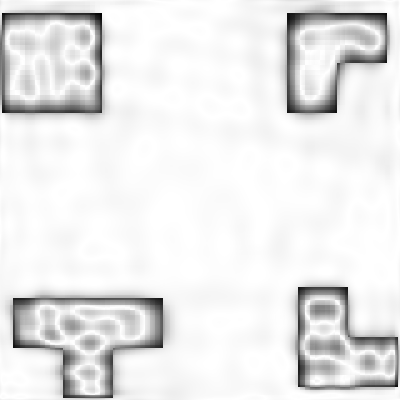} & 
\includegraphics[width=1.15cm]{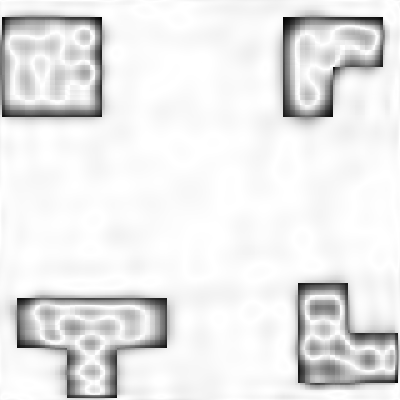} & 
\includegraphics[width=1.15cm]{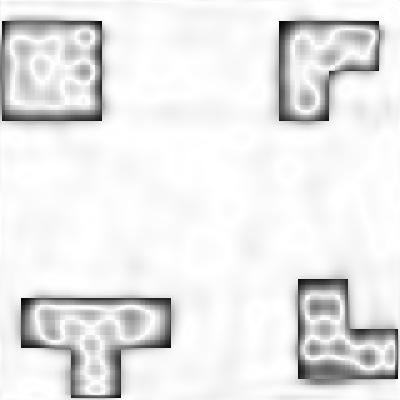} & 
\includegraphics[width=1.15cm]{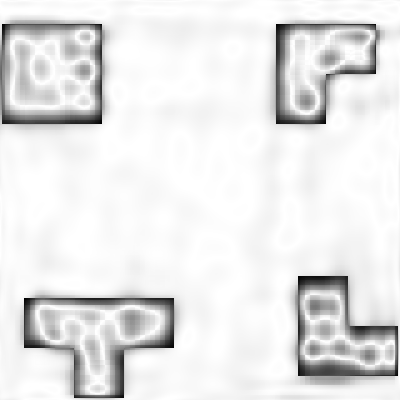} & 
\includegraphics[width=1.15cm]{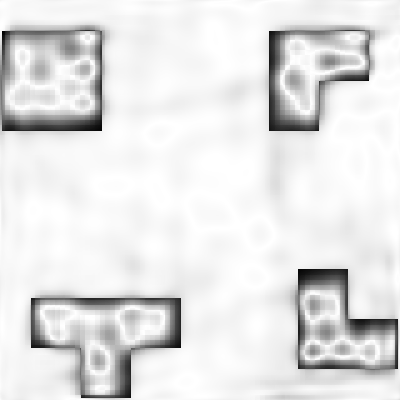} \\

\tiny{\textbf{EMIRKFS}} & 
\includegraphics[width=1.15cm]{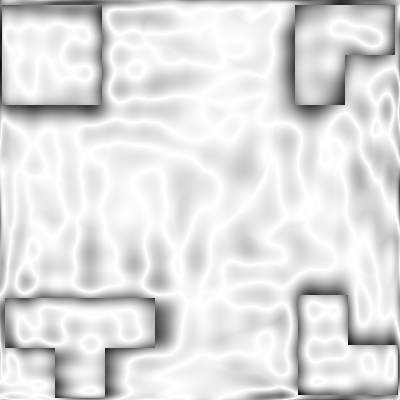} & 
\includegraphics[width=1.15cm]{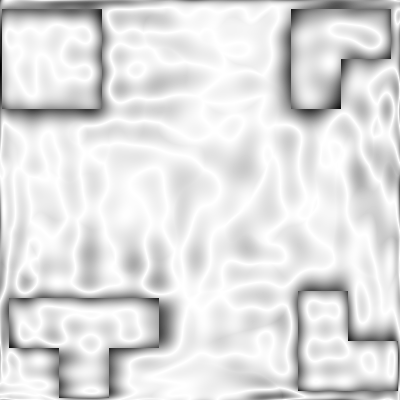} & 
\includegraphics[width=1.15cm]{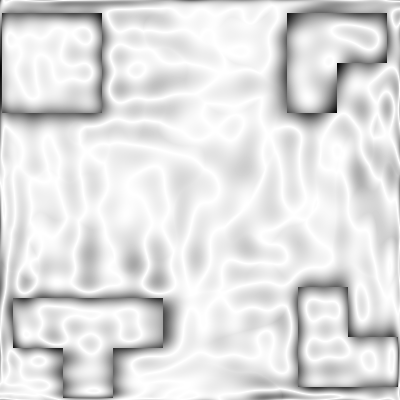} & 
\includegraphics[width=1.15cm]{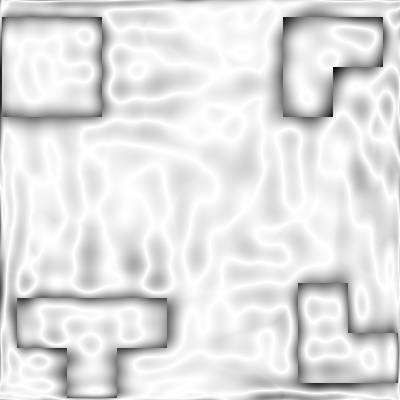} & 
\includegraphics[width=1.15cm]{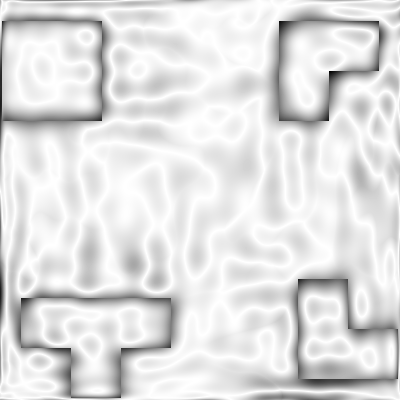} & 
\includegraphics[width=1.15cm]{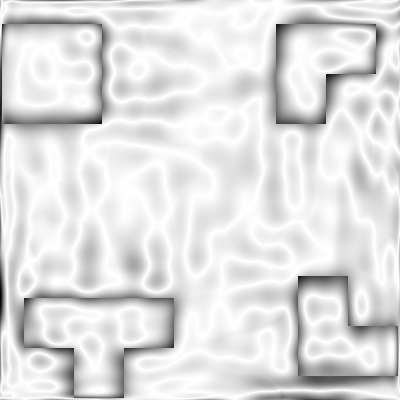} & 
\includegraphics[width=1.15cm]{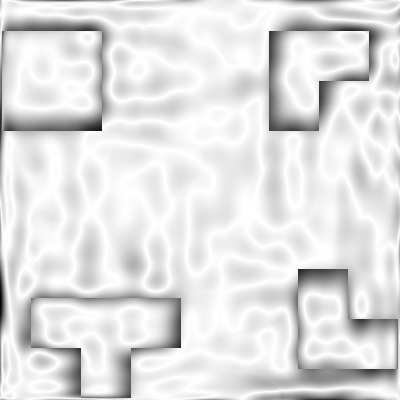} \\

\tiny{\textbf{EMIRKFS-M1}} & 
\includegraphics[width=1.15cm]{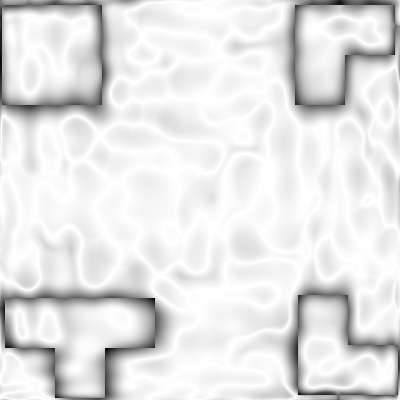} & 
\includegraphics[width=1.15cm]{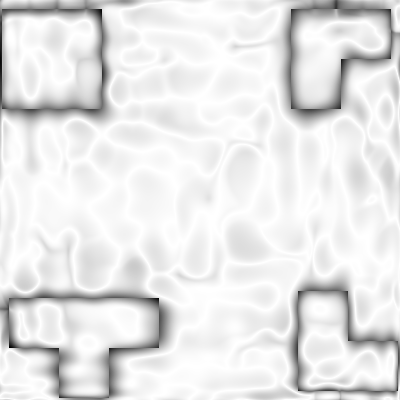} & 
\includegraphics[width=1.15cm]{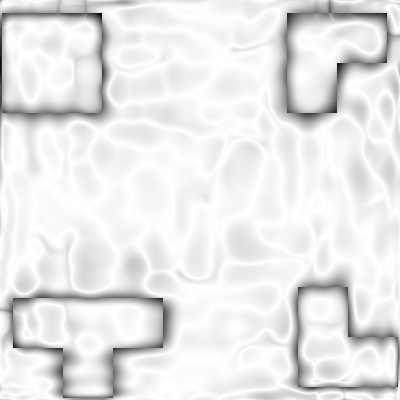} & 
\includegraphics[width=1.15cm]{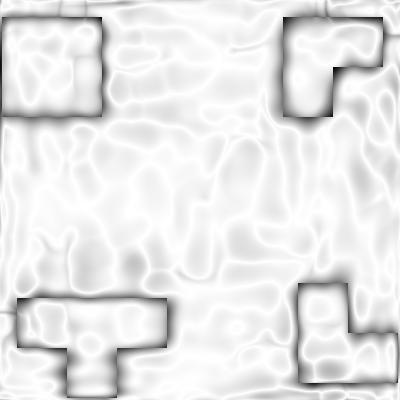} & 
\includegraphics[width=1.15cm]{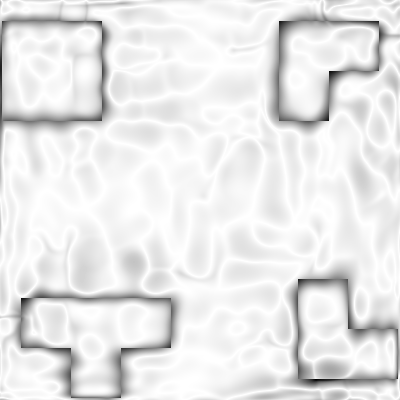} & 
\includegraphics[width=1.15cm]{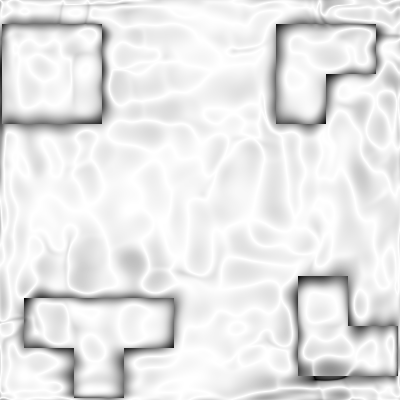} & 
\includegraphics[width=1.15cm]{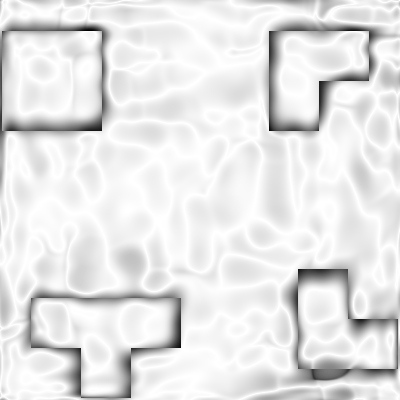} \\

\tiny{\textbf{EMIRKFS-M2}} & 
\includegraphics[width=1.15cm]{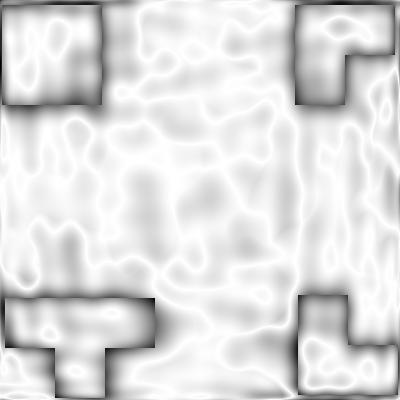} & 
\includegraphics[width=1.15cm]{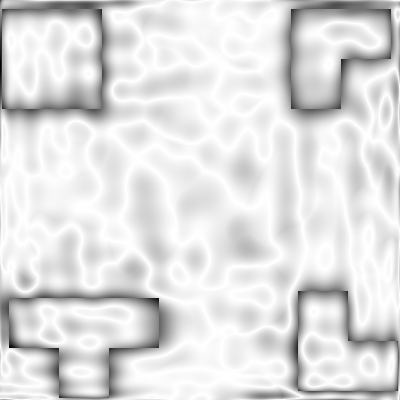} & 
\includegraphics[width=1.15cm]{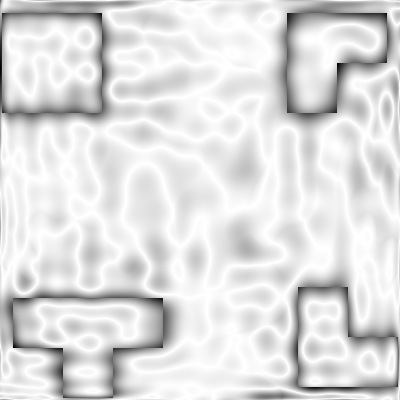} & 
\includegraphics[width=1.15cm]{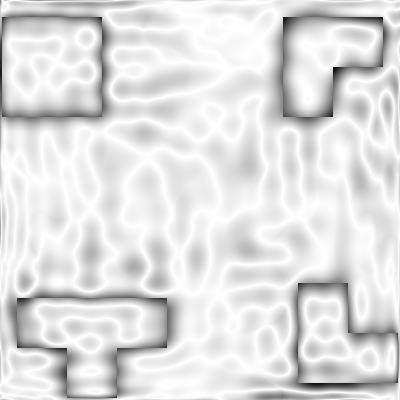} & 
\includegraphics[width=1.15cm]{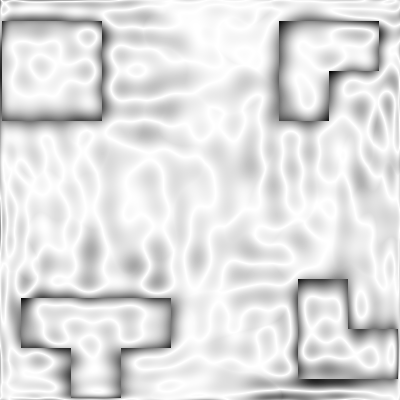} & 
\includegraphics[width=1.15cm]{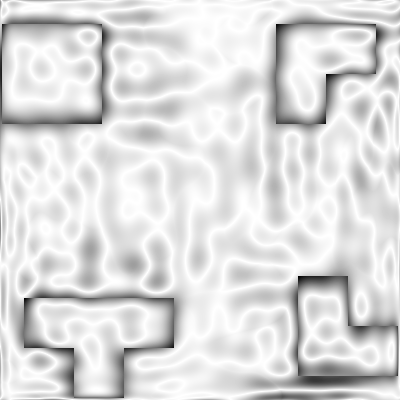} & 
\includegraphics[width=1.15cm]{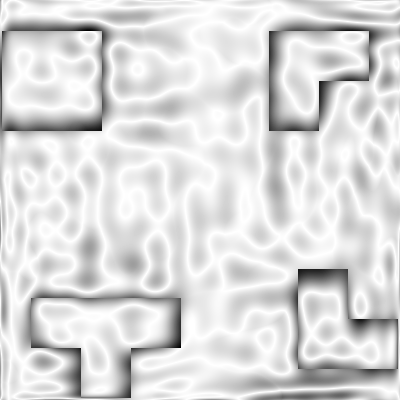} \\

\tiny{\textbf{EMIRKFS-M3}} & 
\includegraphics[width=1.15cm]{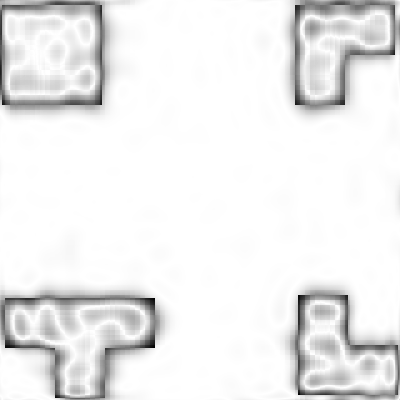} & 
\includegraphics[width=1.15cm]{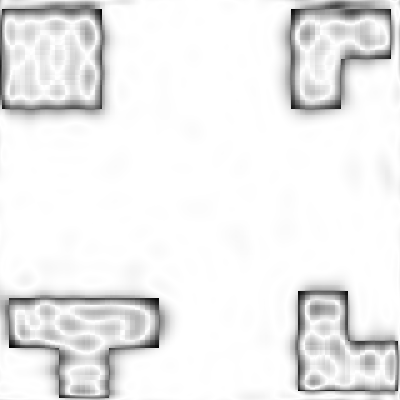} & 
\includegraphics[width=1.15cm]{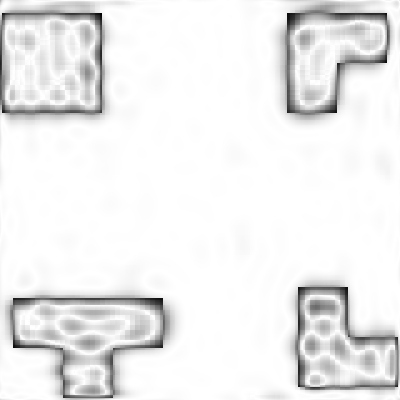} & 
\includegraphics[width=1.15cm]{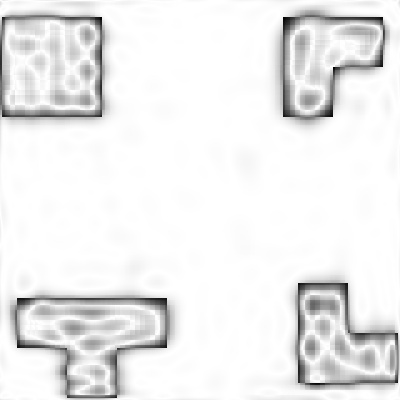} & 
\includegraphics[width=1.15cm]{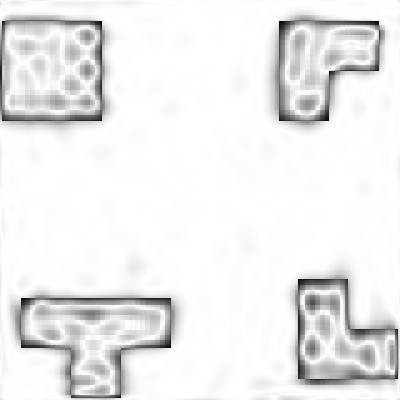} & 
\includegraphics[width=1.15cm]{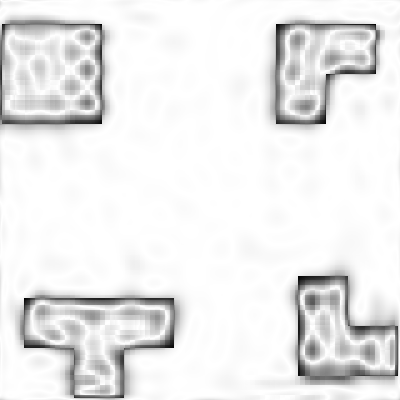} & 
\includegraphics[width=1.15cm]{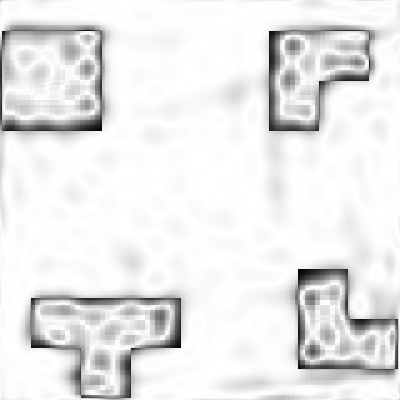} \\
\end{tabular}
\caption{Experiment 1: Reverse grayscale error image reconstructions of methods: AAO, AAO-ST, AAO-OF, IRKFS, IRKFS-M1, IRKFS-M2, IRFKS-M3, EMIRFKS, EMIRFKS-M1, EMIRKFS-M2, and EMIRFKS-M3 (Rows 1-11 from top to bottom) at time-steps $i\in\{3,7,11,15,19,22,29\}$ (Columns 1-7 from left to right).}
\label{fig:method_comparison2}
\end{figure}

\begin{figure}[ht]
\centering
\begin{tabular}{cc} 
\includegraphics[width=5.5cm]{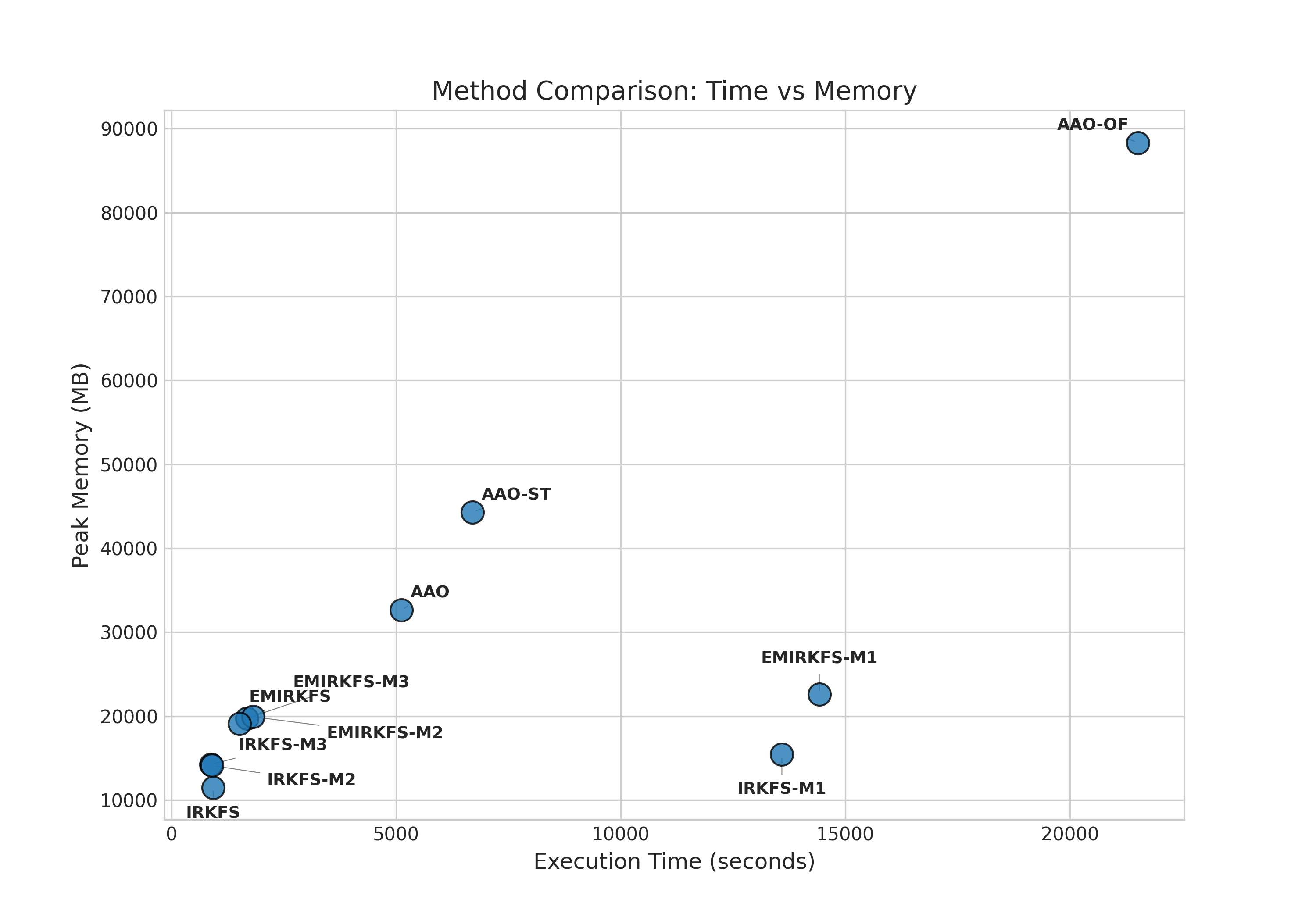} & 
\includegraphics[width=5.5cm]{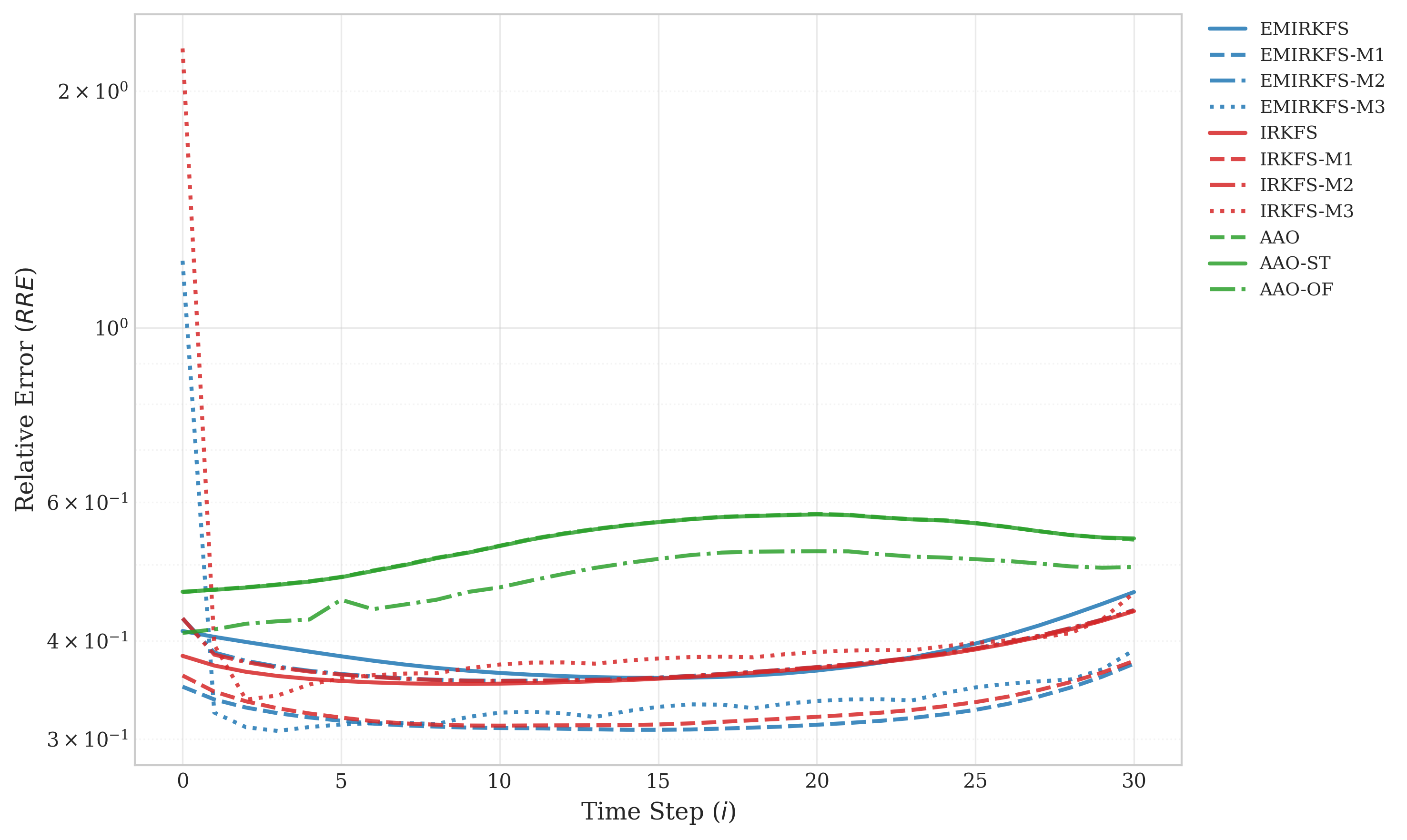}\\  
\end{tabular}
 \caption{Experiment 1: Memory vs. Time (left) and RRE (\ref{RRE}) (right) comparison across all methods (\ref{algorithms})}
 \label{fig:method_comparison3}
\end{figure}

We assess the reconstruction performance of the iterative methods with parameters fixed as follows: $r = 650$, $\alpha = 0.28$, $\ell = 11$, and $n_{\text{iter}}=2$. We initialize with $\boldsymbol{\Psi}^{\text{est}}_0=\mathbf{I}_r$, $\mathbf{Q}^{(0)}_i=1000\alpha^2\mathbf{I}_{n_s}$, $\mathbf{M}^{(0)}_i=\mathbf{I}_{n_s}$, $\mathbf{R}^{(0)}_i=\alpha^2\mathbf{I}_{m_t}$,  and as in \cite{hakkarainen2019undersampled} we set $\mathbf{x}^{\text{est}}_0=\mathbf{P}_r((\mathbf{H}_0\mathbf{P}_r)^\top\mathbf{H}_0\mathbf{P}_r+\frac{1}{\alpha^2}\mathbf{P}_r^\top\mathbf{P}_r)^{-1}\mathbf{y}_0$.
In addition to reconstruction quality, we compare the computational efficiency and memory requirements of method summarized in Algorithm \ref{Alg: EMIRKFS-M} to other alternative benchmark approaches summarized in \ref{algorithms}. We also plot negative image of the absolute value of the error images across the same methods and time-steps.
\paragraph{Results and Discussion} From visuals of Figures \ref{fig:method_comparison}, \ref{fig:method_comparison2} and \ref{fig:method_comparison3}, we can see that the rows of the AAO variants fail to identify the moving blocks given restriction to just 5 projection angles. While the IRKFS, IRKFS-M1, IRKFS-M2, EMIRKFS, EMIRKFS-M1, and EMIRKFS-M2 methods show improved reconstructions, dominant artifacts and noise still remain in background of their reconstructions. This is a contrast to the sharper reconstructions given by both IRKFS-M3 and EMIRKFS-M3 methods utilizing the M3 motion models, demonstrating the models capacity to smooth background noise. When considering both visual quality and relative error, the EMIRKFS-M3 method is the overall winner. 
In Figure \ref{fig:method_comparison3} we see that iterative methods that exclude motion model M1 have superior computational speed and reduced memory consumption in comparison to all other methods. 

Figure \ref{fig:method_comparison3} demonstrates the computed computational time and memory across all methods along with the RRE for the 31 image estimates across all methods.
\subsection{Experiment 2: MNIST Dataset}
\label{ssec:exp_phantom_2}
We now proceed to evaluate the implementation of the EMIRKFS-M method described in Algorithm \ref{Alg: EMIRKFS-M} using the Dynamic MNIST dataset. This dataset comprises 10,000 video sequences, each consisting of 20 frames with a resolution of $64 \times 64$ pixels. Each frame contains two moving handwritten digits. For the purposes of our study, we extract a single sequence and select 13 grayscale frames from the total 20 as the ground truth reference. The resulting phantom is a 3D array with dimensions of $64\times 64\times 13$.
\paragraph{Forward Model and Projection Geometry}
The forward operator is constructed using the \texttt{Trips-Py} toolbox \cite{pasha2024trips}, employing 11 projection angles per image. Each individual system matrix is defined as $\mathbf{H}_i$ has dimension $450 \times 4096$ with observations $\mathbf{y}_i\in\mathbb{R}^{450}$, and the complete forward model is represented by the global matrix $\mathbf{H}$ with dimension $5850 \times 53248$ and corresponding vectorized sinograms from all time-steps $\mathbf{y}\in\mathbb{R}^{5850}$. To simulate realistic acquisition conditions, the respective sinogram data is corrupted with additive noise $\sigma_{\text{NL}}=0.01$. For both M2 and M3 motion models we fix $\zeta=7$, and for M3 we select a patchsize of $z_x=z_y=2$.

We assess the reconstruction performance of the of the iterative methods with parameters fixed as follows: $r=1000$, $\alpha = 0.375$, $\ell =2.9$, and $n_{\text{iter}}=5$. 
We initialize with $\boldsymbol{\Psi}^{\text{est}}_0=\mathbf{I}_r$, $\mathbf{Q}^{(0)}_i=\alpha^2\mathbf{I}_{n_s}$, $\mathbf{M}^{(0)}_i=\mathbf{I}_{n_s}$, $\mathbf{R}^{(0)}_i=\alpha^2\mathbf{I}_{m_t}$,  and as in \cite{hakkarainen2019undersampled} we set $\mathbf{x}^{\text{est}}_0=\mathbf{P}_r((\mathbf{H}_0\mathbf{P}_r)^\top\mathbf{H}_0\mathbf{P}_r+\frac{1}{\alpha^2}\mathbf{P}_r^\top\mathbf{P}_r)^{-1}\mathbf{y}_0$.

\begin{figure}[ht]
\centering
\begin{tabular}{lccccc}
\textbf{Method} & \textbf{i=2} & \textbf{i=5} & \textbf{i=7} & \textbf{i=9} & \textbf{i=11} \\
\hline
\tiny{\textbf{$\vx_{\text{true}}$}} & 
\includegraphics[width=1.15cm]{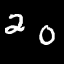} & 
\includegraphics[width=1.15cm]{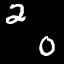} & 
\includegraphics[width=1.15cm]{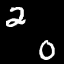} & 
\includegraphics[width=1.15cm]{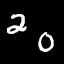} & 
\includegraphics[width=1.15cm]{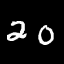} \\

\tiny{\textbf{AAO}} & 
\includegraphics[width=1.15cm]{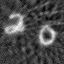} & 
\includegraphics[width=1.15cm]{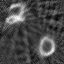} & 
\includegraphics[width=1.15cm]{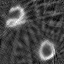} & 
\includegraphics[width=1.15cm]{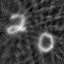} & 
\includegraphics[width=1.15cm]{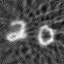} \\

\tiny{\textbf{AAO-ST}} & 
\includegraphics[width=1.15cm]{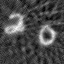} & 
\includegraphics[width=1.15cm]{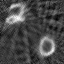} & 
\includegraphics[width=1.15cm]{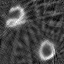} & 
\includegraphics[width=1.15cm]{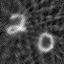} & 
\includegraphics[width=1.15cm]{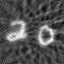} \\

\tiny{\textbf{AAO-OF}} & 
\includegraphics[width=1.15cm]{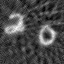} & 
\includegraphics[width=1.15cm]{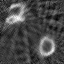} & 
\includegraphics[width=1.15cm]{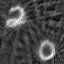} & 
\includegraphics[width=1.15cm]{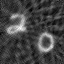} & 
\includegraphics[width=1.15cm]{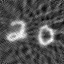} \\

\tiny{\textbf{IRKFS}} & 
\includegraphics[width=1.15cm]{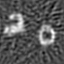} & 
\includegraphics[width=1.15cm]{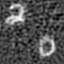} & 
\includegraphics[width=1.15cm]{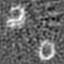} & 
\includegraphics[width=1.15cm]{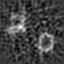} & 
\includegraphics[width=1.15cm]{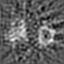} \\

\tiny{\textbf{IRKFS-M1}} & 
\includegraphics[width=1.15cm]{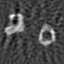} & 
\includegraphics[width=1.15cm]{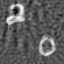} & 
\includegraphics[width=1.15cm]{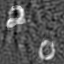} & 
\includegraphics[width=1.15cm]{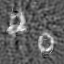} & 
\includegraphics[width=1.15cm]{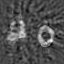} \\

\tiny{\textbf{IRKFS-M2}} & 
\includegraphics[width=1.15cm]{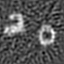} & 
\includegraphics[width=1.15cm]{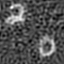} & 
\includegraphics[width=1.15cm]{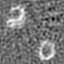} & 
\includegraphics[width=1.15cm]{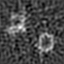} & 
\includegraphics[width=1.15cm]{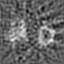} \\

\tiny{\textbf{IRKFS-M3}} & 
\includegraphics[width=1.15cm]{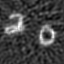} & 
\includegraphics[width=1.15cm]{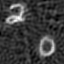} & 
\includegraphics[width=1.15cm]{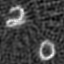} & 
\includegraphics[width=1.15cm]{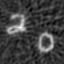} & 
\includegraphics[width=1.15cm]{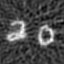} \\

\tiny{\textbf{EMIRKFS}} & 
\includegraphics[width=1.15cm]{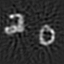} & 
\includegraphics[width=1.15cm]{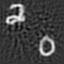} & 
\includegraphics[width=1.15cm]{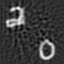} & 
\includegraphics[width=1.15cm]{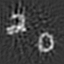} & 
\includegraphics[width=1.15cm]{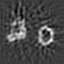} \\

\tiny{\textbf{EMIRKFS-M1}} & 
\includegraphics[width=1.15cm]{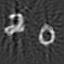} & 
\includegraphics[width=1.15cm]{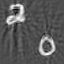} & 
\includegraphics[width=1.15cm]{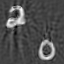} & 
\includegraphics[width=1.15cm]{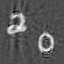} & 
\includegraphics[width=1.15cm]{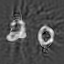} \\

\tiny{\textbf{EMIRKFS-M2}} & 
\includegraphics[width=1.15cm]{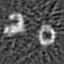} & 
\includegraphics[width=1.15cm]{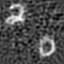} & 
\includegraphics[width=1.15cm]{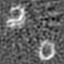} & 
\includegraphics[width=1.15cm]{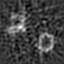} & 
\includegraphics[width=1.15cm]{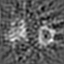} \\

\tiny{\textbf{EMIRKFS-M3}} & 
\includegraphics[width=1.15cm]{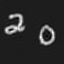} & 
\includegraphics[width=1.15cm]{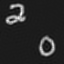} & 
\includegraphics[width=1.15cm]{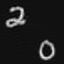} & 
\includegraphics[width=1.15cm]{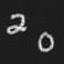} & 
\includegraphics[width=1.15cm]{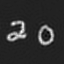} \\

\end{tabular}
\caption{Experiment 2: Truth image (Row 1) comparison to the image reconstruction of methods: AAO, AAO-ST, AAO-OF, IRKFS, IRKFS-M1, IRKFS-M2, IRFKS-M3, EMIRFKS, EMIRFKS-M1, EMIRKFS-M2, and EMIRFKS-M3 (Rows 2-12 from top to bottom) at time-steps $i\in\{2,5,7,9,11\}$ (Columns 1-5 from left to right).}
\label{fig:method_comparison4}
\end{figure}

Again, we compare the reconstruction quality and computational efficiency of the proposed EMIRKFS-M method described in Algorithm \ref{Alg: EMIRKFS-M} to other alternative benchmark approaches of (\ref{algorithms}). We plot negative image of the absolute value of the error images across the same methods and time-steps. Additionally, we consider the accuracy refinements across all iterative methods.
\paragraph{Results and Discussion} The Figures of \ref{fig:method_comparison4} demonstrate the dramatic improvements in the reduction of background artifacts given by the EMIRKFS-M3 method. Reconstructions given via EMIRKFS and EMIRKFS-M1 show slight improvements in background artifacts, while the rest of the methods fail to do so. When considering both visual quality and relative error, the EMIRKFS-M3 method is the overall winner. 
In Figure \ref{fig:method_comparison6} we see that the IRKFS-M3, EMIRKFS, EMIRKFS-M1, and EMIRKFS-M3 methods decrease in relative error over 5 iterations. However, the EMIRKFS-M3 method demonstrates the most significant improvement.

\begin{figure}[ht]
\centering
\begin{tabular}{lccccc}
\textbf{Method} & \textbf{i=2} & \textbf{i=5} & \textbf{i=7} & \textbf{i=9} & \textbf{i=11} \\
\hline

\tiny{\textbf{AAO}} & 
\includegraphics[width=1.15cm]{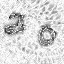} & 
\includegraphics[width=1.15cm]{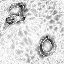} & 
\includegraphics[width=1.15cm]{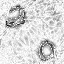} & 
\includegraphics[width=1.15cm]{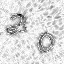} & 
\includegraphics[width=1.15cm]{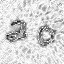} \\

\tiny{\textbf{AAO-ST}} & 
\includegraphics[width=1.15cm]{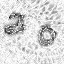} & 
\includegraphics[width=1.15cm]{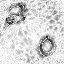} & 
\includegraphics[width=1.15cm]{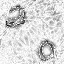} & 
\includegraphics[width=1.15cm]{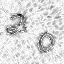} & 
\includegraphics[width=1.15cm]{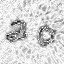} \\

\tiny{\textbf{AAO-OF}} & 
\includegraphics[width=1.15cm]{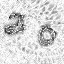} & 
\includegraphics[width=1.15cm]{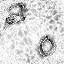} & 
\includegraphics[width=1.15cm]{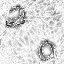} & 
\includegraphics[width=1.15cm]{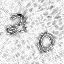} & 
\includegraphics[width=1.15cm]{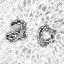} \\

\tiny{\textbf{IRKFS}} & 
\includegraphics[width=1.15cm]{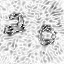} & 
\includegraphics[width=1.15cm]{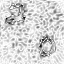} & 
\includegraphics[width=1.15cm]{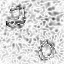} & 
\includegraphics[width=1.15cm]{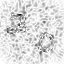} & 
\includegraphics[width=1.15cm]{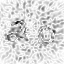} \\

\tiny{\textbf{IRKFS-M1}} & 
\includegraphics[width=1.15cm]{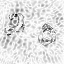} & 
\includegraphics[width=1.15cm]{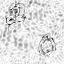} & 
\includegraphics[width=1.15cm]{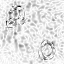} & 
\includegraphics[width=1.15cm]{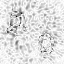} & 
\includegraphics[width=1.15cm]{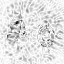} \\

\tiny{\textbf{IRKFS-M2}} & 
\includegraphics[width=1.15cm]{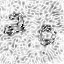} & 
\includegraphics[width=1.15cm]{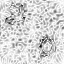} & 
\includegraphics[width=1.15cm]{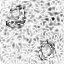} & 
\includegraphics[width=1.15cm]{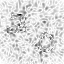} & 
\includegraphics[width=1.15cm]{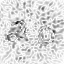} \\

\tiny{\textbf{IRKFS-M3}} & 
\includegraphics[width=1.15cm]{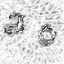} & 
\includegraphics[width=1.15cm]{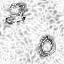} & 
\includegraphics[width=1.15cm]{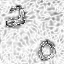} & 
\includegraphics[width=1.15cm]{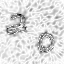} & 
\includegraphics[width=1.15cm]{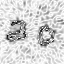} \\

\tiny{\textbf{EMIRKFS}} & 
\includegraphics[width=1.15cm]{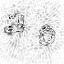} & 
\includegraphics[width=1.15cm]{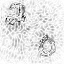} & 
\includegraphics[width=1.15cm]{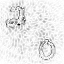} & 
\includegraphics[width=1.15cm]{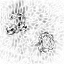} & 
\includegraphics[width=1.15cm]{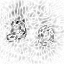} \\

\tiny{\textbf{EMIRKFS-M1}} & 
\includegraphics[width=1.15cm]{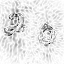} & 
\includegraphics[width=1.15cm]{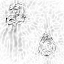} & 
\includegraphics[width=1.15cm]{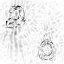} & 
\includegraphics[width=1.15cm]{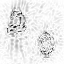} & 
\includegraphics[width=1.15cm]{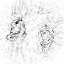} \\

\tiny{\textbf{EMIRKFS-M2}} & 
\includegraphics[width=1.15cm]{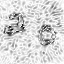} & 
\includegraphics[width=1.15cm]{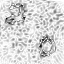} & 
\includegraphics[width=1.15cm]{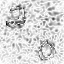} & 
\includegraphics[width=1.15cm]{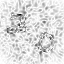} & 
\includegraphics[width=1.15cm]{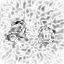} \\

\tiny{\textbf{EMIRKFS-M3}} & 
\includegraphics[width=1.15cm]{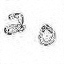} & 
\includegraphics[width=1.15cm]{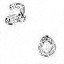} & 
\includegraphics[width=1.15cm]{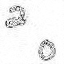} & 
\includegraphics[width=1.15cm]{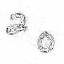} & 
\includegraphics[width=1.15cm]{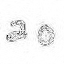}
\end{tabular}
\caption{Experiment 2: Reverse grayscale error image reconstructions of methods: AAO, AAO-ST, AAO-OF, IRKFS, IRKFS-M1, IRKFS-M2, IRFKS-M3, EMIRFKS, EMIRFKS-M1, EMIRKFS-M2, and EMIRFKS-M3 (Rows 1-11 from top to bottom) at time-steps $i\in\{2,5,7,9,11\}$ (Columns 1-5 from left to right).}
\label{fig:method_comparison5}
\end{figure}

\begin{figure}[ht]
\centering
\begin{tabular}{cc} 
\includegraphics[width=5.5cm]{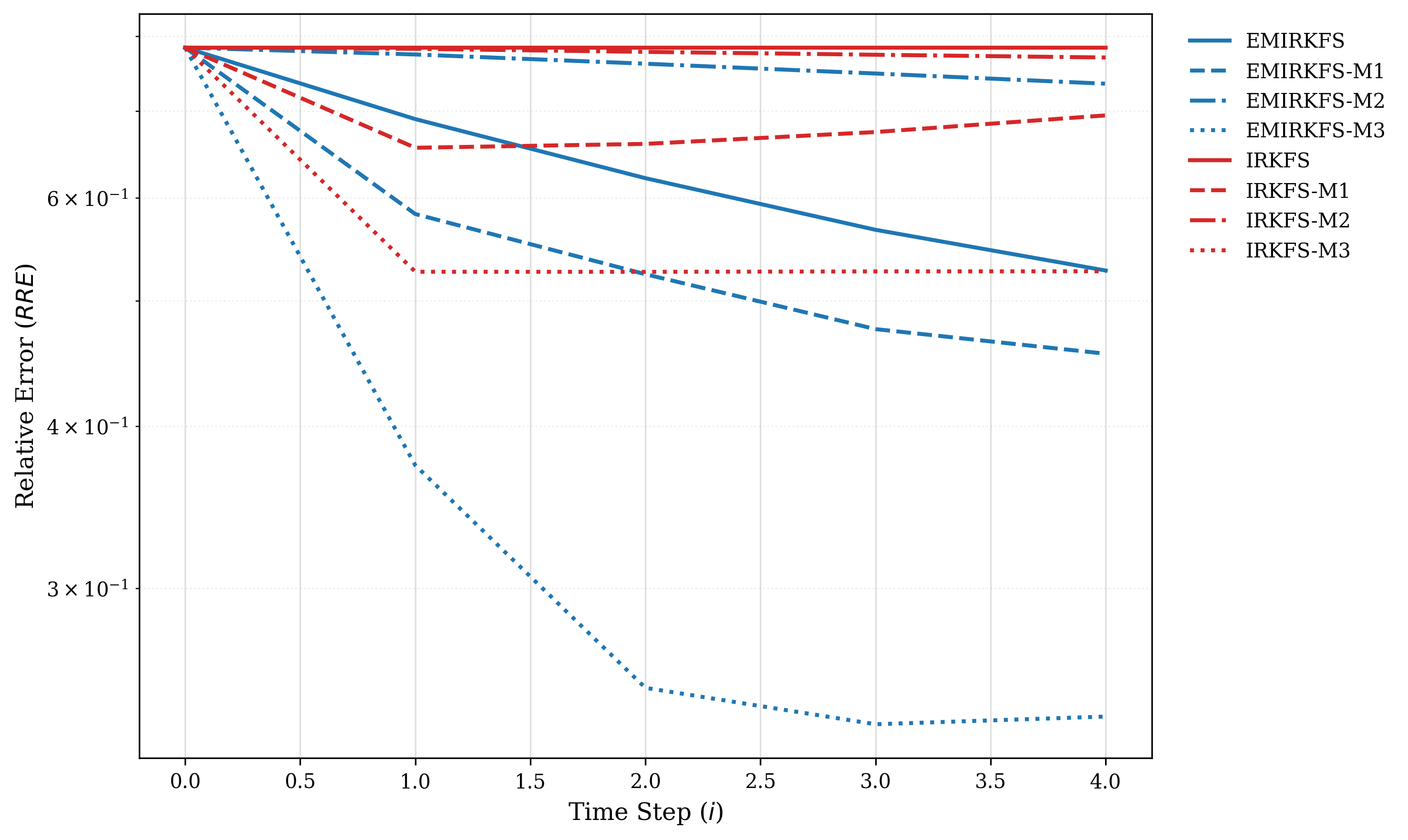} & 
\includegraphics[width=5.5cm]{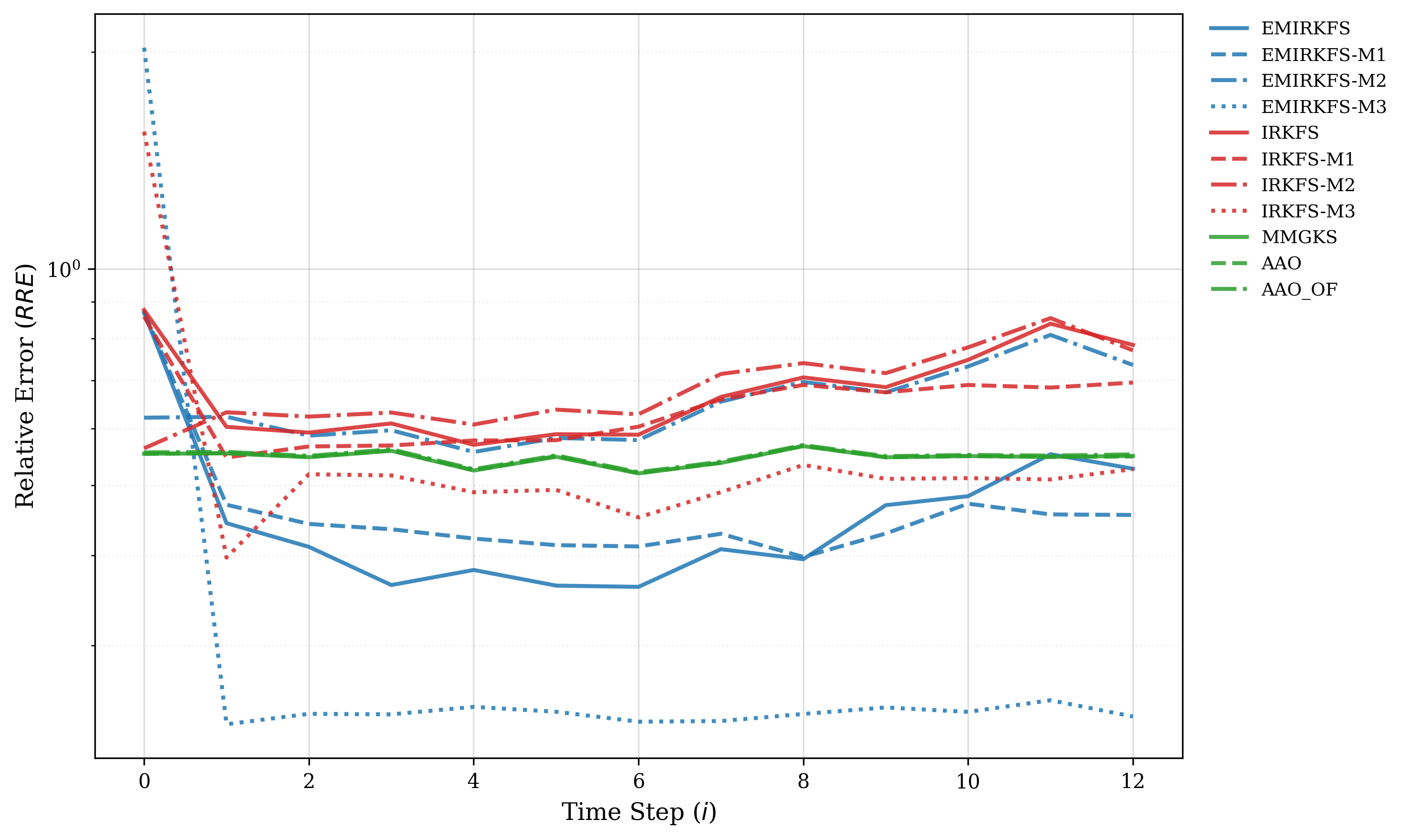} \\ 

\end{tabular}
\caption{Experiment 2: Memory vs. Time (left) and RRE (\ref{RRE}) (right) comparison across all methods (\ref{algorithms})}
\end{figure}\label{fig:method_comparison6}
\subsection{Experiment 3: Emoji Dataset}\label{ssec:exp_phantom_3}
We now proceed to evaluate the implementation of algorithm (\ref{Alg: EMIRKFS-M}) using the Dynamic Emoji dataset. This dataset considers real data of an “emoji” phantom measured at the University of Helsinki \cite{pasha2023computational}. Here the data is modified as in \cite{okunola2025efficient} to determine a limited angle problem by limiting the number of projection angles to 10. The resulting phantom is a
3D array with dimensions of $128\times 128\times 33$.
\paragraph{Forward Model and Projection Geometry}
The forward operator is constructed using the \texttt{Trips-Py} toolbox \cite{pasha2024trips}, employing 10 projection angles per image. Each individual system matrix is defined as $\mathbf{H}_i$ has dimension $2170\times 16384$ with observations $\mathbf{y}_i\in\mathbb{R}^{2170}$, and the complete forward model is represented by the global matrix $\mathbf{H}$ with dimension $71610\times 540672$ and corresponding vectorized sinograms from all time-steps $\mathbf{y}\in\mathbb{R}^{71610}$. The respective sinogram data is corrupted with additional additive noise $\sigma_{\text{NL}}=0.01$. For both M2 and M3 motion models we fix $\zeta=3$, and for M3 we select a patchsize of $z_x=z_y=2$.  
We assess the reconstruction performance of the of the iterative methods with parameters fixed as follows: $r=1000$, $\alpha = 0.01$,  $\ell =4.8$, and $n_{\text{iter}}=3$. 
We initialize with $\boldsymbol{\Psi}^{\text{est}}_0=\mathbf{I}_r$, $\mathbf{Q}^{(0)}_i=\alpha^2\mathbf{I}_{n_s}$, $\mathbf{M}^{(0)}_i=\mathbf{I}_{n_s}$, $\mathbf{R}^{(0)}_i=\alpha^2\mathbf{I}_{m_t}$,  and as in \cite{hakkarainen2019undersampled} we set $\mathbf{x}^{\text{est}}_0=\mathbf{P}_r((\mathbf{H}_0\mathbf{P}_r)^\top\mathbf{H}_0\mathbf{P}_r+\frac{1}{\alpha^2}\mathbf{P}_r^\top\mathbf{P}_r)^{-1}\mathbf{y}_0$.

\begin{figure}[ht]
\centering
\begin{tabular}{lccccccc}
\footnotesize{\textbf{Method}} & \footnotesize{\textbf{i=3}} & \footnotesize{\textbf{i=7}} & \footnotesize{\textbf{i=11}} & \footnotesize{\textbf{i=15}} & \footnotesize{\textbf{i=19}} & \footnotesize{\textbf{i=22}} & \footnotesize{\textbf{i=29}} \\
\hline
\tiny{\textbf{AAO}} & 
\includegraphics[width=1.15cm]{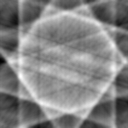} & 
\includegraphics[width=1.15cm]{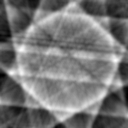} & 
\includegraphics[width=1.15cm]{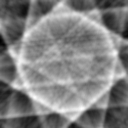} & 
\includegraphics[width=1.15cm]{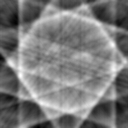} & 
\includegraphics[width=1.15cm]{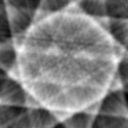} & 
\includegraphics[width=1.15cm]{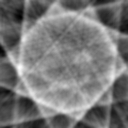} & 
\includegraphics[width=1.15cm]{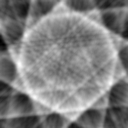} \\

\tiny{\textbf{AAO-ST}} & 
\includegraphics[width=1.15cm]{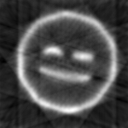} & 
\includegraphics[width=1.15cm]{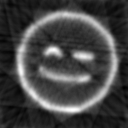} & 
\includegraphics[width=1.15cm]{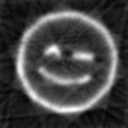} & 
\includegraphics[width=1.15cm]{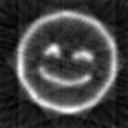} & 
\includegraphics[width=1.15cm]{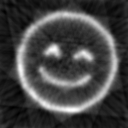} & 
\includegraphics[width=1.15cm]{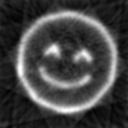} & 
\includegraphics[width=1.15cm]{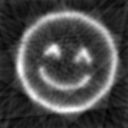} \\

\tiny{\textbf{AAO-OF}} & 
\includegraphics[width=1.15cm]{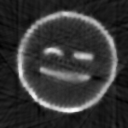} & 
\includegraphics[width=1.15cm]{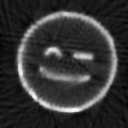} & 
\includegraphics[width=1.15cm]{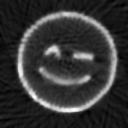} & 
\includegraphics[width=1.15cm]{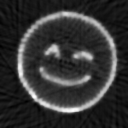} & 
\includegraphics[width=1.15cm]{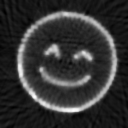} & 
\includegraphics[width=1.15cm]{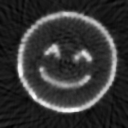} & 
\includegraphics[width=1.15cm]{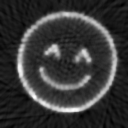} \\
\tiny{\textbf{IRKFS}} & 
\includegraphics[width=1.15cm]{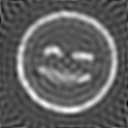} & 
\includegraphics[width=1.15cm]{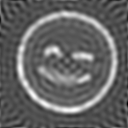} & 
\includegraphics[width=1.15cm]{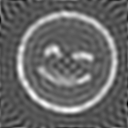} & 
\includegraphics[width=1.15cm]{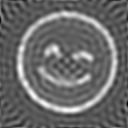} & 
\includegraphics[width=1.15cm]{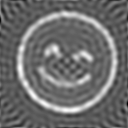} & 
\includegraphics[width=1.15cm]{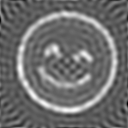} & 
\includegraphics[width=1.15cm]{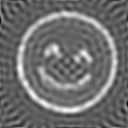} \\

\tiny{\textbf{IRKFS-M1}} & 
\includegraphics[width=1.15cm]{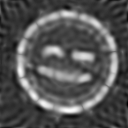} & 
\includegraphics[width=1.15cm]{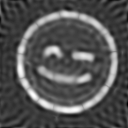} & 
\includegraphics[width=1.15cm]{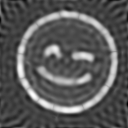} & 
\includegraphics[width=1.15cm]{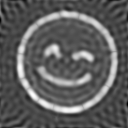} & 
\includegraphics[width=1.15cm]{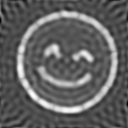} & 
\includegraphics[width=1.15cm]{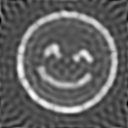} & 
\includegraphics[width=1.15cm]{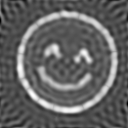} \\

\tiny{\textbf{IRKFS-M2}} & 
\includegraphics[width=1.15cm]{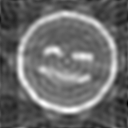} & 
\includegraphics[width=1.15cm]{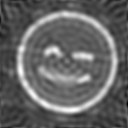} & 
\includegraphics[width=1.15cm]{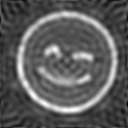} & 
\includegraphics[width=1.15cm]{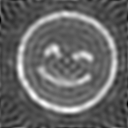} & 
\includegraphics[width=1.15cm]{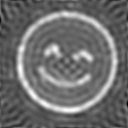} & 
\includegraphics[width=1.15cm]{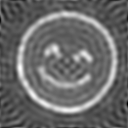} & 
\includegraphics[width=1.15cm]{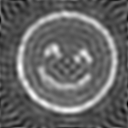} \\

\tiny{\textbf{IRKFS-M3}} & 
\includegraphics[width=1.15cm]{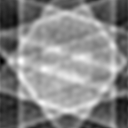} & 
\includegraphics[width=1.15cm]{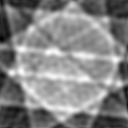} & 
\includegraphics[width=1.15cm]{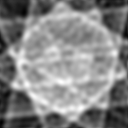} & 
\includegraphics[width=1.15cm]{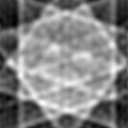} & 
\includegraphics[width=1.15cm]{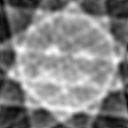} & 
\includegraphics[width=1.15cm]{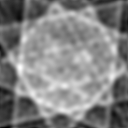} & 
\includegraphics[width=1.15cm]{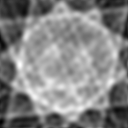} \\

\tiny{\textbf{EMIRKFS}} & 
\includegraphics[width=1.15cm]{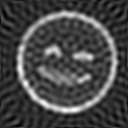} & 
\includegraphics[width=1.15cm]{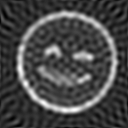} & 
\includegraphics[width=1.15cm]{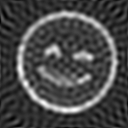} & 
\includegraphics[width=1.15cm]{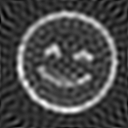} & 
\includegraphics[width=1.15cm]{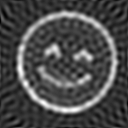} & 
\includegraphics[width=1.15cm]{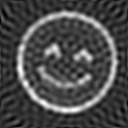} & 
\includegraphics[width=1.15cm]{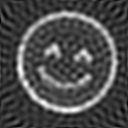} \\

\tiny{\textbf{EMIRKFS-M1}} & 
\includegraphics[width=1.15cm]{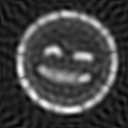} & 
\includegraphics[width=1.15cm]{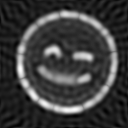} & 
\includegraphics[width=1.15cm]{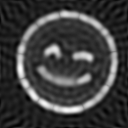} & 
\includegraphics[width=1.15cm]{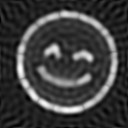} & 
\includegraphics[width=1.15cm]{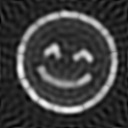} & 
\includegraphics[width=1.15cm]{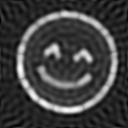} & 
\includegraphics[width=1.15cm]{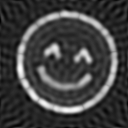} \\

\tiny{\textbf{EMIRKFS-M2}} & 
\includegraphics[width=1.15cm]{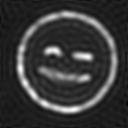} & 
\includegraphics[width=1.15cm]{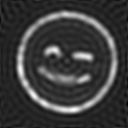} & 
\includegraphics[width=1.15cm]{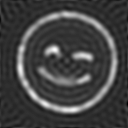} & 
\includegraphics[width=1.15cm]{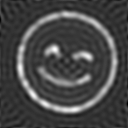} & 
\includegraphics[width=1.15cm]{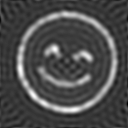} & 
\includegraphics[width=1.15cm]{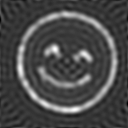} & 
\includegraphics[width=1.15cm]{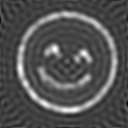} \\

\tiny{\textbf{EMIRKFS-M3}} & 
\includegraphics[width=1.15cm]{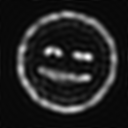} & 
\includegraphics[width=1.15cm]{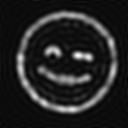} & 
\includegraphics[width=1.15cm]{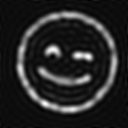} & 
\includegraphics[width=1.15cm]{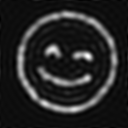} & 
\includegraphics[width=1.15cm]{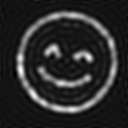} & 
\includegraphics[width=1.15cm]{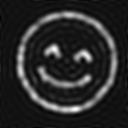} & 
\includegraphics[width=1.15cm]{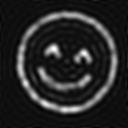}
\end{tabular}
\caption{Experiment 3: Truth image (Row 1) comparison to the image reconstruction of methods: AAO, AAO-ST, AAO-OF, IRKFS, IRKFS-M1, IRKFS-M2, IRFKS-M3, EMIRFKS, EMIRFKS-M1, EMIRKFS-M2, and EMIRFKS-M3 (Rows 2-12 from top to bottom) at time-steps $i\in\{3,7,11,15,19,22,29\}$ (Columns 1-7 from left to right).}
\label{fig:method_comparison7}
\end{figure}

\begin{figure}[ht]
\centering
\includegraphics[width=10.5cm]{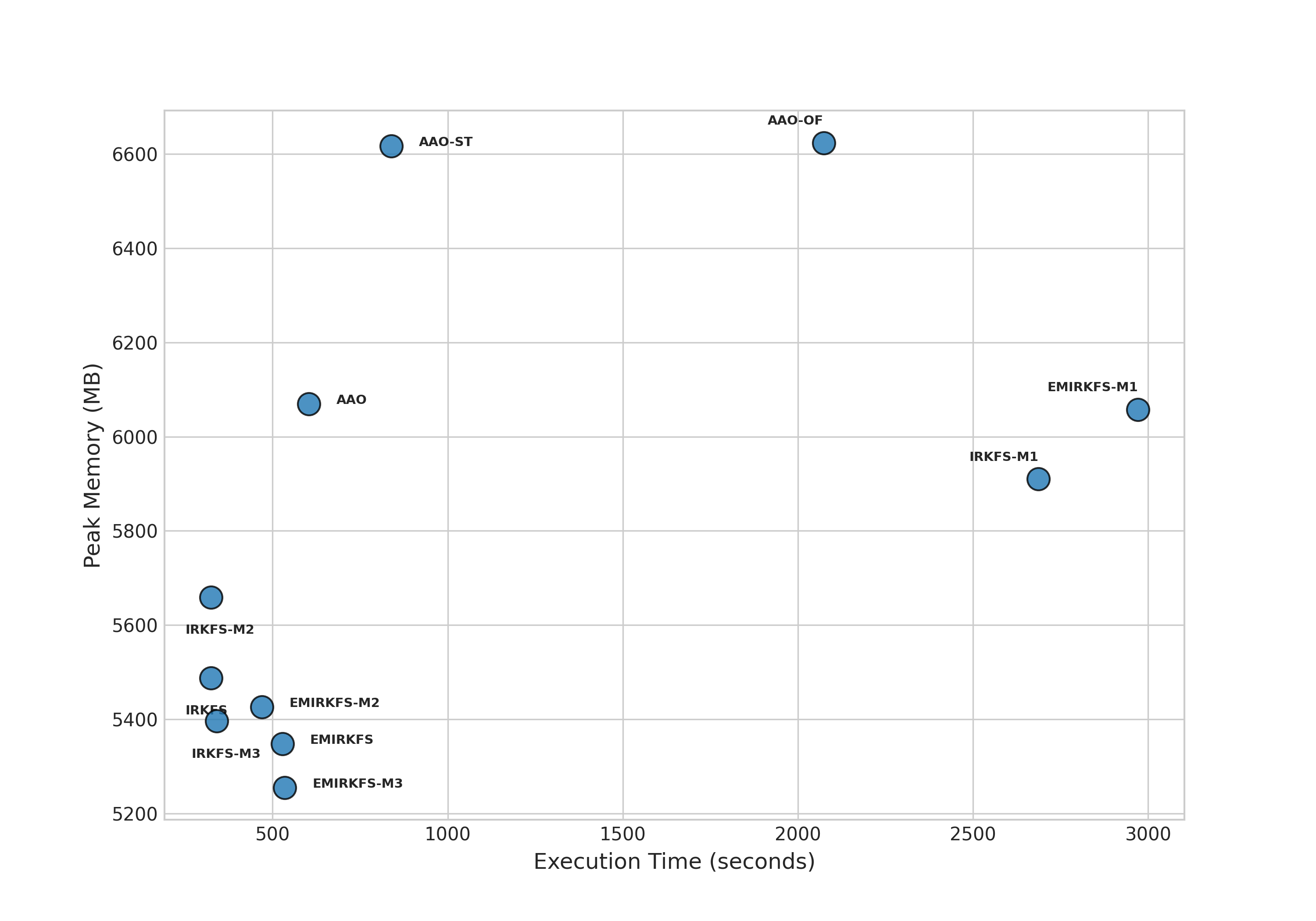}
\caption{Experiment 3: Memory vs. Time comparison across all methods (\ref{algorithms})}
\label{fig:method_comparison8}
\end{figure}

\paragraph{Results and Discussion} From Figure \ref{fig:method_comparison7}, it is evident that most methods demonstrate sharp reconstructions besides that of AAO and IRKFS-M3 methods. The highest quality reconstructions are given via 
AAO-OF and EMIRKFS-M3 methods, however, from Figure \ref{fig:method_comparison8}, the EMIRKFS-M3 method is again superior in both memory consumption and computational time. When considering both visual quality and these computational aspects, either EMIRKFS-M1 or EMIRKFS-M3 methods are superior. 
\section{Conclusion and Outlook}\label{sec: Conclusions}
In this work, we introduced \emph{EMIRKFS-M}, an iterative, dimension-reduced Kalman filtering and smoothing framework with automatic parameter selection. The proposed method builds on the reduced Kalman filtering and smoothing (RKFS) approach for dynamic inverse problems \cite{hakkarainen2019undersampled} by integrating expectation–maximization and motion estimation techniques. This integration enables the covariance and motion matrices of linear Gaussian systems to be updated automatically and iteratively, while preserving the computational efficiency that makes RKFS attractive for large-scale problems.

Within this framework, we considered three different approaches for estimating the motion model. One approach is based on optical flow, which provides accurate motion estimation but is computationally expensive. The remaining two approaches are computationally cheaper alternatives, designed to reduce overall cost while still capturing the dominant motion dynamics. This comparison highlights the trade-offs between accuracy and efficiency when incorporating motion information into the filtering and smoothing process.

To evaluate the performance of the proposed approach, we systematically examined the contributions of the individual motion models described in Subsections \ref{sub: M1}, \ref{sub: M2}, and \ref{sub: M3} and the expectation–maximization procedure described in Section \ref{sec: Proposed} using three computed tomography test problems with both low- and high-angle projection geometries. In addition, we compared EMIRKFS-M (Algorithm \ref{Alg: EMIRKFS-M}) against several all-at-once reconstruction methods, as well as the original RKFS framework \cite{hakkarainen2019undersampled}. The comparison included qualitative reconstructions, relative reconstruction error (RRE) metrics, and computational cost analyses. Across all experiments, EMIRKFS-M consistently produced sharper reconstructions and lower relative errors, while retaining the favorable memory usage and runtime performance of RKFS relative to AAO variants.

Several directions for future work remain. A key extension is the formulation of the underlying linear Gaussian system within a more general non-Gaussian prior setting, for example through the use of ensemble Kalman filtering with Besov priors. Additionally, incorporating time-dependent state covariance matrices $\boldsymbol{\Sigma}_i$ would further enhance the flexibility of EMIRKFS-M, enabling it to better capture non-stationary dynamics in time-dependent inverse problems.

\section*{Acknowledgments}
AK and MP acknowledge support from NSF DMS 2410699. Any opinions, findings, conclusions, or recommendations expressed in this material are those of the authors and do not necessarily reflect the views of the National Science Foundation. MP would like to thank Youssef Marzouk for discussions about the sequential setting for dynamic inverse problems and Ricardo Baptista for the discussions on expectation maximization that lead to improvements of this work. The authors acknowledge Advanced Research Computing at Virginia Tech for providing computational resources and technical support that have contributed to the results reported within this paper. URL: https://arc.vt.edu/.
\newpage
\appendix
\vspace*{1em}
\renewcommand{\thesection}{A.\arabic{section}}
\titleformat{\section}{\normalfont\large\bfseries}{\thesection}{1em}{} 
\titlespacing*{\section}{0pt}{1ex plus .2ex minus .2ex}{0.5ex plus .1ex}
\section{Dimension-Reduced Kalman Filter}\label{sub: App_A}
\vspace{-0.5em}

\begin{algorithm}[H]
\caption{Dimension-Reduced Kalman Filter (RKF)}
\label{Alg:RKF}
\footnotesize
{\begin{algorithmic}[1]
\Require Covariance parameters $\alpha$, $\ell$; reduction dimension $r$; initial estimates $\mathbf{x}_0^{{\text{\text{est}}}}$, $\bm{\Psi}_0^{{\text{\text{est}}}}$; motion operators $\{\mathbf{M}_i\}_{i=1}^{T}$; noise covariances $\{\mathbf{Q}_i\}_{i=1}^{T}$, $\{\mathbf{R}_i\}_{i=1}^{T}$; observations $\{\mathbf{y}_i\}_{i=1}^{T}$
\Ensure Estimated states $\{\mathbf{x}_i^{{\text{\text{est}}}}\}_{i=0}^{T}$, covariances $\{\bm{\Psi}_i^{{\text{\text{est}}}}\}_{i=0}^{T}$, computation matrices $\{\mathbf{B}_i\}_{i=1}^{T}$

\vspace{1mm}
\State \textbf{Initialization:}
\Statex \hspace{\algorithmicindent} Compute covariance matrix: $\mathbf{\Sigma}_{ij} = \alpha^2 \exp\left(-\frac{d(x_i,x_j)^2}{2\ell^2}\right)$
\Statex \hspace{\algorithmicindent} Perform SVD: $\mathbf{\Sigma} = \mathbf{U} \mathbf{S} \mathbf{U}^\top$
\Statex \hspace{\algorithmicindent} Construct projection matrix:
\[
\mathbf{U}_r = \mathbf{U}(:,1\!:\!r), \quad \mathbf{S}_r = \mathbf{S}(1\!:\!r,1\!:\!r), \quad \mathbf{P}_r = \mathbf{U}_r \mathbf{S}_r^{1/2}
\]

\For{$i = 1$ to $T$}
    \State Predict state: $\mathbf{x}_i^p = \mathbf{M}_i \mathbf{x}_{i-1}^{{\text{\text{est}}}}$
    \State Decompose covariance: $\bm{\Psi}_{i-1}^{{\text{\text{est}}}} = \mathbf{A}_i \mathbf{A}_i^\top$
    \State Compute: $\mathbf{B}_i = \mathbf{M}_i \mathbf{P}_r \mathbf{A}_i$
    \State Compute inverse covariance:
    \[
    (\mathbf{C}_i^p)^{-1} \mathbf{P}_r = \mathbf{Q}_i^{-1} \mathbf{P}_r - \mathbf{Q}_i^{-1} \mathbf{B}_i 
    \left(\mathbf{B}_i^\top \mathbf{Q}_i^{-1} \mathbf{B}_i + \mathbf{I}\right)^{-1}
    \mathbf{B}_i^\top \mathbf{Q}_i^{-1} \mathbf{P}_r
    \]
    \State Compute posterior covariance:
    \[
    \bm{\Psi}_i^{{\text{\text{est}}}} = \left( (\mathbf{H}_i \mathbf{P}_r)^\top \mathbf{R}_i^{-1} (\mathbf{H}_i \mathbf{P}_r) + \mathbf{P}_r^\top (\mathbf{C}_i^p)^{-1} \mathbf{P}_r \right)^{-1}
    \]
    \State Estimate reduced state:
    \[
    \bm{\alpha}_i^{{\text{est}}} = \bm{\Psi}_i^{{\text{est}}} (\mathbf{H}_i \mathbf{P}_r)^\top \mathbf{R}_i^{-1} (\mathbf{y}_i - \mathbf{H}_i \mathbf{x}_i^p)
    \]
    \State Update full state estimate: $\mathbf{x}_i^{{\text{est}}} = \mathbf{x}_i^p + \mathbf{P}_r \bm{\alpha}_i^{{\text{est}}}$
\EndFor
\end{algorithmic}}
\end{algorithm}

\section{Dimension-Reduced Kalman Smoother}\label{sub: App_B}
\vspace{-0.5em}
\begin{algorithm}[H]

\caption{Dimension-Reduced Kalman Smoother (RKS)}
\label{Alg:RKS}
\footnotesize
\begin{algorithmic}[1]
\Require Projection matrix $\mathbf{P}_r$;  state estimates $\{\mathbf{x}_i^{ {\text{est}}}\}_{i=0}^{T}$; covariance estimates $\{\bm{\Psi}_i^{ {\text{est}}}\}_{i=0}^{T-1}$; motion operators $\{\mathbf{M}_i\}_{i=1}^{T}$; noise covariances $\{\mathbf{Q}_i\}_{i=1}^{T}$; auxiliary matrices $\{\mathbf{B}_i\}_{i=1}^{T}$
\Ensure Smoothed state estimates $\{\mathbf{x}_i^{ {\text{sm}}}\}_{i=0}^{T}$; smoothed covariance estimates $\{\boldsymbol{\Psi}_i^{ {\text{sm}}}\}_{i=0}^{T} $

\State \textbf{Initialization:} Set $\mathbf{x}_T^{ {\text{sm}}} = \mathbf{x}_T^{ {\text{est}}}$; (Optional)  Set $\boldsymbol{\Psi}_T^{ {\text{sm}}} = \boldsymbol{\Psi}_T^{ {\text{est}}}$

\For{$i = T$ down to $1$}
    \State Compute covariances:
    \[
    \mathbf{C}_{i-1}^{ {\text{est}}} = \mathbf{P}_r\, \bm{\Psi}_{i-1}^{ {\text{est}}}\, \mathbf{P}_r^\top, \quad
    \mathbf{C}_{i}^{p} = \mathbf{B}_i\, \mathbf{B}_i^\top + \mathbf{Q}_i
    \]
    \State Update smoothed state:
    \[
    \mathbf{x}_{i-1}^{ {\text{sm}}} = \mathbf{x}_{i-1}^{ {\text{est}}} + \mathbf{C}_{i-1}^{ {\text{est}}}\, \mathbf{M}_i^\top\, (\mathbf{C}_{i}^{p})^{-1} \left(\mathbf{x}_i^{ {\text{sm}}} - \mathbf{M}_i\, \mathbf{x}_{i-1}^{ {\text{est}}}\right)
    \]
    \State (Optional) Update smoothed covariance $\boldsymbol{\Psi}_{i-1}^{\text{sm}}$ via (\ref{eq: reduced_cov})
\EndFor

\end{algorithmic}
\end{algorithm}

\section{Linear Update Operator Solver}\label{sub: App_C}
\vspace{-0.5em}
\begin{algorithm}[H]
\footnotesize
\caption{Linear Update Operator Solver}
\label{alg:linear_update_operator}
\begin{algorithmic}[1]
\Require Current and previous state vectors $\mathbf{x}_{i-1}, \mathbf{x}_i$
\Ensure Approximate linear update operator $\mathbf{M}_i$ such that $\mathbf{x}_i \approx \mathbf{M}_i \mathbf{x}_{i-1}$
\Function{UpdateOperator}{$\mathbf{x}_{i-1}, \mathbf{x}_i$}
    \State Compute $\mathbf{V}(\mathbf{x}_{i-1})$ and transformation $\mathbf{T}(t_{i-1})$
    \State Construct regularization matrix:
    \[
    \mathbf{L} = \begin{bmatrix}
        \mathbf{I}_{n_x} \otimes \mathbf{L}_y \\
        \mathbf{L}_x \otimes \mathbf{I}_{n_y}
    \end{bmatrix}
    \]
    \State Solve inverse problem using MMGKS:
    \[
    \mathbf{s}(t_{i-1}) = \texttt{MMGKS}(\mathbf{V}(\mathbf{x}_{i-1}), [\mathbf{L} \quad \mathbf{L}], -\mathbf{T}(t_{i-1}))
    \]
    \State Construct update operator: $\mathbf{M}_i = \mathbf{M}(-\mathbf{s}(t_{i-1}))$
    \State \Return $\mathbf{M}_i$
\EndFunction
\end{algorithmic}
\end{algorithm}

\newpage
\bibliographystyle{plain}
\bibliography{Arxiv}

@article{buccini2025krylov,
  title={Krylov Subspace Based FISTA-Type Methods for Linear Discrete Ill-Posed Problems},
  author={Buccini, Alessandro and Chen, Fei and Pasha, Mirjeta and Reichel, Lothar},
  journal={Numerical Linear Algebra with Applications},
  volume={32},
  number={1},
  pages={e2610},
  year={2025},
  publisher={Wiley Online Library}
}

@article{lindbloom2025priorconditioned,
  title={Priorconditioned sparsity-promoting projection methods for deterministic and {B}ayesian linear inverse problems},
  author={Lindbloom, Jonathan and Pasha, Mirjeta and Glaubitz, Jan and Marzouk, Youssef},
  journal={arXiv preprint arXiv:2505.01827},
  year={2025}
}

@article{shumway1982approach,
  title={An approach to time series smoothing and forecasting using the EM algorithm},
  author={Shumway, Robert H and Stoffer, David S},
  journal={Journal of time series analysis},
  volume={3},
  number={4},
  pages={253--264},
  year={1982},
  publisher={Wiley Online Library}
}

@article{buccini2020modulus,
  title={Modulus-based iterative methods for constrained $\ell_p-\ell_q$ minimization},
  author={Buccini, Alessandro and Pasha, Mirjeta and Reichel, Lothar},
  journal={Inverse Problems},
  volume={36},
  number={8},
  pages={084001},
  year={2020},
  publisher={IOP Publishing}
}

@inproceedings{cite12,
  author    = {Bubba, T. A. and M{\"a}rz, M. and Purisha, Z. and Lassas, M. and Siltanen, S.},
  title     = {Shearlet-based regularization in sparse dynamic tomography},
  booktitle = {Proceedings of Wavelets and Sparsity XVII},
  series    = {Proc. SPIE},
  volume    = {10394},
  pages     = {103940Y},
  year      = {2017},
  doi       = {10.1117/12.2274485}
}

@article{cite13,
  author    = {Burger, M. and Dirks, H. and Frerking, L. and Hauptmann, A. and Helin, T. and Siltanen, S.},
  title     = {A variational reconstruction method for undersampled dynamic {X}-ray tomography based on physical motion models},
  journal   = {Inverse Problems},
  volume    = {33},
  number    = {12},
  pages     = {124008},
  year      = {2017},
  doi       = {10.1088/1361-6420/aa99cf}
}

@article{cite14,
  author    = {Niemi, E. and Lassas, M. and Kallonen, A. and Harhanen, L. and H{\"a}m{\"a}l{\"a}inen, K. and Siltanen, S.},
  title     = {Dynamic multisource {X}-ray tomography using a spacetime level set method},
  journal   = {J. Comput. Phys.},
  volume    = {291},
  pages     = {218--237},
  year      = {2015},
  doi       = {10.1016/j.jcp.2015.03.010}
}

@article{emkal1,
  author    = {Khan, Mohammad Emtiyaz and Dutt, Deshpande Narayana},
  title     = {{An Expectation-Maximization Algorithm Based Kalman Smoother Approach for Event-Related Desynchronization (ERD) Estimation from EEG}},
  journal   = {IEEE Transactions on Biomedical Engineering},
  volume    = {54},
  number    = {7},
  pages     = {1191--1198},
  year      = {2007},
  doi       = {10.1109/TBME.2007.894827},
  issn      = {0018-9294},
  pmid      = {17605350}
}

@article{emkal2,
  author    = {Shi, Zhuangwei},
  title     = {{Incorporating Transformer and LSTM to Kalman Filter with EM algorithm for state estimation}},
  journal   = {Engineering Applications of Artificial Intelligence},
  volume    = {106},
  pages     = {104523},
  year      = {2021},
  doi       = {10.1016/j.engappai.2021.104523},
  publisher = {Elsevier}
}

@inproceedings{emkal3,
  author    = {Ghahramani, Zoubin and Roweis, Sam T.},
  title     = {Learning Nonlinear Dynamical Systems using an {EM} Algorithm},
  booktitle = {Advances in Neural Information Processing Systems 12 ({NIPS} 1999)},
  pages     = {431--437},
  year      = {2000},
  editor    = {S. A. Solla and T. K. Leen and K.-R. M{\"u}ller}
}

@misc{emoji,
      title={Tomographic {X}-ray data of {3D} emoji}, 
      author={Alexander Meaney and Zenith Purisha and Samuli Siltanen},
      year={2018},
      eprint={1802.09397},
      archivePrefix={arXiv},
      primaryClass={physics.med-ph},
      url={https://arxiv.org/abs/1802.09397}, 
}

@article{pasha2023recycling,
  title={Recycling MMGKS for large-scale dynamic and streaming data},
  author={Pasha, Mirjeta and de Sturler, Eric and Kilmer, Misha E},
  journal={arXiv preprint arXiv:2309.15759},
  year={2023}
}

@article{yang2020ground,
  title={Ground moving target tracking and refocusing using shadow in video-SAR},
  author={Yang, Xiaqing and Shi, Jun and Zhou, Yuanyuan and Wang, Chen and Hu, Yao and Zhang, Xiaoling and Wei, Shunjun},
  journal={Remote Sensing},
  volume={12},
  number={18},
  pages={3083},
  year={2020},
  publisher={MDPI}
}

@article{CT1,
  author    = {Ambrose, J. and Hounsfield, G.},
  title     = {Computerized transverse axial tomography},
  journal   = {British Journal of Radiology},
  year      = {1973},
  month     = {February},
  volume    = {46},
  number    = {542},
  pages     = {148--149},
  pmid      = {4686818}
}

@article{CT2,
  author    = {Hounsfield, G. N.},
  title     = {Computerized transverse axial scanning (tomography). 1. Description of system},
  journal   = {British Journal of Radiology},
  year      = {1973},
  month     = {December},
  volume    = {46},
  number    = {552},
  pages     = {1016--1022},
  doi       = {10.1259/0007-1285-46-552-1016},
  pmid      = {4757352}
}

@article{FBP1,
  author    = {Ramachandran, G. N. and Lakshminarayanan, A. V.},
  title     = {Three-dimensional reconstruction from radiographs and electron micrographs: application of convolutions instead of Fourier transforms},
  journal   = {Proceedings of the National Academy of Sciences of the United States of America},
  year      = {1971},
  month     = {September},
  volume    = {68},
  number    = {9},
  pages     = {2236--2240},
  doi       = {10.1073/pnas.68.9.2236},
  pmid      = {5289381},
  pmcid     = {PMC389392}
}

@article{Kalman1960,
  author    = {Kalman, R. E.},
  title     = {A New Approach to Linear Filtering and Prediction Problems},
  journal   = {Journal of Basic Engineering},
  volume    = {82},
  number    = {1},
  pages     = {35--45},
  year      = {1960},
  doi       = {10.1115/1.3662552},
  note      = {Transactions of the ASME, Series D}
}

@article{EM2,
  author    = {Dempster, A. P. and Laird, N. M. and Rubin, D. B.},
  title     = {Maximum Likelihood from Incomplete Data via the EM Algorithm},
  journal   = {Journal of the Royal Statistical Society: Series B (Methodological)},
  year      = {1977},
  volume    = {39},
  number    = {1},
  pages     = {1--38}
}

@article{OF1,
  author    = {Horn, B. K. P. and Schunck, B. G.},
  title     = {{Determining Optical Flow}},
  journal   = {Artificial Intelligence},
  year      = {1981},
  volume    = {17},
  number    = {1--3},
  pages     = {185--203},
  doi       = {10.1016/0004-3702(81)90024-2}
}

@book{Covariance,
  author    = {Rasmussen, Carl Edward and Williams, Christopher K. I.},
  title     = {Gaussian Processes for Machine Learning},
  publisher = {The MIT Press},
  address   = {Cambridge, MA},
  year      = {2006},
  isbn      = {0-262-18253-X}
}

@article{RKFS,
  author = {Farrell and Ioannou},
  title = {State estimation using a reduced order Kalman filter},
  year = {2001},
  journal = {J. Atmos. Sci.},
  volume = {58},
  pages = {3666-3680},
  publisher = {J. Atmos. Sci.},
  language = {eng},
}

@misc{ARCVT,
  title = {Facilities, Equipment, and Other Resources Statement},
  author = {{Advanced Research Computing at Virginia Tech (ARC)}},
  year = {2025},
  url = {https://www.docs.arc.vt.edu/pi_info/fer.html},
  note = {Accessed: November 4, 2025}
}

@inproceedings{OF2,
  author    = {Lucas, B. D. and Kanade, T.},
  title     = {An Iterative Image Registration Technique with an Application to Stereo Vision},
  booktitle = {Proceedings of the 7th International Joint Conference on Artificial Intelligence (IJCAI '81)},
  pages     = {674--679},
  year      = {1981}
}

@inproceedings{OF3,
  author    = {Brox, T. and Bregler, C. and Malik, J.},
  title     = {Large Displacement Optical Flow},
  booktitle = {Proceedings of the IEEE Conference on Computer Vision and Pattern Recognition (CVPR)},
  pages     = {41--48},
  year      = {2009},
  doi       = {10.1109/CVPR.2009.5206767}
}

@article{dim_red1,
  author    = {Solonen, A. and Cui, T. and Hakkarainen, J. and Marzouk, Y.},
  title     = {On dimension reduction in Gaussian filters},
  journal   = {Inverse Problems},
  volume    = {32},
  number    = {4},
  pages     = {045003},
  year      = {2016},
  month     = mar,
  doi       = {10.1088/0266-5611/32/4/045003},
  publisher = {IOP Publishing}
}

@article{FBP2,
  author    = {Shepp, L. A. and Logan, B. F.},
  title     = {{The Fourier Reconstruction of a Head Section}},
  journal   = {IEEE Transactions on Nuclear Science},
  year      = {1974},
  volume    = {21},
  number    = {3},
  pages     = {21--43},
  doi       = {10.1109/TNS.1974.6499235}
}

@article{buccini2021linearized,
  title={{Linearized Krylov subspace Bregman iteration with nonnegativity constraint}},
  author={Buccini, Alessandro and Pasha, Mirjeta and Reichel, Lothar},
  journal={Numerical Algorithms},
  volume={87},
  number={3},
  pages={1177--1200},
  year={2021},
  publisher={Springer}
}

@article{lan2023spatiotemporal,
  title={{Spatiotemporal Besov priors for Bayesian inverse problems}},
  author={Lan, Shiwei and Pasha, Mirjeta and Li, Shuyi and Shen, Weining},
  journal={Journal of the American Statistical Association},
  pages={1--15},
  year={2025},
  publisher={Taylor \& Francis}
}

@article{hansen1994regularization,
  title={Regularization tools: A Matlab package for analysis and solution of discrete ill-posed problems},
  author={Hansen, Per Christian},
  journal={Numerical algorithms},
  volume={6},
  number={1},
  pages={1--35},
  year={1994},
  publisher={Springer}
}

@article{chung2018efficient,
  title={Efficient generalized {G}olub--{K}ahan based methods for dynamic inverse problems},
  author={Chung, Julianne and Saibaba, Arvind K and Brown, Matthew and Westman, Erik},
  journal={Inverse Problems},
  volume={34},
  number={2},
  pages={024005},
  year={2018},
  publisher={IOP Publishing}
}

@article{huang2017majorization,
  title={Majorization--minimization generalized {K}rylov subspace methods for $\ell_p-\ell_q$ optimization applied to image restoration},
  author={Huang, G and Lanza, A and Morigi, S and Reichel, L and Sgallari, F},
  journal={BIT Numerical Mathematics},
  volume={57},
  number={2},
  pages={351--378},
  year={2017},
  publisher={Springer}
}

@article{okunola2025efficient,
  title={Efficient Dynamic Image Reconstruction with motion estimation},
  author={Okunola, Toluwani and Pasha, Mirjeta and Kilmer, Misha and Freitag, Melina},
  journal={arXiv preprint arXiv:2501.12497},
  year={2025}
}

@book{anderson2005optimal,
  title={Optimal filtering},
  author={Anderson, Brian DO and Moore, John B},
  year={2005},
  publisher={Courier Corporation}
}

@article{hakkarainen2019undersampled,
  title={Undersampled dynamic {X}-ray tomography with dimension reduction Kalman filter},
  author={Hakkarainen, Janne and Purisha, Zenith and Solonen, Antti and Siltanen, Samuli},
  journal={IEEE Transactions on Computational Imaging},
  volume={5},
  number={3},
  pages={492--501},
  year={2019},
  publisher={IEEE}
}

@article{pasha2024trips,
  title={{TRIPs-Py}: Techniques for regularization of inverse problems in {P}ython},
  author={Pasha, Mirjeta and Gazzola, Silvia and Sanderford, Connor and Ugwu, Ugochukwu O},
  journal={Numerical Algorithms},
  pages={1--38},
  year={2024},
  publisher={Springer}
}

@article{gazzola2019flexible,
  title={Flexible {GMRES} for total variation regularization},
  author={Gazzola, Silvia and Sabat{\'e} Landman, Malena},
  journal={BIT Numerical Mathematics},
  volume={59},
  pages={721--746},
  year={2019},
  publisher={Springer}
}

@article{lanza2015generalized,
  title={A generalized {K}rylov subspace method for $\ell_p$-$\ell_q$ minimization},
  author={Lanza, Alessandro and Morigi, Serena and Reichel, Lothar and Sgallari, Fiorella},
  journal={SIAM Journal on Scientific Computing},
  volume={37},
  number={5},
  pages={S30--S50},
  year={2015},
  publisher={SIAM}
}

@article{pasha2023computational,
  title={A computational framework for edge-preserving regularization in dynamic inverse problems},
  author={Pasha, Mirjeta and Saibaba, Arvind K and Gazzola, Silvia and Espa{\~n}ol, Malena I and de Sturler, Eric},
  journal={Electronic Transactions on Numerical Analysis},
  volume={58},
  pages={486--516},
  year={2023},
  publisher={Kent State University}
}

@article{stuart2010inverse,
  title={Inverse problems: A {B}ayesian perspective},
  author={Stuart, Andrew M},
  journal={Acta Numerica},
  volume={19},
  pages={451--559},
  year={2010},
  publisher={Cambridge University Press}
}

@book{calvetti2023bayesian,
  title={Bayesian Scientific Computing},
  author={Calvetti, Daniela and Somersalo, Erkki},
  volume={215},
  year={2023},
  publisher={Springer Nature}
}

\end{document}